\documentclass[12pt]{amsart}

\usepackage[T1]{fontenc}
\usepackage[cp1250]{inputenc}

\usepackage[arrow, matrix, curve]{xy}
\usepackage{amsmath,amsthm,amssymb}
\usepackage{amscd,graphics}

\usepackage{graphicx}

\newenvironment{Proof}{\textbf{Proof.}}{$\qquad \blacksquare$\par}
\newenvironment{Proof of}[1]{\textbf{Proof #1.}}{$\qquad \blacksquare$\par}

\newcommand{\mantysap}{\textbf{\{}} \newcommand{\mantysak}{ \textbf{\}}}
\DeclareMathOperator{\hull}{hull}
\renewcommand{\H}{\mathbb H}

\newcommand{\B}{\mathcal B}

\newcommand{\M}{\widetilde M}

\newcommand{\x}{\widetilde x}

\newcommand{\y}{\widetilde y}

\newcommand{\U}{\mathbb U}

\newcommand{\TPsi}{\widetilde \Psi}
\newcommand{\TPhi}{\widetilde \Phi}
\newcommand{\TDelta}{\widetilde \Delta}
\newcommand{\tdelta}{\widetilde \delta}

\newcommand{\tomega}{\widetilde \omega}
\newcommand{\tal}{\widetilde \alpha}

\newcommand{\al}{\alpha}

\newcommand{\A}{\mathcal A}
\newcommand{\TA}{\widetilde{\mathcal A}}
\renewcommand{\AA}{\mathbb A}

\newcommand{\C}{\mathbb C}
\newcommand{\R}{\mathbb R}
\newcommand{\Q}{\mathbb Q}
\newcommand{\Z}{\mathbb Z}
\newcommand{\N}{\mathbb N}

\newtheorem{thm}{Theorem}[section]
\newtheorem{lem}[thm]{Lemma}
\newtheorem{prop}[thm]{Proposition}

\theoremstyle{definition}

\newtheorem{defn}[thm]{Definition}
\newtheorem{ex}[thm]{Example}
\newtheorem{rem}[thm]{Remark}

\begin{document}
 \thispagestyle{empty}
   \title[$C^*$-algebras associated with  extensions of logistic maps]{$C^*$-algebras associated with reversible extensions of logistic maps}

\author{B. K.  Kwa\'sniewski}
%\thanks{This work was in part supported by Polish Ministry of Science and High Education grant number N N201 382634}

\address{Institute of Mathematics,  University  of Bialystok\\
ul. Akademicka 2,  PL-15-267  Bialystok,   Poland}
\email{bartoszk@math.uwb.edu.pl}
\urladdr{http://math.uwb.edu.pl/~zaf/kwasniewski}


\thanks{The author expresses his gratitude to A. V. Lebedev for his active interest in preparation of this paper. This work was in part supported by Polish Ministry of Science and High Education grant number N N201 382634}

\begin{abstract} The construction of  reversible extensions of  dynamical systems presented in \cite{maxid} is  enhanced, so that it applies to arbitrary mappings (not necessarily with open range). It is based on  calculating  the maximal ideal space of $C^*$-algebras that extends  endomorphisms  to partial automorphisms
 via partial isometric representations, and  involves a newfound  set of "parameters" (the role of parameters 
 play chosen sets or ideals). %Additionally, it is characterised as a universal object.
 \par
As model examples, we give  a thorough description of reversible extensions of logistic maps, and a classification of systems associated with compression of unitaries generating homeomorphisms of the circle.
\end{abstract}
\maketitle

\setcounter{tocdepth}{1}
 \tableofcontents

\section*{Introduction.}
A general $C^*$-method of construction of reversible extensions of irreversible dynamical systems was developed in  \cite{maxid}, where in particular the complete description of the maximal ideal spaces of the arising $C^*$-algebras was given. These algebras naturally spring out in the spectral analysis of weighted shift operators and transfer operators (see, in particular,  \cite{kwa-phd}, \cite{Anton_Lebed}) and as is shown in \cite{maxid} their maximal ideal spaces are tightly related  to inverse (projective) limits of dynamical systems.
\par
As a $C^*$-algebraic basis of their construction the authors of \cite{maxid} explored the leading concept of    \cite{Leb-Odz}
-- an  algebra whose elements play the role of \emph{Fourier
coefficients} in  partial isometric extensions of $C^*$-algebras. Namely, there was  studied the 
\mbox{$C^*$-algebra} $C^*(\A, U)$ generated by  a  $^*$-algebra $\A\subset L(H)$, $1\in \A$,  and  an operator $U\in L(H)$ under the assumption that $\A$  has    the following %conditions hold%
three properties:
\begin{align}
	\A \ni a &\to \delta(a):  =UaU^*\in \A,\label{covariant condition1}	 \\
		\A \ni a &\to \delta_*(a): =U^*aU\in \A,\label{covariant condition2}	 \\
	Ua   & =\delta(a)U, \qquad a\in \A. \label{covariant condition3}	
\end{align}
In this event,  $\A$ is called \emph{coefficient algebra} of $C^*(\A, U)$, cf. \cite[Prop. 2.4]{Leb-Odz}.
One has to stress that such objects are the major structural elements of the most successful  crossed product constructions, like the ones  developed by J. Cuntz and W. Krieger \cite{cuntz,cuntz-krieger}, W. L. Paschke
\cite{Paschke},  G. J. Murphy
\cite{Murphy}, R. Exel \cite{exel1} and others.
In general,  conditions \eqref{covariant condition1}, \eqref{covariant condition2},  \eqref{covariant condition3} %, together with  $1\in \A$,
 imply that  $\delta(\cdot)  =U (\cdot) U^*$  is an \emph{endomorphism} of $\A$ %and if $\delta$  is multiplicative on $\A$, then
  (then  $U$ is necessarily a \emph{partial isometry}) and $\delta_*(\cdot)  =U (\cdot) U^*$ is a unique non-degenerate \emph{transfer operator}  for $\delta:\A\to\A$ (in the sense of \cite{exel2}), see \cite{Bakht-Leb}. %, \cite{kwa-trans}.
  The general crossed-product based on relations \eqref{covariant condition1}, \eqref{covariant condition2},  \eqref{covariant condition3}  is developed in \cite{Ant-Bakht-Leb}, \cite{kwa-leb3}. As is shown in \cite{kwa-trans}, it could be viewed as the crossed product by a Hilbert bimodule \cite{aee}, and hence it is one of the fundamental models of relative Cuntz-Pimsner algebras \cite{ms},  $C^*$-algebras associated with $C^*$-correspondences \cite{katsura}, \cite{katsura1},  and Doplicher-Roberts algebras  \cite{dr}, see \cite{kwa-dr}.
\par
In this article we develop the  $C^*$-formalism of \cite{maxid} and apply  it to  a series of classical dynamical systems, in order to get the description of maximal ideal spaces of $C^*$-algebras associated with their reversible extensions.
We recall that the starting point of \cite{maxid} was a commutative unital $C^*$-subalgebra $\A\subset L(H)$ and an endomorphism  $\delta:\A\to \A$  such that
\begin{equation}\label{conditions for maxid}
	\A \ni a \to \delta(a):  =UaU^*\in \A, \qquad U^*U\in \A,
\end{equation}
for a certain  $U\in L(H)$.
Then   \eqref{covariant condition1} and  \eqref{covariant condition3} hold, and the $C^*$-algebra
$$
\B=C^*\big(\bigcup_{n=0}^\infty U^{*n}\A U^{n}\big)
$$
generated by  $\bigcup_{n=0}^\infty U^{*n}\A U^{n}$ is the smallest (still commutative) coefficient \mbox{$C^*$-algebra} of $C^*(\A,U)$ such that $\A\subset \B$. The passage from $\A$ to $\B$ corresponds to passage from irreversible to reversible dynamics. Namely, endomorphisms $\delta:\A\to\A$ and $\delta:\B\to\B$ are given, via  Gelfand transform, by   partial dynamical systems $(M,\al)$ and $(\M,\tal)$ where  $(\M,\tal)$ may be viewed as a universal reversible extension of  $(M,\al)$ (we will make the latter  statement   precise in Theorem~\ref{universality}). The authors  of \cite{maxid} gave a complete description of  $(\M,\tal)$  in terms of $(M,\al)$ and noticed that   $(\M,\tal)$ contains, as a subsystem, the inverse limit of  $(M,\al)$. This indicates  that  the structure of the $C^*$-algebra $\B$ is related to  hyperbolic attractors,   \cite{Williams2},  \cite{Brin}, \cite{Devaney} (such as  solenoids or  horseshoes of Smale);
irreversible continua   \cite{Nadler} (the most known are Brouwer-Janiszewski-Knastera continuum, or Knaster's  pseudoarc);
 and systems associated with classical substitution tilings \cite{Anderson Putnam} (these include tilings of Penrose, Amman, Fibonnaci, Morse, etc.).
\par
The results of \cite{maxid}, however, have   one  drawback. The  only seemingly  technical  assumption  $U^*U\in \A$  implies that the image of $\al$ is necessarily open, which in turn   excludes  a great deal of  important examples. As the first step in the present paper we  eliminate this inconvenience.
The key hint on how to overcome the mentioned drawback is  given in \cite[Rem. 3.7]{maxid}. Namely, one has to pass from the \mbox{$C^*$-algebra} $\A$ to  the \mbox{$C^*$-algebra}
 $\A_+:=C^*(\A,U^*U)$
  generated by $\A$ and the projection $U^*U$ and apply the $C^*$-method of  the reversible extension construction to the dynamical system generated on $\A_+$. In the starting Sections~\ref{preliminaries section} and~\ref{sect-2} we provide the corresponding analysis and as a  result  obtain
a description of the extended system $(\M,\tal)$ under the conditions  \eqref{covariant condition1}, \eqref{covariant condition3}
which  are weaker than \eqref{conditions for maxid}. As we show in Theorem \ref{takie zgrabne takie uadne thm}, these axioms embrace all endomorphisms of $\A$, and  thereby all partial dynamical systems  $(M,\al)$. The principal novelty here is uncovering of the fact that $(\M,\tal)$ depends not only on $(M,\al)$ but also on a certain set of parameters $Y \subset X$ (or ideals in $\A$). This new observation  has a number of interesting consequences. For instance, we  get nontrivial results implementing our method  to (already) reversible dynamical systems, such as homeomorphisms of the circle. Nevertheless,  our primary example and one of our  main goals is  a depictive presentation of $C^*$-algebras associated with reversible extensions of the family of \emph{logistic maps }$\al_\lambda:[0,1]\to [0,1]$:
\begin{equation}\label{odwozorowanie logistyczne label}
\al_\lambda(x)=4\lambda x(1-x), \qquad \,\,\,\, 0<\lambda \leq 1,
\end{equation}
 In the process of portraying the  maximal ideal spaces  of arising $C^*$-algebras, apart from the developed formalism, we  take advantage of the results concerning inverse limits of logistic maps  \cite{Barge-Ingram}, \cite{Barge-Martin}. In particular, we discuss in detail  how the extended systems are influenced by such phenomena as \emph{bifurcations} or  \emph{chaos}. % (existence of  periodic orbit of period not being a power of $2$).
\par
The paper is organized as follows. In Section \ref{preliminaries section} we introduce  notation and generalize or adapt from \cite{maxid} the basic concepts  required for presentation of the results of the  article. Here we also  discuss a general $C^*$-method of extending partial dynamical systems. The main result of \cite{maxid}, description of  maximal ideal spaces of $C^*$-algebras corresponding to reversible extensions of $C^*$-dynamical systems,
 is refined in Section  \ref{natural reversible extensions}, where
 apart from giving a purely topological definition  we  characterize  such systems as  universal objects. Section \ref{logistic section} is devoted to presentation of reversible extensions of the logistic family, and finally, in Section \ref{homeomorphisms of a circle section} we classify the $C^*$-algebras  associated with homeomorphisms of the circle  via their rotation numbers.

\section{Preliminaries. Endomorphisms of commutative $C^*$-algebras, dynamical systems and their extensions}\label{preliminaries section}
Throughout the article we let $\A$ be a  unital \emph{commutative} $C^*$-algebra. By using the Gelfand transform  we  assume the identification $\A=C(M)$,  where $M=M(\A)$ is the \emph{maximal ideal space} (also called   \emph{spectrum}) of the algebra $\A$.

\subsection{Endomorphisms and partial dynamical systems}
It is well known (see, for example, \cite[Thm 2.2]{maxid}) that  every   endomorphism   $\delta:\A\to \A$  is of the form
$$
\delta(a)=\left\{ \begin{array}{ll} a(\al(x))& ,\ \ x\in
\Delta\\
0 & ,\ \ x\notin \Delta \end{array}\right. ,\qquad a\in \A=C(M).
$$
where $\al:\Delta\to M$ is a continuous mapping defined on a closed and open (briefly clopen) subset $\Delta\subset M$. Namely, treating points of $M$ as  functionals on $\A$ we have
$$
\Delta=\{x\in M: x(\delta(1))\neq 0\}\quad\textrm{and}\quad \al=\delta^*|_{\Delta}
$$
where $\delta^*$ is the dual operator to $\delta:\A\to \A$. Therefore we will refer to $\al:\Delta\to M$  as to  a \emph{mapping dual to endomorphism} $\delta$.\par

To start with we describe the necessary for our further presentation objects related to algebras and endomorphisms.

For every subset $I\subset \A$ the set
 $$
 \hull(I):=\{x\in M: x(I)=0\},
 $$
 is a closed subset  of $M$ and if $I$ is an ideal in $\A$, then  $I=C_{\hull(I)}(M)$ where $C_K(M)$ stands for the set of continuous functions on $M$  vanishing on $K\subset M$.
Plainly, in the above notation we have
$$
\hull(\ker\delta)= \al(\Delta).
$$
 Since $\Delta$ is closed and $\al$ is continuous it follows that $\al(\Delta)$ is closed as well, and evidently   $\al(\Delta)$ is open iff the characteristic function of $\al(\Delta)$ belongs to $C(M)$ (in which case it is a unit in $\ker\delta=C_{\al(\Delta)}(M)$). Accordingly, we arrive at
\begin{prop}\label{openness of the image}
 The image  of the mapping $\al$ dual to an endomorphism $\delta:\A\to \A$ is open if and only if the kernel of $\delta$ is unital.
\end{prop}
The  \emph{annihilator} of an ideal $I$ in $\A$ is the set
    $$
   I^\bot:=\{a\in \A: a I=\{0\}\}.
   $$
Clearly, $I^\bot$  is an ideal and it
    could be equivalently defined as the largest ideal in $\A$ such that $I\cap I^{\bot}=\{0\}$.
In particular, if the kernel of $\delta:\A\to \A$ is unital, then  $\A$ admits  decomposition into the following direct sum of ideals
\begin{equation}\label{decomposition 1}
   \A=\ker\delta\oplus (\ker\delta)^\bot,
\end{equation}
  and $\delta$ yields an isomorphism between the ideal $(\ker\delta)^\bot$ and the subalgebra \mbox{$\delta(\A)\subset\A$.}
   \begin{defn}\label{czesciowy automorfizm defn}
An endomorphism $\delta:\A\to \A$ with unital kernel $\ker\delta$ and  the image $\delta(\A)$ being an ideal in $\A$ will be called a  \emph{partial automorphism}  of  $\A$.
   \end{defn}
Since $\delta(\A)$ is an ideal in $\A$ iff $\delta(\A)=\delta(1)\A$, it follows that   $\delta$ is a partial automorphism of  $\A$ iff   the algebra $\A$ admits two decompositions into direct sum of ideals
   $$
   \A=\ker\delta\oplus (\ker\delta)^\bot, \qquad\qquad  \A=\delta(\A)\oplus \delta(\A)^\bot.
   $$
If this is the case, we denote by $\delta_*: \delta(\A) \to(\ker\delta)^\bot$ the  inverse to the isomorphism $\delta:(\ker\delta)^\bot\to \delta(\A)$ and prolong $\delta_*$ onto $\A$ by putting $\delta_*|_{\delta(\A)^\bot}\equiv 0$.  Clearly, $\delta_*:\A\to \A$ is a  partial automorphism and
\begin{equation}\label{generalized inverses relations}
\delta\circ \delta_* \circ \delta= \delta, \qquad \delta_* \circ \delta \circ \delta_*= \delta_*.
\end{equation}
Hence $\delta_*$ is the so-called \emph{generalized inverse} to $\delta$.

The next proposition shows a deeper relation between the objects introduced above.
\begin{prop}\label{reversibility of endomorphisms}
Endomorphism $\delta:\A\to \A$  is a partial automorphism  if and only if its  dual map  $\al:\Delta\to \al(\Delta)$ is a homeomorphism and $\al(\Delta)$ is clopen in $M$.
\par
Moreover, if $\delta$  is a partial automorphism, then  there is a  unique partial automorphism  $\delta_*:\A \to \A$ which is  a generalized inverse for $\delta$, and it is given by
$$
\delta_*(a)=\left\{ \begin{array}{ll} a(\al^{-1}(x))& ,\ \ x\in
\al(\Delta)\\
0 & ,\ \ x\notin \al(\Delta) \end{array}\right. ,\qquad a\in \A=C(M).
$$
\end{prop}
\begin{Proof}
If $\delta:\A\to \A$  is a partial automorphism, then $\al(\Delta)$ is clopen (by Proposition \ref{openness of the image}) and  $\al:\Delta\to \al(\Delta)$ is a homeomorphism since it is a mapping dual to the isomorphism $\delta:(\ker\delta)^\bot\to \delta(\A)$. The converse implication is straightforward.  Suppose now that $\delta$ and $\delta_*$ are partial automorphisms   satisfying  \eqref{generalized inverses relations}. Then  $\delta= \delta\circ \delta_* \circ \delta$ implies that $\ker\delta_*\cap  \delta(\A)=\{0\}$, equivalently  $\delta(\A)\subset (\ker\delta^*)^\bot$, and    $\delta_*= \delta_*\circ\delta\circ  \delta_*$  implies  $\delta_*(\delta(\A))=\delta_*(\A)$. But since $\delta_*$ is a partial automorphism the latter relation gives  $(\ker\delta^*)^\bot\subset \delta(\A)$ and therefore $
(\ker\delta_*)^\bot=\delta(\A)$. By symmetry we also have   $(\ker\delta)^\bot=\delta_*(\A)
$, and thus it follows that $\delta_*: \delta(\A) \to(\ker\delta)^\bot$ coincides with the  inverse for $\delta:(\ker\delta)^\bot\to \delta(\A)$. In particular, $\delta_*$ is uniquely determined by $\delta$ and its dual mapping is $\al^{-1}$.
\end{Proof}
The above consideration makes it natural to adopt the following definitions, cf. \cite[Def. 2.4, 2.6]{maxid}.
\begin{defn}\label{definition of partial systems}
By a \emph{(partial) dynamical system} we will mean a triple $(M,\Delta,\al)$, where $M$ is a compact Hausdorff space, $\Delta$ a clopen subset of $M$, and $\al:\Delta\to M$ a continuous map. Unless a misunderstanding can arise,  we will simply write $(M,\al)$.
\end{defn}

\begin{defn}
We will say that  a partial dynamical system $(M,\Delta,\al)$ is \emph{reversible} if $\al(\Delta)$ is an open subset of $M$ and the map $\al:\Delta\to \al(\Delta)$ is a homeomorphism (so that the triple $(M,\al(\Delta),\al^{-1})$ is also a partial dynamical system).
\end{defn}

\subsection{Extensions of partial dynamical systems and endomorphisms}
The main concept of the $C^*$-method of construction of reversible extensions of irreversible dynamical systems is given in  \cite{maxid}. The present and the next subsections \ref{1.3}, \ref{1.4} and Section \ref{sect-2} are devoted to  description of the main structural blocks of this construction. Moreover we also give a refinement of the construction
that is necessary for complete  analysis of the partial dynamical systems under investigation in the paper, and establish universality and minimality of natural reversible extensions (Theorem \ref{universality}).
\begin{defn}\label{extension of system definition}
Let  $(M_\al,\Delta_\al,\al)$ and  $(M_\beta,\Delta_\beta,\beta)$ be partial dynamical systems. We say that a surjective continuous map $\Psi:M_\beta\to M_\alpha$ is a
\emph{semiconjugacy} (or a \emph{factor map}) if and only if the following conditions hold
\begin{equation}\label{semicon1}
\Psi^{-1}(\Delta_\al)=\Delta_\beta,
\end{equation}
\begin{equation}\label{semicon2}
\al(\Psi(x))= \Psi(\beta(x)),\qquad \quad x\in \Delta_\beta.
\end{equation}
If there is a semiconjugacy from $(M_\beta,\beta)$ to $(M_\al,\al)$ we say that $(M_\beta,\beta)$ is an \emph{extension} of $(M_\al,\al)$ and
 $(M_\al,\al)$ is a \emph{factor}  of $(M_\beta,\beta)$.  If  additionally the system $(M_\beta,\beta)$ is reversible we call it a \emph{reversible extension} of $(M_\al,\al)$. If the factor map is one-to-one, then its inverse is also a factor map,   and we say that $(M_\beta,\beta)$ and $(M_\al,\al)$ are \emph{conjugated}, or \emph{equivalent}.
  \end{defn}

The next proposition clarifies  the role of the objects introduced in the above definition and shows that   the  class of partial dynamical systems  with factor maps  as morphisms form a category dual to the category of endomorphisms of unital $C^*$-algebras where the role of morphisms is played by  unital monomorphisms that conjugate endomorphisms.
\begin{prop}\label{duality of categories}
Let $\delta:\A\to \A$ and  $\gamma:\B\to \B$ be endomorphisms of unital commutative $C^*$-algebras $\A$, $\B$, and let $(M(\A),\al)$ and $(M(\B),\beta)$ be the corresponding dual partial dynamical systems. Endomorphism $\gamma$ extends $\delta$ in the sense that there exits a unital monomorphism  $T:\A\to \B$ such that
\begin{equation}\label{intertwine}
T\circ \delta = \gamma\circ T,
\end{equation}
if and only if the partial dynamical system $(M(\B), \beta )$ is an extension of the system $(M(\A),\al)$. Furthermore, the factor map $\Psi:M(\B) \to  M(\A)$ is  the  dual map to the unital monomorphism $T$ satisfying \eqref{intertwine}: $
T(a) =a \circ \Psi$, $a\in \A=C(M(\A))$.
\end{prop}
\begin{Proof}
A mapping dual  to a unital monomorphism $T:\A\to \B$ maps $M(\B)$ onto $M(\A)$ and hence $T$ is an operator of composition with a surjection $\Psi:M(\B) \to  M(\A)$ where $\Psi=T^*|_{M(\B)}$. For $a$ we have
$$
(\gamma \circ T)(a)(x)=\begin{cases}
a(\Psi(\beta(x)), & x \in \Delta_\beta,\\
0, & x \notin \Delta_\beta,
\end{cases}\quad (T \circ \delta)(a)(x)=\begin{cases}
a(\al(\Psi(x)), & \Psi(x) \in \Delta_\al,\\
0, & \Psi(x) \notin \Delta_\al.
\end{cases}
$$
Using these formulas one checks that \eqref{intertwine} holds iff $\Psi$ satisfies \eqref{semicon1} and  \eqref{semicon2}.
\end{Proof}

  \begin{rem}
  \label{1}
  We have to stress that the extension described in  Definition \ref{extension of system definition} is  a slightly
   different (weaker) notion than the corresponding one described in \cite[Def. 2.7]{maxid}. Namely,   note that \eqref{semicon1} is equivalent to  two relations
  $$
  \Psi(\Delta_\beta)=\Delta_\al \ \ \text{and}\ \  \Psi(M_\beta\setminus \Delta_\beta)=M_\al\setminus \Delta_\al,
  $$
  and  \eqref{semicon2}  implies that
\begin{equation}\label{semicon3}
  \Psi(M_\beta\setminus \beta(\Delta_\beta))\supset M_\al \setminus \al (\Delta_\al).
\end{equation}
 However, unlike  \cite{maxid}, we allow $\al (\Delta_\al)$ not to be open (this is vital for the applications considered in this article). Moreover, as the  further part of the article show the principle situation of interest is   the case   when  $\beta(\Delta_\beta)$ is open, thus we can not require to have  equality in  \eqref{semicon3} (as in \cite[Def. 2.7]{maxid}).

  Continuing the above remark  one may consider  the  following consequence  of \eqref{semicon3}
\begin{equation}\label{semicon4}
 \Psi(\overline{ M_\beta\setminus \beta(\Delta_\beta)})\supset \overline{M_\al \setminus \al (\Delta_\al)}.
\end{equation}
 In particular, as it will be shown in the article, by refining (fixing) the left hand part of this inclusion one can obtain various  extensions naturally arising in the analysis of dynamical systems.
\end{rem}

\begin{defn}\label{set Y distinction}
Let $(M_\al,\Delta_\al,\al)$ and $(M_\beta,\Delta_\beta,\beta)$  be  partial dynamical systems and let $Y\subset M_\al$. If  there is a semiconjugacy $\Psi:M_\beta\to M_\alpha$
such that
\begin{equation}\label{semicon5}
 \Psi(\overline{M_\beta\setminus \beta(\Delta_\beta)})=Y
\end{equation}
 we  say that   $(M_\beta,\Delta_\beta,\beta)$  is an \emph{extension of $(M_\al,\Delta_\al,\al)$ associated with }$Y$.
\end{defn}
\begin{rem}
Plainly, the set $Y$ satisfying \eqref{semicon5} is necessarily closed and contains $M_\al \setminus \al (\Delta_\al)$. Conversely, for any such set $Y$ one easily constructs an  extension of $(M_\al,\Delta_\al,\al)$ associated with $Y$, see for instance Fig. \ref{uzwarcenie rysunek}.
\end{rem}

%Interpretation of Definition \ref{set Y distinction} in terms of endomorphisms  is related with annihilators of the kernels. More precisely,

Let us describe now the relation between extensions of dynamical systems associated with $Y$ and extensions of endomorphisms associated with ideals.

Suppose that $\delta:\A\to \A$ and  $\gamma:\B\to \B$ are endomorphisms conjugated by $T$, as in Proposition \ref{duality of categories}. Then \eqref{intertwine} implies the relations
\begin{equation}\label{relations for endomorphism extension}
\ker\delta =T^{-1}(\ker\gamma), \qquad (\ker\delta)^\bot \supset T^{-1}((\ker\gamma)^\bot),
\end{equation}
where by  taking hulls the former yields  \eqref{semicon1} and the latter gives \eqref{semicon4}. In particular,
$$
\hull((\ker\delta)^\bot)=\overline{M_\al \setminus \al (\Delta_\al)}\ \ \text{and} \  \ \hull \big(T^{-1}((\ker\gamma)^\bot)\big)=\Psi(\overline{ M_\beta\setminus \beta(\Delta_\beta)}),
$$
where $\Delta_\al$ and $\Delta_\beta$ are the domains of $\al$ and $\beta$ respectively, and $\Psi$ is the dual mapping to $T$.
Accordingly, we arrive at %Thus we get %the following improvement of Proposition \ref{duality of categories}.%In particular,  we get
\begin{thm}\label{duality of categories2}
 Let   $J$ be an ideal  in $\A$ and set  $Y=\hull(J)$. Under the notation of Proposition \ref{duality of categories}, the system $(M(\B), \beta )$ is an extension of the system $(M(\A),\al)$ associated with $Y$ if and only if  there exits a unital monomorphism  $T:\A\to \B$ such that
\begin{equation}\label{intertwine2}
T\circ \delta = \gamma\circ T, \quad \textrm{ and } \quad T^{-1}( (\ker\gamma)^\bot)=J.
\end{equation}
If this is the case, then  $J\subset (\ker\delta)^\bot$ (equivalently $Y\supset M(\A) \setminus \al (\Delta_\al)$),
so that endomorphism $\gamma$ extending $\delta$ "shrinks the annihilator of its kernel" to $J$.\end{thm}
\begin{rem}
Assuming the identification $\A \subset \B \ (\A \cong T(\A)\subset \B)$ \ we have $J=(\ker\gamma)^{\bot}\cap \A$ and relations \eqref{relations for endomorphism extension} tell that $\gamma$ may enlarge the kernel of $\delta$ but only outside $\A$, that is $(\ker\gamma)\cap \A=\ker\delta$. In other words, the only reason why  $J$ may be strictly smaller than $(\ker\delta)^\bot$ is that  $\ker\gamma \setminus \A $  is non-empty.  Thus it seems natural to look for extensions where  $J=(\ker\delta)^\bot$  (equivalently $Y=\overline{M(\A) \setminus \al (\Delta_\al)}$), cf. \cite[Def. 2.7]{maxid} and Remark \ref{1}.
\end{rem}
\subsection{$C^*$-dynamical systems and partial dynamical systems}
\label{1.3} Partial dynamical systems are naturally associated with $C^*$-dynamical systems and representations of the latter ones  in turn are defined by actions of partial isometries in  Hilbert spaces. Discussion of this
 relationship is the theme of the present subsection.
%The basic idea  of \cite{maxid} and the present paper is that extensions of dynamical systems arise naturally  when they are realized via operators  in a Hilbert space.
\begin{defn}\label{dynamical system definition}
By a \emph{(concrete) $C^*$-dynamical system} we mean a pair $(\A,U)$, where $\A$ is a commutative $C^*$-subalgebra of the algebra $L(H)$ of  bounded operators in a Hilbert space $H$, $1\in \A$, and $U\in L(H)$ is a partial isometry such that
\begin{equation}\label{conditions for system}
U\A U^*\subset \A, \qquad U^*U\in \A'
\end{equation}
where $\A'$ is the commutator of $\A$ in $L(H)$.

Plainly, relations  \eqref{conditions for system} imply that
\begin{equation}\label{definition of delta}
\delta(a):=UaU^*, \qquad a\in \A,
\end{equation}
is an endomorphism of $\A$, and we will say that $\delta:\A\to \A$ is an \emph{endomorphism generated by }$U$.
Similarly, if $(M,\Delta,\al)$ is the partial dynamical system dual  to $\delta$ we will  say that $(M,\Delta,\al)$  is a \emph{partial dynamical system  generated by} $U$. %on the maximal ideal space of $\A$.
\end{defn}
\begin{rem} The above definition extends \cite[Def. 2.11]{maxid}, see relations~\eqref{conditions for maxid}. In particular,
relations \eqref{conditions for system} and assumption that $U$ is a partial isometry are equivalent to \eqref{covariant condition1} and \eqref{covariant condition3}, cf.  \cite[Prop. 2.2]{Leb-Odz}.
\end{rem}
\begin{rem}\label{remark on representativity}
We will show  (in Theorem \ref{takie zgrabne takie uadne thm}) that for an arbitrary endomorphism \mbox{$\delta:\A\to \A$}  the $C^*$-algebra $\A$ may be identified with an algebra of operators acting in  a certain Hilbert space $H$ in such a way that $\delta$ is generated by an operator $U\in L(H)$.  Hence every  \emph{abstract $C^*$-dynamical system} $(\A,\delta)$ can be represented by a \emph{concrete one}.
\end{rem}
As it is indicated by Theorem \ref{duality of categories2},
in order to investigate extensions of an endomorphism  generated by $U\in L(H)$ we need to identify  the annihilator of its kernel  in terms of $(\A,U)$. To this end, let us consider  the  sets
$$
(1- U^*U)\A \cap \A=\{a\in \A: U^*Ua=0\},\qquad U^*U\A \cap \A=\{a\in \A: U^*Ua=a\},
$$
which  are (mutually orthogonal)  ideals in $\A$. The next proposition shows the role of these sets in the description of endomorphisms.
 \begin{prop}\label{introductory proposition}
Let $(\A,U)$ be a $C^*$-dynamical system, $\delta:\A\to \A$  an endomorphism generated by $U$, $(M,\Delta,\al)$  a partial dynamical system dual to $\delta:\A\to \A$, and $Y=\hull(U^*U\A \cap \A)$. Then
$$
\ker\delta=(1-U^*U)\A \cap \A\quad\textrm{and}\quad  (\ker\delta)^\bot \supset U^*U\A \cap \A,
$$
that is $Y \supset M\setminus \al(\Delta)$. If additionally $U^*U\in \A$, then
 $$\ker\delta=(1-U^*U)\A\quad \textrm{ and }\quad  (\ker\delta)^\bot =U^*U\A,
 $$
 that is $Y = M\setminus \al(\Delta)$. In particular,
 \begin{itemize}
 \item[i)] $U^*U\in \A$ $\,\,\Longrightarrow\,\,$ $\ker\delta$ is unital  ($\al(\Delta)$ is clopen);
 \item[ii)]  $U$ is an isometry $\Longrightarrow$ $\delta:\A\to \A$ is a monomorphism  ($\al(\Delta)=M$);
%  \item[iii)] $U^*$ is an isometry $\,\,\Longleftrightarrow\,\,$ $\delta:\A\to \A$ is unital, and  $\Delta=M$;
 \item[iii)]  $U$ is unitary $\,\,\Longrightarrow\,\,$ $\delta:\A\to \A$ is a unital monomorphism ($\Delta=M$ and $\al:M \to M$ is surjective).
  \end{itemize}
  \end{prop}
  \begin{Proof}
  To see that  $(1-U^*U)\A \cap \A$ coincides with $\ker\delta$ let $a\in \A$ and note that
  \begin{align*}
  U^*Ua =0\,\,\,& \Longrightarrow \,\,\, \delta(a)=U aU^* =U(U^*U a)U^*=0,\\
    U^*Ua \neq 0\,\,\,& \Longrightarrow  \,\,\, U^*\delta(a)U= U^*U a U^*U = U^*Ua \neq 0 \,\,\,\Longrightarrow \,\,\, \delta(a)\neq 0.
   \end{align*}
Now, since
$
\ker\delta\cap  \big(U^*U\A \cap \A\big)=\{0\}
$ we have $U^*U\A \cap \A\subset (\ker\delta)^\bot$, and if $U^*U\in \A$, the projection $(1-U^*U)$ is the unit for $\ker\delta=(1-U^*U)\A$, and consequently $U^*U$ is the unit for $(\ker\delta)^\bot$. Thus the remaining part of proposition is  straightforward, cf. \cite[Prop. 2.3]{maxid}.
  \end{Proof}
In the next Lemma \ref{lemma on carriers} and Theorem \ref{takie sobie stw} we give a  description of  $(\ker\delta)^\bot$ and  obtain restraint for  $U^*U$    in terms of the objects related to  central carriers of elements in  von Neumann algebras (for completeness of presentation we include in Lemma \ref{lemma on carriers} the known properties i) and ii) of carriers).
  \begin{defn}
By a \emph{carrier of a   $C^*$-subalgebra} $K\subset L(H)$ we mean the orthogonal  projection $Q\in L(H)$ onto the subspace $K H \subset H$.
  \end{defn}
  \begin{lem}\label{lemma on carriers}
Let  $Q\in L(H)$ be a carrier of an ideal $I$ in a (not necessarily commutative) $C^*$-algebra $\A\subset L(H)$. Then
\begin{itemize}
\item[i)] $Q=\textrm{s-}\lim_{\lambda_\in \Lambda} \mu_\lambda$ where $\{\mu_\lambda\}_{\lambda\in \Lambda}$ is an approximate  unit for $I$ (the  limit is taken in the strong operator topology).
\item[ii)] $Q\in \A'$.
\item[iii)] $ I\subset Q\A\cap \A$ and $I^{\bot}=(1-Q)\A\cap \A$.
\end{itemize}
  \end{lem}
 \begin{Proof}
By Vigier's Theorem  the limit  $Q:=\textrm{\emph{s}-}\lim_{\lambda_\in \Lambda} \mu_\lambda$ exists. For $a\in\A$ and $h\in H$ we have
$$
  Qah=\lim_{\lambda \in \Lambda} \mu_\lambda ah=\lim_{\lambda\in \Lambda} \mu_\lambda (\lim_{\lambda'\in \Lambda} a \mu_{\lambda'}  )h= QaQh,
  $$
  and similarly
  $$
  aQh=\lim_{\lambda\in \Lambda} \mu_\lambda ah=\lim_{\lambda\in \Lambda} \mu_\lambda (\lim_{\lambda'\in \Lambda} a \mu_{\lambda'}  )h= QaQh.
  $$
Using these relations  one deduces that $Q$ is a carrier of $I$ and  $Q\in \A'$. Clearly,  $I\subset QA\cap \A$ and $I^\bot \subset (1-Q)\A\cap \A$. To see that $I^\bot \supset (1-Q)\A\cap \A$ let $a\notin \A\setminus I^\bot$. Then there is $b\in I$ and $h\in H$ such that $ab h\neq 0$. Hence
  $
  Qab h=abh\neq 0
  $ and therefore  $Qa\neq 0$, which is equivalent to $(1-Q)a\neq a$.
  \end{Proof}
By item iii) in the above lemma  we see that  complements of  carriers are 'born' to deal  with annihilators of ideals.
  \begin{thm}\label{takie sobie stw}
Under the notation of Proposition \ref{introductory proposition}, let $P=1-Q\in L(H)$ be the complement of the carrier $Q$ of the ideal  $(1-U^*U)\A\cap \A=\ker\delta$ in $\A$. Then
$$
U^*U \leq P, \qquad  P \in \A'\qquad \textrm{ and }\qquad (\ker\delta)^\bot= P\A \cap \A.
$$
Moreover, if $(\ker\delta)^\bot= U^*U\A \cap \A$ (equivalently $Y=\overline{M\setminus \al(\Delta)}$), which automatically holds  when  $U^*U=P$,  then the implications in items i)-iii) in Proposition \ref{introductory proposition}  are  in fact equivalences.
\end{thm}
 \begin{Proof}
 By Lemma \ref{lemma on carriers},  $P \in \A'$  and $(\ker\delta)^\bot= P\A \cap \A$. By definition of $Q$,    $Q  \leq 1-U^*U$ and thus   $U^*U\leq P$. Now suppose that $(\ker\delta)^\bot= U^*U\A \cap \A$; by Proposition \ref{introductory proposition}, condition $U^*U\in \A$ gives even more, namely   $(\ker\delta)^\bot =U^*U\A$.
We  show  the converses  to implications  in items i)-iii) of  Proposition \ref{introductory proposition}.
 \par
 i) If $\ker\delta$ is unital, then $Q$ is the unit  for $\ker\delta$ and consequently $P$ is the unit  for $(\ker\delta)^\bot= U^*U\A \cap \A$. This implies that  $U^*U=P\in \A$.
\par
ii) If $\ker\delta=\{0\}$, then $\A=(\ker\delta)^\bot= U^*U\A \cap \A$, that is $\A=U^*U\A$ and therefore $U$ is an isometry.
\par
iii) It suffices to combine items i) and ii).
   \end{Proof}

The main theme of \cite{maxid} is a description of the $C^*$-method of construction of reversible $C^*$-dynamical systems, that is systems such that   not only $(\A,U)$ but also $(\A,U^*)$ is a $C^*$-dynamical system, and 
thus  both of the mappings
 $$
\delta(a):=UaU^*,\qquad \delta_*(a):=U^*aU,\qquad a\in \A
$$
  are endomorphisms of $\A$. We adopt  an equivalent version of  \cite[Def. 2.15]{maxid}.
\begin{defn}
By a \emph{reversible $C^*$-dynamical system} we mean a pair $(\A,U)$
 where $\A\subset L(H)$ is commutative, $1\in \A$, and $U\in L(H)$ is a partial isometry such that
\begin{equation}\label{conditions for  reversible system}
U\A U^*\subset \A, \qquad U^*\A U\subset \A,
\end{equation}
Clearly,  relations \eqref{conditions for  reversible system} are equivalent to the condition  that  both $(\A,U)$ and $(\A,U^*)$ are $C^*$-dynamical systems.
\end{defn}
      \begin{prop}\label{reversible C-dynamical systems}
If $(\A,U)$ is a reversible $C^*$-dynamical system, then  endomorphisms $\delta:\A\to \A$ and $\delta_*:\A\to \A$ generated by $U$ and $U^*$  are mutually generalized inverse partial automorphisms and
   $$
   (\ker\delta)^\bot =\delta_*(\A)=U^*U\A,\qquad \delta(\A)=(\ker\delta_*)^\bot=UU^*\A.
   $$ In particular, the partial dynamical system $(M(\A),\al)$ generated by $U$ is reversible.
   \end{prop}
\begin{Proof}
By Proposition  \ref{introductory proposition} and the symmetry between $\delta$ and $\delta_*$ it suffices to show that $\delta(\A)=UU^*\A$ which follows because
$
\delta(\A)=U\A U^*=UU^* U\A U^*\subset UU^*\A,
$
and
$
 UU^*\A=UU^*\A UU^*= \delta(\delta_*(\A))\subset \delta(\A).
$
\end{Proof}

   \subsection{$C^*$-method of extending partial dynamical systems}\label{1.4}
   Given a concrete $C^*$-dynamical system $(\A,U)$ we  have  at our disposal two mappings
   that are defined on the whole of the $C^*$-algebra $L(H)$:
\begin{equation}\label{delta defn}
\delta(a):=UaU^*,\qquad \delta_*(a):=U^*aU,\qquad a\in L(H).
\end{equation}
Moreover, $\delta$ restricted to $\A$ is an endomorphism, and $\delta_*$ restricted to $\A$ is an endomorphism iff   $(\A,U)$ is reversible. Obviously $\delta$ (and similarly $\delta_*$) may define endomorphisms on many different subalgebras of $L(H)$, and if additionally  such a subalgebra, say $\B$, contains $\A$ it yields a natural extension $\delta:\B\to \B$ of the initial endomorphism $\delta:\A\to \A$, cf. \cite[Def. 2.17]{maxid}. This in turn can be rewritten  in the language of partial dynamical systems.
   \begin{thm}\label{takie sobie stw1}
Let $(\A,U)$ and $(\B,U)$ be  $C^*$-dynamical systems  and let $(M(\A),\al)$ and  $(M(\B),\beta)$ be partial dynamical systems generated by $U$ on the maximal ideal spaces of $\A$ and $\B$, respectively.  If  $\A\subset \B$, then $(M(\B),\beta)$ is an extension of $(M(\A),\al)$ associated with the set
$$
Y=\hull((\ker\delta|_{\B})^\bot
\cap \A)
=\hull(P\B
\cap \A)
$$
where $(\ker\delta|_{\B})^\bot$ is the annihilator  of the kernel of $\delta:\B\to \B$, and $P$ is the complement of the carrier of  $(1-U^*U)\B\cap \B$. Moreover,
\begin{itemize}
\item[i)] $U^*U\in \B \,\,\, \Longrightarrow \,\,\, Y=\hull(U^*U\A\cap \A)$,
\item[ii)] $U^*U\in \A\,\,\, \Longrightarrow \,\,\, Y=M(\A)\setminus \al(\Delta)$.
\end{itemize}
\end{thm}
\begin{Proof}
Since $\delta:\B \to \B$ is an extension of $\delta:\A\to \A$ in the sense of Proposition~\ref{duality of categories}, where $T=id$, the first part of assertion follows from Theorem \ref{duality of categories2} and Proposition \ref{takie sobie stw}. If additionally $U^*U\in \B$, then in view of  Proposition \ref{introductory proposition}  we have $(\ker\delta|_{\B})^\bot
\cap \A= U^*U\B \cap \A=U^*U\A\cap \A$. Similarly, if  $U^*U\in \A$, then $(\ker\delta|_{\B})^\bot
\cap \A=U^*U\A$ and $\hull(U^*U\A)=M(\A)\setminus \al(\Delta)$, cf. Proposition \ref{openness of the image}.
\end{Proof}
Using the above method  one can always obtain a reversible extension  which  (in the context of coefficients algebras) was first noticed in \cite{Leb-Odz},  cf. \cite{maxid}.
   \begin{thm}\label{takie sobie stw2}
Suppose that, under the notation of Theorem  \ref{takie sobie stw1},  $\A\subset \B$ and $(\B,U)$ is reversible. Then $(M(\B),\beta)$ is a reversible extension of $(M(\A),\al)$ associated with the set
$$
Y=\hull(U^*U\A\cap \A).
$$
Moreover, for any $C^*$-dynamical system $(\A,U)$ there exists a minimal reversible $C^*$-dynamical system $(\B,U)$ such that $\A\subset \B$. Namely,
the $C^*$-algebra
\begin{equation}\label{extended algebra}
\B=C^*\big(\bigcup_{n=0}^\infty U^{*n}\A U^{n}\big)
\end{equation}
 generated by $\bigcup_{n=0}^\infty U^{*n}\A U^{n}$ is commutative, and it is the smallest $C^*$-algebra that contains $\A$ and satisfies \eqref{conditions for  reversible system}.
\end{thm}
\begin{Proof}
The first part of assertion follows from Proposition \ref{reversible C-dynamical systems} and Theorem~\ref{takie sobie stw1} (since  $U^*U=U^* 1 U \in \B$). Commutativity of the   algebra $\B$, given by \eqref{extended algebra} and the property that  both   $\delta:\B\to \B$ and $\delta_*:\B\to \B$ are endomorphisms was established in \cite[Prop. 4.1]{Leb-Odz}. The remaining part is straightforward.
\end{Proof}

\section{Natural reversible  extensions of dynamical systems}\label{natural reversible extensions}\label{sect-2}
In this section we give a complete purely topological description of    partial dynamical systems  corresponding to   minimal reversible extensions of  $C^*$-dynamical systems introduced in Theorem \ref{takie sobie stw2}.
Additionally,
we  characterize such systems as universal objects and  discuss their relation  to the notion of  inverse limit, which therefore we recall now.
\begin{defn}\label{definicja granicy odwrotnej}
If  $\al:M\to M$ is a continuous mapping of a topological space $M$, then the \emph{inverse limit}\label{granica odwrotna defn} of the inverse sequence $
M\stackrel{\al}{\longleftarrow} M
\stackrel{\al}{\longleftarrow}M  \stackrel{\al}{\longleftarrow}...\,
$
is the topological space of the form
$$
\underleftarrow{\,\,\lim\,} (M,\al):=\{(x_0,x_1,...)\in \prod_{n\in \N} M : \al(x_{n+1})=x_n, \, n \in \N\}$$
equipped with the product topology inherited from $\prod_{n\in \N} M$. Furthermore, on the space  $\underleftarrow{\,\,\lim\,} (M,\al)$ we have a naturally defined homeomorphism
$$
\tal(x_0,x_1,...)=(\al(x_0),x_0,x_1,...), \qquad (x_0,x_1,...)\in \underleftarrow{\,\,\lim\,} (M,\al),
$$
called a \emph{homeomorphism induced by the mapping} $\al:M\to M$.
\end{defn}

\subsection{$C^*$-dynamical approach}
Throughout this section we let  $(\A,U)$ be  a $C^*$-dynamical system,
 $\delta$ and $\delta_*$  the mappings given by  \eqref{delta defn}, and $\B$ the  $C^*$-algebra from Theorem \ref{takie sobie stw2}. We denote by  $(M,\Delta,\al)$ the partial dynamical system dual to $\delta:\A\to \A$ and  by  $(\M,\TDelta,\tal)$ the reversible partial dynamical system dual to $\delta:\B\to \B$:
 $$
 \A\cong C(M),\qquad \B=C^*\big(\bigcup_{n=0}^\infty U^{*n}\A U^{n}\big)\cong C(\M).
 $$
The dynamical system $(\M,\TDelta,\tal)$ plays the principal role in the paper. As was noted in \cite[Rem. 3.7]{maxid} to obtain description of $(\M,\TDelta,\tal)$ in terms of $(M,\Delta,\al)$  one has to apply  the main result of \cite{maxid} to
the $C^*$-algebra
 $$\A_+:=C^*(\A,U^*U)$$
  generated by $\A$ and the projection $U^*U$.  Therefore we need to analyze the algebra $\A_+$
   and the partial dynamical system dual to endomorphism $\delta:\A_+\to \A_+$. Hereafter we proceed to the discussion of these objects.
\subsubsection{Adjoining  the projection $U^*U$ to the algebra $\A$}\label{podrozdzial o zbiorzez Y}
\begin{figure}[htb]
  \begin{center}
\setlength{\unitlength}{1.mm}
\begin{picture}(125,41)(0,0)
\put(-8,-3){\includegraphics[angle=0, scale=0.32]{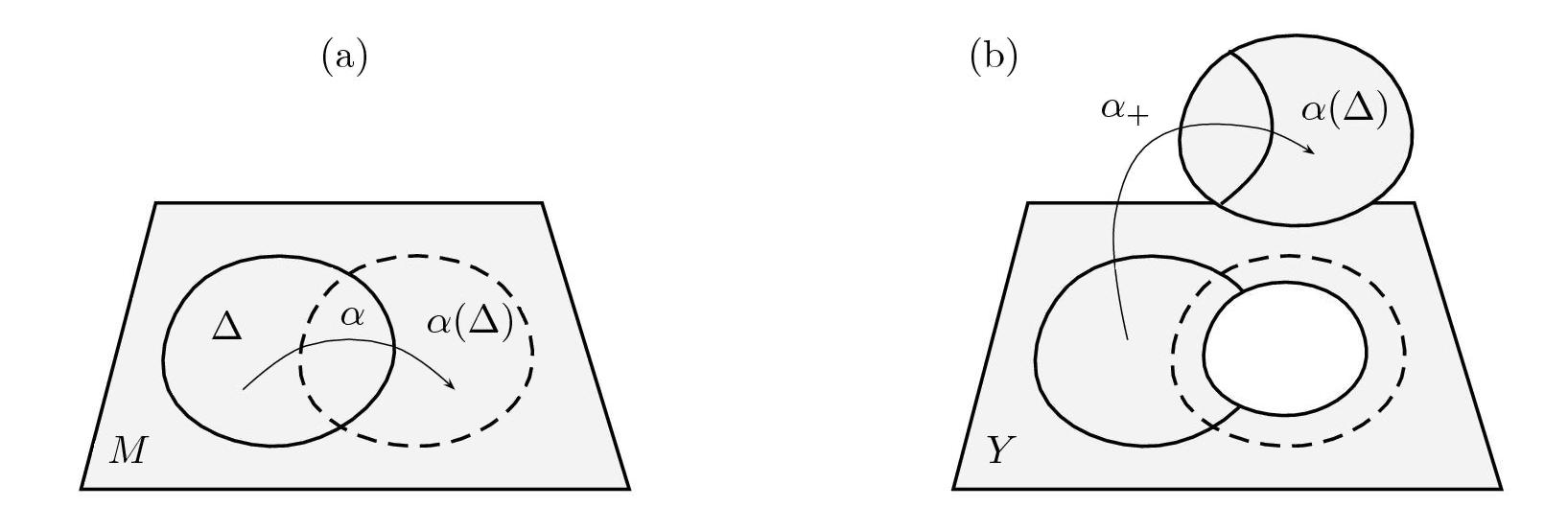}}
 \end{picture}
    \end{center}
   \caption{Partial dynamical system $(M,\al)$ generated by $U$ on $\A$ (a); partial dynamical system  $(M_+,\al_+)$ generated by $U$ on $\A_+$ (b). \label{uzwarcenie rysunek}}
 \end{figure}

Plainly, the $C^*$-algebra $\A_+$ is the direct sum of ideals
\begin{equation}\label{algebra A_+}
\A_+= U^*U\A \oplus (1-U^*U)\A
\end{equation}
where $(1-U^*U)\A=\ker (\delta|_{\A_+})$ is the kernel of $\delta:\A_+\to \A_+$, see Proposition \ref{introductory proposition}. As a simple consequence, cf. for instance \cite[Lem. 10.1.6]{kadison},  we get
 \begin{prop}\label{stw algebra A_+}
Algebra $\A_+$ is isomorphic the direct sum of quotient algebras
\begin{equation}\label{algebra A_+1}
\A_+\cong \A/\ker(\delta|_{\A})\, \oplus \,\A/\big(U^*U\A \cap \A\big)
\end{equation}
 where  $\ker(\delta|_{\A})=(1-U^*U)\A \cap \A$ is the kernel of $\delta:\A\to \A$ and
$U^*U\A \cap \A $
is an ideal  contained in the annihilator  $\ker(\delta|_{\A})^\bot$ of $\ker(\delta|_{\A})$. Under the isomorphism \eqref{algebra A_+1}   endomorphism $\delta:\A_+\to \A_+$ takes the form
$$a + \ker(\delta|_{\A})\, \oplus\, b + \big(U^*U\A \cap \A\big)\longmapsto \delta(a) +\ker(\delta|_{\A})\,\oplus \,\delta(a) +\big(U^*U\A \cap \A\big).$$
 \end{prop}
 \begin{rem} The  summand
$
(1-U^*U)\A\cong \A/\big(U^*U\A \cap \A\big)
$ is a unital $C^*$-algebra which coincides with the kernel of $\delta:\A_+\to \A_+$. Thus  one may interpret adjoining the projection $U^*U$ to the algebra $\A$ as a unitization of  the kernel of $\delta:\A\to \A$ (equivalently compactification of the complement of the image of $\al$).
 \end{rem}
Let $(M_+,\Delta_+,\al_+)$ denote the partial dynamical system generated by $U$ on the maximal ideal space of $\A_+$. In view of Theorem \ref{takie sobie stw1},  $(M_+,\al_+)$ is an extension of $(M,\al)$ associated with the set
 %\begin{equation}\label{omega defn equation}
 $$
 Y= \hull(U^*U\A \cap \A).
$$
% \end{equation}
By passing in Proposition \ref{stw algebra A_+}  to duals we get the complete description of the partial dynamical system $(M_+,\Delta_+,\al_+)$.
  \begin{prop}\label{uklad (M_+, al_+) prop}
Under the above notation  $Y$ is the closed set containing   $M\setminus \al(\Delta)$,
and the spectrum  $M_+$ of the algebra $\A_+$ realizes as the following  (topological) direct sum
 $$
 M_+=\al(\Delta) \sqcup Y.
  $$
The partial mapping   $\al_+$ dual to  $\delta:\A_+\to \A_+$  is defined on the set
 $$
\Delta_+=(\al(\Delta)\cap \Delta)\,\sqcup\, (Y\cap \Delta)
 $$
and  attains values in the first summand of  $M_+=\al(\Delta) \sqcup Y$ acting by the formula
 $$
 \al_+(x)=\al(x), \qquad x\in \Delta_+,
$$
see Fig. \ref{uzwarcenie rysunek}. In particular,  $U^*U\in \A$, that is  $\A=\A_+$, if and only if  $Y=M\setminus \al(\Delta)$, and then   $\al(\Delta)$ is clopen.
 \end{prop}
 \begin{Proof}
 Since the maximal ideal spaces of the quotient algebras $\A/ \ker(\delta|_{\A})$  and $\A/\big(U^*U\A \cap \A\big)$ may be identified with the sets $\hull(\ker(\delta|_{\A}))=\al(\Delta)$ and $\hull(U^*U\A \cap \A)=Y$, the first part of assertion  follows from Proposition \ref{stw algebra A_+}. For the second part  apply  Theorem \ref{takie sobie stw}.
\end{Proof}
Pictorially speaking, the space $M_+$   arise from $(M,\al)$ by
"cutting out" $\al(\Delta)$  and "replacing"  the set
$M\setminus \al(\Delta)$  with  the closed set $Y$ that contains $M\setminus \al(\Delta)$. The image of   $\al_+$  is the set
$\al(\Delta)$ that has been "cut out" and thus  may be identified with the image of $\al$. The domain     of $\al_+$
is enlarged with respect to the domain of $\al$  by the set
 $Y\cap (\Delta  \cap \al(\Delta))$.

Now, as we have obtained the description of the system  $(M_+,\Delta_+,\al_+)$ we pass to the main object of this subsection -- the system $(\M,\TDelta,\tal)$.

\subsubsection{Description of the  system $(\M,\tal)$ dual to  $\delta:\B\to \B$}\label{subsection opis spektralny}
 %Now we are approach the main problem of the section.% As a first application we will get  Theorem \ref{takie zgrabne takie uadne thm}.
Let \mbox{$\x\in \M$} be  a multiplicative linear functional on  $\B$.  We associate to it a sequence of functionals on $\A$ given by the formula
$$
x_n(a):=\x ( \delta_*^n( a)), \qquad a\in \A,\,\, n\in \N.
$$
Since  $\delta_*:\B\to \B$ is an endomorphism, the functionals   $x_n: \A \to \C$ are  multiplicative linear, and thus
$$
x_n \in M\quad \textrm{ or }\quad  x_n\equiv 0.
$$
Moreover, since  $\B=C^*(\bigcup_{n=0}^\infty \delta_*^n(\A))$ the sequence
$(x_0,x_1,...)$ determines  $\x$ uniquely. Consequently we have the injective mapping
\begin{equation}\label{isomor homeomor}
\x\longmapsto (x_0,x_1,...,x_n,...)
\end{equation}
that embeds the space   $\M$ into the Cartesian product  $\prod_{n\in \N} (M\cup \{0\})$ of the countable number of copies of the space $M\cup \{0\}$.

The next theorem is a refinement of the main result of \cite{maxid}  up to  objects under consideration.
  \begin{thm}[Description of the reversible system $(\M,\tal)$]\label{opis kosmiczakow}
The maximal ideal space of the  algebra $\B$ may be identified via the mapping \eqref{isomor homeomor} with the following topological space
$$
\M=\bigcup_{N=0}^{\infty}M_N\cup M_\infty
$$
where
$$
M_N=\{(x_0,x_1,...,x_N,0,...): x_n\in \Delta,\, \al(x_{n})=x_{n-1},\,\, n=1,...,N,\,x_N\in
Y\},
$$
$$
M_\infty=\{(x_0,x_1,...): x_n\in \Delta,\,
\al(x_{n})=x_{n-1},\, \,n\geqslant 1\},
$$
are equipped with the product topology inherited from $\prod_{n\in \N} (M\cup \{0\})$, where $\{0\}$ is a clopen singleton, and  $Y$ is a closed set containing $M\setminus \Delta$ (namely $Y=\hull(U^*U\A\cap \A)$). Furthermore the partial homeomorphisms  dual to   partial automorphisms $\delta, \delta_*:\B\to \B$ are defined respectively on the clopen sets
$$
\TDelta=\{(x_0,x_1,...)\in \M: x_0 \in \Delta\},\qquad \tal(\TDelta)=\{(x_0,x_1,...)\in \M: x_1 \neq 0\}
$$
and  act according to formulae
\begin{equation}\label{czesciowe homeomorfizmy rozszserzeone}
\tal(x_0,x_1,...)=(\al(x_0),x_0,x_1,...),\qquad \tal^{-1}(x_0,x_1,...)=(x_1,...).
\end{equation}
\end{thm}
\begin{Proof} It suffices to apply \cite[Thm. 3.5, 4.1]{maxid} to the $C^*$-dynamical system $(\A_+,U)$ and then use  Proposition \ref{uklad (M_+, al_+) prop}, cf. also \cite[Rem. 3.7]{maxid}.
\end{Proof}
\begin{rem}\label{mapping Phi remark}
The mapping $\Phi: \M\to M$ dual to the inclusion $\A\subset \B$ is given by the formula
\begin{equation}\label{rzut na zerowo wspolrzedno}
\Phi(x_0,x_1,...)=x_0.
\end{equation}
It is a factor map establishing that $(\M,\tal)$ is a reversible extension of $(M,\al)$ associated with $Y$. %, cf. Theorem \ref{takie sobie stw2}.
\end{rem}
The next result shows that in the situations when $U$ belongs to a series of commonly exploited classes of operators  the structure of  $(\M,\tal)$  becomes more transparent.
\begin{thm}\label{szczególne przypadki}
Under the  notation of Theorem \ref{opis kosmiczakow}, we have $M\setminus \al(\Delta)\subset Y$ and
\begin{itemize}
\item[i)] If $U^*U\in \A$, then  $Y=M\setminus \al(\Delta)$, and hence a sequence  $(x_0,x_1,...,x_N,0,...)$ is an element of  $M_N$ iff
$$
x_N\notin\al(\Delta)\quad  \textrm{ (i.e. }  x_N  \textrm{ does not have a preimage)}
$$
and  $x_n\in \Delta$, $\al(x_{n})=x_{n-1}$ for $n=1,...,N. $
\item[ii)]  If $U$ is an isometry, then  $\al:\Delta\to M$ is a surjection, and thus
$$
\M= M_\infty.
$$
\item[iii)]  If  $U$ is unitary, then  $\al:M\to M$ is a surjection,
$$
\M=\underleftarrow{\,\,\lim}(M,\al)
$$
  and $\tal:\M\to \M$ is a homeomorphism induced by  $\al:M\to M$.
\end{itemize}
\end{thm}
\begin{Proof}
In view of Theorem  \ref{opis kosmiczakow}, item  i) follows from Proposition \ref{uklad (M_+, al_+) prop}, whereas items ii),  iii) follow by Proposition \ref{introductory proposition}.
\end{Proof}

\subsection{Construction of operators generating arbitrary partial dynamical system $(M,\al)$} \label{Construction of operators generating}
Let $(M,\Delta,\al)$ be a partial dynamical system and let $Y\subset M$ be a closed set containing $M\setminus \al(\Delta)$. Theorem~\ref{opis kosmiczakow} leads to a natural question: does there exist a $C^*$-dynamical system $(\A,U)$  generating the partial dynamical system  $(M,\Delta,\al)$ and such that its reversible extension is associated with $Y$. In other words,  whether all the objects described in Theorem~\ref{opis kosmiczakow} are realizable. The answer is: yes.    Namely, following \cite[4.2]{maxid}, we  use the explicit description of the reversible extension $(\M,\tal)$ associated with $Y$   to give  a simple construction of operators generating $(M,\al)$.

Let
$$\label{partial isometries za pomoca ukladu rozszerzonego}
(Uf)(\x)=\begin{cases}
f(\tal(\x)),& \x \in \TDelta \\
0, & \x \notin \TDelta
\end{cases}
$$
act in the Hilbert space  $H=\ell^2(\M)$  (which may be treated as  $L^2_\mu(\M)$ where $\mu$ is the counting measure), and let  $\A\subset L(H)$ be the algebra consisting of operators of multiplication by function from $C(\M)$ dependent only on the zeroth coordinate:
$$\A=\{a\in C(\M): a(\x)=a(x_0), \textrm{ where }\x=(x_0,...)\in \M\}.  $$
Clearly,  $\A\cong C(M)$ and   $U$  generates on $M$ the partial mapping $\al$.  The projection  $U^*U$ is the operator of multiplication by a characteristic function of $\tal(\TDelta)$, and thus
$
U^*U\A \cap \A=\{ a\in \A: U^*U a=a\}\cong C_{Y}(M).
$
Accordingly, we get
\begin{thm}\label{takie zgrabne takie uadne thm}
Let $(M,\al)$ be an arbitrary partial dynamical system  and let $Y$ be arbitrary closed set containing $M\setminus \al(\Delta)$. Then there is a Hilbert space  $H$, a \mbox{$C^*$-algebra} $\A\subset L(H)$ whose spectrum is homeomorphic to $M$, and  a partial isometry $U\in L(H)$ such that
\begin{itemize}
\item[i)]  on the spectrum of $\A$ operator $U$ generates the partial mapping $\al$,
\item[ii)]  on the spectrum of $\B:=C^*(\bigcup_{n\in\N}U^{*n}\A U^n)$ operator $U$ generates the  reversible extension $(\M,\tal)$ of $(M,\al)$ associated with $Y$.
\end{itemize}
\end{thm}
\begin{rem} It follows from \cite[Ex. 1.6]{fmr},  see \cite{kwa-leb1}, \cite{kwa-dr}, that concrete $C^*$-dynamical systems  generating a fixed "abstract"   endomorphism $\delta:\A\to\A$ are in one-to-one correspondence with representations of a certain $C^*$-correspondence (Hilbert bimodule) $X$  constructed from $\delta$. Moreover, for  an ideal $J$ in $\A$, the $C^*$-dynamical systems  $(\A,U)$ for which $U^*U\A \cap \A \subset J$ correspond to  \emph{representations of $X$ coisometric on $J$}. %
Consequently, one could construct a pair $(\A,U)$ enjoying the properties described in Theorem \ref{takie zgrabne takie uadne thm} by taking a quotient of the Fock representation of $X$, see \cite{ms}, which in comparison to our construction  is a much more involved approach.
In view of the foregoing observation one could  call \mbox{$Y=\hull(U^*U\A \cap \A)$} a set of \emph{cosurjectivity} for $(\A,U)$. \end{rem}

\subsection{Topological definition,  relation with inverse limits  and universality}\label{relation to inverse limits  and universality}
We finish the section with a discussion of universality of the reversible extension described in Theorem~\ref{opis kosmiczakow}. This theorem shows, in particular,   that  the system  $(\M,\tal)$   is independent of its operator theoretical origin. Therefore   we adopt the following
\begin{defn}\label{naturalne rozszerzeni definicja} Let $(M,\al)$ be a partial dynamical system and  $Y$  a closed subset of $M$ containing  $M\setminus \al(\Delta)$.  The partial dynamical system  $(\M,\tal)$ described in Theorem \ref{opis kosmiczakow} will be called the \emph{natural reversible extension of $(M,\al)$ associated with~$Y$}.
\end{defn}

\begin{rem}\label{naturalne rozszerzeni dyskusja}
In the case when $\al:M\to M$ is defined on
 the whole of $M$ the natural reversible extension  $(\M,\tal)$ associated
 with $Y\supset M\setminus\al(M)$ has the following structure:
\begin{itemize}
\item[-] the system $(M_\infty,\tal)$ coincides with the  inverse limit  system $(\underleftarrow{\,\,\lim\,}(M,\al), \tal)$,
\item[-] for each   $N\in \N$,
  the set  $M_N$ is homeomorphic to $Y$ and   $\tal$ carries homeomorphically  $M_N$ onto $M_{N+1}$:
  $$
 M_0\stackrel{\tal}{\longrightarrow} M_1
\stackrel{\tal}{\longrightarrow} M_2  \stackrel{\tal}{\longrightarrow}\,...\,\,\, (M_\infty,\tal)=
(\underleftarrow{\,\,\lim\,}(M,\al),\tal) .
$$
\end{itemize}
\end{rem}
\begin{rem}
The system   $(\M,\tal)$ has an advantage over the system $(\underleftarrow{\,\,\lim\,}(M,\al), \tal)$. Namely,
   $(\M,\tal)$ is always an extension of  $(M,\al)$ whereas $(\underleftarrow{\,\,\lim\,}(M,\al), \tal)$ may degenerate. The next example illustrates this observation.
 \end{rem}
 \begin{ex}\label{przyklad z odwzorowaniem stalym}
Let $\al:M\to M$ be a constant map with the only value being a non-isolated point   $p\in M$, see Fig. \ref{przyklad z odwzorowaniem stalym1} (a).
  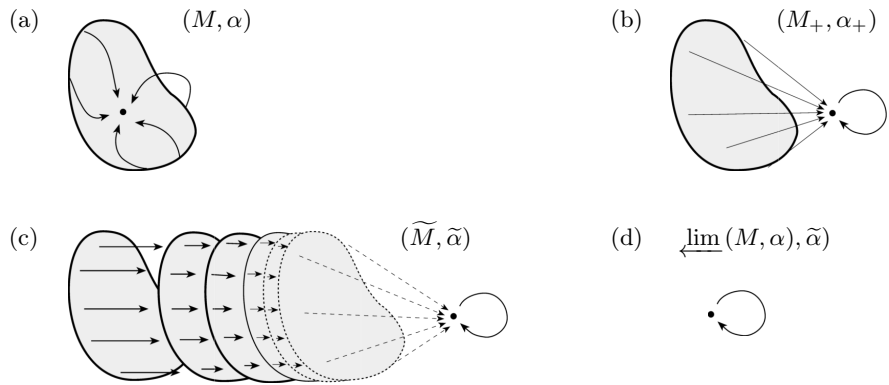
\begin{figure}[htb]
  \begin{center}
\setlength{\unitlength}{1.mm}
\begin{picture}(115,47)(3,-3)
\scriptsize
\put(8,43){\includegraphics[angle=-90, scale=0.8]{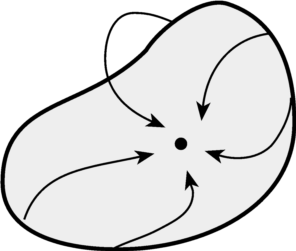}}
\put(88,43){\includegraphics[angle=-90, scale=0.8]{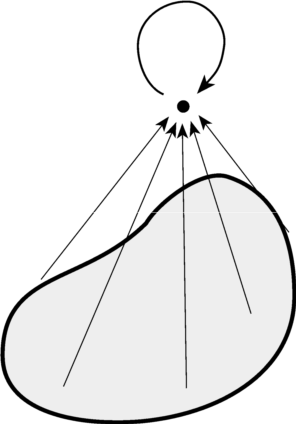}}
\put(8,15){\includegraphics[angle=-90, scale=0.8]{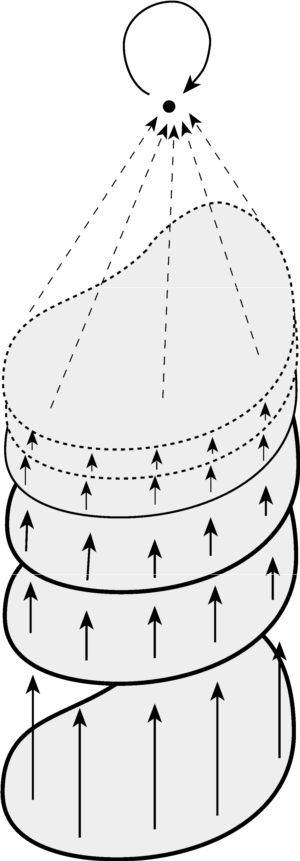}}
\put(93,7){\includegraphics[angle=-90, scale=0.8]{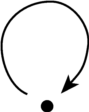}}
%\put(8,43){\includegraphics[angle=-90, scale=0.2]{zbior_x.png}}
%\put(88,43){\includegraphics[angle=-90, scale=0.2]{zbior_x2.png}}
%\put(8,15){\includegraphics[angle=-90, scale=0.2]{zbior_x3.png}}
%\put(93,7){\includegraphics[angle=-90, scale=0.2]{zbior_x4.png}}
%
%Rysunek a)
%podpisy
\put(0,42){(a)} \put(80,42){(b)}
\put(0,13){(c)} \put(80,13){(d)}
\put(23,42){$(M,\al)$}\put(102,42){$(M_+,\al_+)$}
\put(52,13){$(\M,\tal)$}\put(89,13){$\underleftarrow{\,\,\lim\,}(M,\al), \tal)$}
 \end{picture}
\end{center}\caption{Systems associated with a constant map.\label{przyklad z odwzorowaniem stalym1}}
 \end{figure}
The only closed set containing $M\setminus \al(\Delta)=M\setminus\{p\}$ is $Y=M$. For this set the system  $(M_+,\al_+)$ described in Proposition \ref{uklad (M_+, al_+) prop} arise from $(M,\al)$ by adjoining a clopen copy of the singleton $\{p\}$, Fig. \ref{przyklad z odwzorowaniem stalym1} (b). The space   $\M$ is a countable family of copies  of  $M$ compactified
with a single  point   $M_\infty =\{p\}$, Fig. \ref{przyklad z odwzorowaniem stalym1} (c).
Finally, the space $(\underleftarrow{\,\,\lim\,}(M,\al), \tal)$
  degenerates to a single point, Fig. \ref{przyklad z odwzorowaniem stalym1} (d).
\end{ex}
A hard piece of evidence  that $(\M,\tal)$ is an appropriate reversible counterpart of $(M,\al)$   is its identification as a universal minimal object that we now present.
\begin{thm}[Universality of natural reversible extension]\label{universality}
Let $(M,\Delta,\al)$ be a partial dynamical system, $(\M,\TDelta, \tal)$ its  natural reversible extension  associated with a  set $Y\subset M$, and  $\Phi$ the corresponding factor map (\ref{rzut na zerowo wspolrzedno}).
\begin{itemize}
\item[i)] If $(M_\beta,\Delta_\beta,\beta)$ is a reversible extension  of  $(M,\al)$ associated with $Y$, and $\Psi$ is the corresponding semiconjugacy, then there is a unique  semiconjugacy $\TPsi$ from $(M_\beta,\beta)$ to $(\M,\tal)$ such that $\Psi=\Phi \circ \TPsi$, i.e. the diagram
$$
\xymatrix{
(M_\beta,\beta) \ar[rr]^{\TPsi} \ar[rd]_\Psi  & &  (\M,\tal)   \ar[ld]^\Phi \\
 &  (M,\al) &
}
$$
commutes, and    $\TPsi(M_\beta\setminus \beta(\Delta_\beta))= \M \setminus \tal (\TDelta)$, so that $(M_\beta,\beta)$ is  an extension of $(\M,\tal)$ associated with $\M \setminus \tal (\TDelta)$.
\item[ii)] If  a partial dynamical system $(M_\beta,\beta)$  possess the property of $(\M,\tal)$ described in item i), then  $(M_\beta,\beta)$ and $(\M,\tal)$ are  equivalent and equivalence is established by means of $\TPsi$.
\end{itemize}
\end{thm}
\begin{Proof}
 i) Denote by $M_N^\beta$, $N\in \N$, the set of points that has  exactly $N$ preimages under $\beta$, that is  $y\in M_N^\beta$ iff $y, \beta^{-1}(y), ..., \beta^{-(N-1)}(y)\in \beta(\Delta_\beta)$ and $ \beta^{-N}(y)\notin \beta(\Delta_\beta)$. Similarly, we denote by $M_\infty^\beta$  the set of points that have infinitely many preimages under $\beta$. We put
$$
M_N^\beta \ni y \longmapsto \TPsi(y):=(\Psi(y),\Psi(\beta^{-1}(y)),...,\Psi(\beta^{-N}(y)),0,...)\in M_N,
$$
$$
M_\infty^\beta \ni y \longmapsto \TPsi(y):=(\Psi(y),\Psi(\beta^{-1}(y)),...,\Psi(\beta^{-n}(y)),...)\in M_{\infty}.
$$
By \eqref{semicon2} and  \eqref{semicon5} this yields a well defined mapping $\TPsi:M_\beta\to \M$. It is continuous, as for an open set $\widetilde{U}\subset \M$ of the form
$$
\widetilde{U}=\{(x_0,x_1,...)\in \M: x_n \in U\},\quad \textrm{ where } U \textrm{ is open in } M,
$$
the set $\TPsi^{-1}(\widetilde{U})=\beta^{n}(\Psi^{-1}(U))$ is also open.
  Condition \eqref{semicon5} implies that $\TPsi$ maps $M_N^\beta$ onto $M_N$. In particular, $\TPsi(M_\beta\setminus \beta(\Delta_\beta))=\TPsi(M_0^\beta)= M_0=\M \setminus \tal (\TDelta)$. To show  that $\TPsi$ maps $M_\infty^\beta$ onto $M_\infty$ we fix $\x=(x_0,x_1,x_2,...)\in M_\infty$ and  put
$$
D_n=\TPsi^{-1}(\{\y=(y_0,y_1,...) \in \M: y_n=x_n\} ),\qquad n\in \N.
$$
It is evident that $\{D_n\}_{n\in \N}$ forms a decreasing sequence of non-empty compact sets, and therefore $\bigcap_{n\in \N} D_n\neq \emptyset$. Taking $y \in \bigcap_{n\in \N} D_n$ we have $\TPsi(y)=\x$, which proves surjectivity of $\TPsi$.
Relations $
\TPsi^{-1}(\TDelta)=\Delta_\beta
$ and $\Psi=\Phi \circ \TPsi$ are straightforward. \\
For the uniqueness of $\TPsi$ we note that reversibility of the systems $(M_\beta,\beta)$, $(\M,\tal)$ and the equality $\TPsi(M_\beta\setminus \beta(\Delta_\beta))= \M \setminus \tal (\TDelta)$ imply that a semiconjugacy  $\TPsi$  from $(M_\beta,\beta)$ to $(\M,\tal)$ is automatically a semiconjugacy from $(M_\beta,\beta^{-1})$ to $(\M,\tal^{-1})$. This together with relation  $\Psi=\Phi \circ \TPsi$ give a family of relations   $\Psi\circ \beta^{-n}=\Phi \circ \tal^{-n}\circ \TPsi$, $n\in\N$,  understood in the sense that  not only functions but also their natural domains are equal. This forces $\TPsi$ to act according to formulas which we used as a definition in the first part  of the proof.
\\
ii) We have two semiconjugacies $\TPhi:M_\beta\to \M$ and $\TPsi:\M\to M_\beta $ which by the argument from item i) (with the same convention concerning domains) satisfy
$$
\Psi\circ \beta^{-n}=\Phi \circ \tal^{-n}\circ \TPsi,\qquad \Phi \circ \tal^{-n}=\Psi \circ \beta^{-n}\circ \TPhi, \qquad n\in \N.
$$
Thus $(\Phi \circ \tal^{-n})\circ (\TPsi\circ \TPhi)  = \Phi \circ \tal^{-n}$ and by the form of $\M$ it follows that $\TPsi\circ \TPhi=id$.
Therefore $\TPhi$ is a homeomorphism which yields a desired equivalence.
\end{Proof}

Identifying, % like in subsection \ref{Construction of operators generating}, 
$\A=C(M)$ with the  subalgebra $\{a\circ \Phi: a\in C(M)\}$ of $\TA:=C(\M)$,  denoting by $\delta$ and $\tdelta$ endomorphisms corresponding to $\al$ and $\tal$ respectively, and putting  $J=C_Y(M)$,  one may interpret the above result in terms of endomorphisms as follows.
\begin{thm}\label{universality1}
Every partial automorphism $\gamma:\B\to \B$ that extends $\delta:\A\to \A$ in such a way that $(\ker\gamma)^\bot\cap \A= J$
 automatically extends  the partial automorphism $\tdelta:\TA\to\TA$ in such a way that $ (\ker\gamma)^\bot\cap \TA=(\ker\tdelta)^\bot
$. Moreover, this property characterizes the pair $(\TA,\tdelta)$ (up to isomorphisms conjugating endomorphisms).
\end{thm}

\section{Reversible extensions of logistic maps}\label{logistic section}

Now as the necessary preparatory work is implemented and the required \mbox{$C^*$-objects} are described we pass to presentation of calculation of concrete examples of reversible extensions of dynamical systems.

      \begin{figure}[htb]\label{wykres gamma}
      \begin{center}
\setlength{\unitlength}{1.mm}
\begin{picture}(125,65)(0,0)
\put(-7,1){\includegraphics[angle=0, scale=0.32]{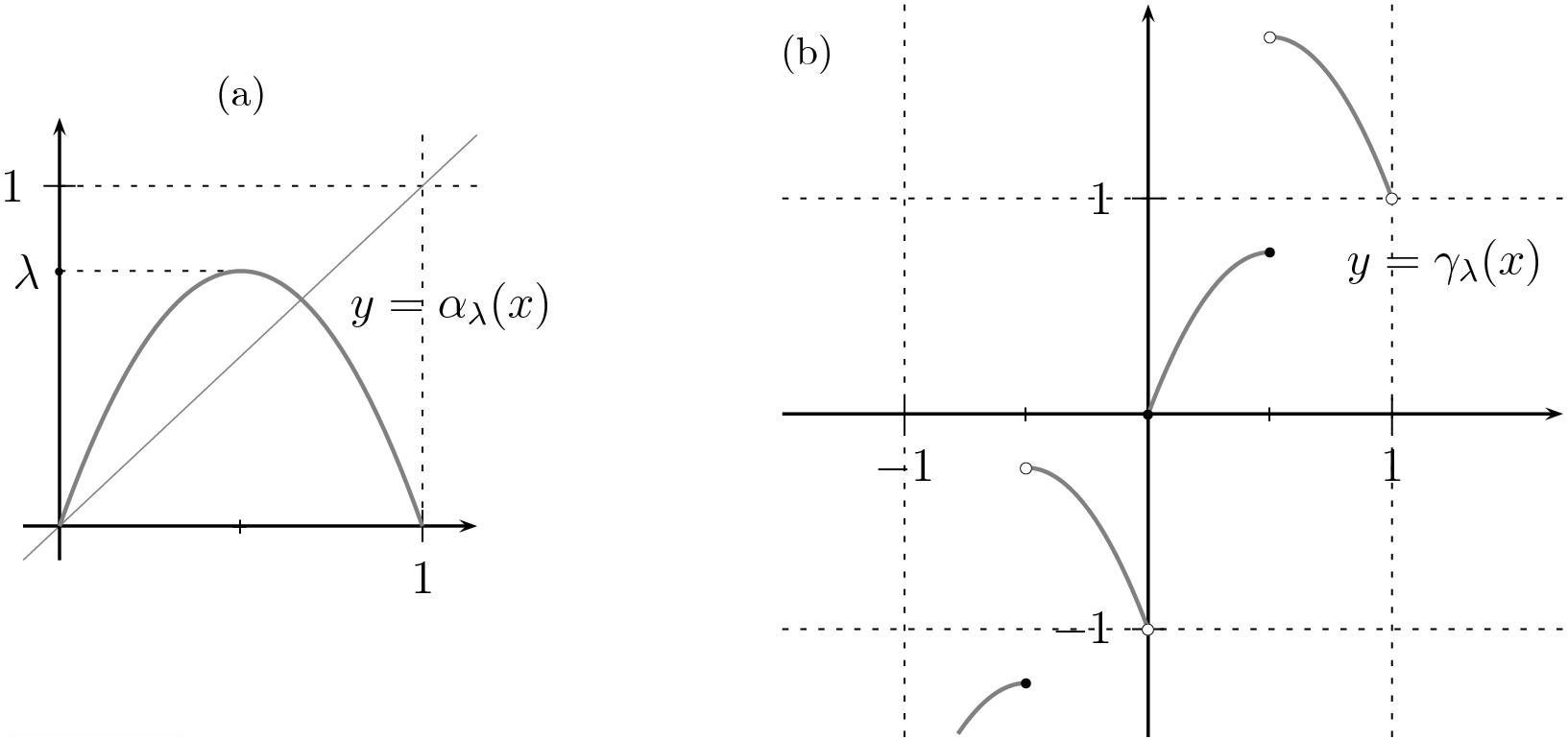}}
%\put(-7,1){\includegraphics[angle=0, scale=0.24]{3.jpg}}
 \end{picture}
   \end{center}
\caption{Graph of the logistic map  $\al_\lambda:[0,1]\to [0,1]$  (a); graph of the injective map $\gamma_\lambda:\R\to \R$ arising from $\al_\lambda$ (b).}
 \end{figure}

 In this section, we conduct  a thorough  analysis  of reversible extensions of  the family of logistic maps  $\{\al_\lambda\}_{\lambda\in (0,1]}$, where by a  \emph{logistic map with parameter $\lambda\in (0,1]$} we mean a quadratic map  $\al_\lambda:[0,1]\to [0,1]$ given by \eqref{odwozorowanie logistyczne label}, cf. Fig. \ref{wykres gamma} (a).  For better illustration of our $C^*$-method we define the operators generating $\al_\lambda$, $\lambda\in (0,1]$, in a concrete fashion, even though we already know that  such operators always exist,  see Theorem \ref{takie zgrabne takie uadne thm}. To this end,  we let $H=L^2(\R)$ and consider the $C^*$-algebra $\A$\label{algebra of periodic functions}  consisting of operators of multiplication by periodic functions $a(t)$ of period $1$, continuous on $[0,1)$, and possessing a  limit in $1$ from below, i.e.:
\begin{equation}\label{warunki na widmo [0,1]}
a(t+1)=a(t),\quad  a|_{[0,1)} \in C([0,1))   \,\, \textrm{ and there exists }\lim_{t\to 1^-}a(t).
\end{equation}
Plainly, $\A$ is  isomorphic to $C([0,1])$ and  we will identify its spectrum  $M$ with the unit interval $[0,1]$:
$$
M=[0,1].
$$
 We fix  $\lambda\in (0,1]$ and define a piecewise continuous function  $\gamma_\lambda:\R\to \R$ where
\begin{equation}\label{gammapierwsza}
\gamma_\lambda(x)=\begin{cases}
4\lambda t(1-t) +2k, & t\in [k,k+\frac{1}{2}]\\
4\lambda t(1-t) +2k+1, & t\in [k+\frac{1}{2},k+1)
\end{cases},\,\,\, k\in \Z.
\end{equation}
The graph of $\gamma_\lambda$ arises from the graph of $\al_\lambda$ by a "propagation"  of the halves of the parabola  $y=\al_\lambda(x)$  on $\R^2$, so  that one obtains a graph of an injective mapping on  $\R$, see Fig. \ref{wykres gamma} (b). The mapping $\gamma_\lambda$ is injective, and  it is   bijective iff $\lambda=1$.

We let   $U_\lambda:L^2(\R)\to L^2(\R)$ be a normalized operator of composition with  $\gamma_\lambda$:%, that is we put
\begin{equation}\label{T_lambda operator defn}
(U_\lambda f)(t)=\sqrt{|\gamma_\lambda'(t)|}f(\gamma_\lambda(t))=2\sqrt{|\lambda(1-2t)|}f(\gamma_\lambda(t))
,\qquad f\in L^2(\R),
\end{equation}
so that $U_\lambda$ is a coisometry where the  adjoint isometry is given by the formula
$$
(U_\lambda^* f)(t)=\begin{cases}
\sqrt{|(\gamma_\lambda^{-1})'(t)|}f(\gamma_\lambda^{-1}(t)), & t\in \gamma_\lambda(\R)\\
0, & t\notin \gamma_\lambda(\R)
\end{cases},\qquad f\in L^2(\R).
$$
In particular,  $U_\lambda$ is unitary if and only if  $\lambda=1$.
\begin{prop}\label{otoczka logistycznych operatorow}
For each $\lambda\in (0,1]$ the operator $U_\lambda$  generates on the spectrum of the $C^*$-algebra $\A$ the logistic map  $\al_\lambda$ %( Definition \ref{definicja operatora kompozycji}).
Moreover the hull of the ideal
$$
J=(U_\lambda^* U_\lambda)\A\cap \A=\{ a\in \A: U_\lambda^* U_\lambda a = a\}
$$
is the smallest possible, cf. Proposition \ref{introductory proposition}, that is
$$
Y=\hull(J)=\overline{M\setminus \al_\lambda(M)}=\begin{cases}
[\lambda,1] & \textrm{ if } \lambda <1,
\\
\emptyset  & \textrm{ if } \lambda =1.
\end{cases}
$$
\end{prop}
\begin{Proof}
Since $U_\lambda^* U_\lambda$ is the operator of multiplication by the characteristic function $\chi_{\gamma_\lambda(\R)}$ of the set $\gamma_\lambda(\R)$, we have $U_\lambda^* U_\lambda \in \A' $ and the hull of $J$
is the set $
 Y=  \hull(\{ a\in \A: \chi_{\gamma_\lambda(\R)}a=a\})=[\lambda,1]$ for $\lambda <1$  and  $Y=\emptyset$  for $\lambda =1$. For any  operator $a$ of multiplication by  $a(t)$,  $U_\lambda a U_\lambda^*$ is an operator of multiplication by
$$
(U_\lambda a U_\lambda^*)(t)=a (\gamma_\lambda(t))=a(\mantysap\gamma_\lambda(t)\mantysak)=a(\al_\lambda(\textbf{\{}t\textbf{\}})),
$$
where $\textbf{\{}t\textbf{\}}\in [0,1)$  denotes the fractional part of a number $t\in \R$. Hence one sees that $U_\lambda$ generates  on $[0,1]$ the mapping $\al_\lambda$.
\end{Proof}
In view of the above  the first of the following operations
  $$
 \delta_\lambda(a):=U_\lambda a U_\lambda^*,\qquad  \delta_{*,\lambda}(a):=U_\lambda^* aU_\lambda,\qquad a\in L(L^2(\R))
  $$
preserves the algebra $\A$ and the logistic map $\al_\lambda$ is its dual map.  In particular, we may adopt the identifications:
$$
\A=C([0,1]),\qquad \delta_\lambda(a)=a\circ \al_\lambda,\quad a\in \A.
$$
% $$ (\A,\delta_\lambda)\, \longleftrightarrow \,([0,1],\al_\lambda).$$
On the other hand, for $a\in \A$, $\delta_{*,\lambda}(a)$ is the operator of multiplication by the function
$$
\delta_{*,\lambda}(a)(t)=(U_\lambda^* a U_\lambda)(t) =\begin{cases}
a(\gamma_\lambda^{-1}(t)), & x\in \gamma_\lambda(\R)\\
0, & x\notin \gamma_\lambda(\R)
\end{cases},
$$
which is periodic but in general its period is two, not one. Therefore  the mappings
$
 \delta_{*,\lambda}
  $, $\lambda \in (0,1]$,
do not preserve the algebra  $\A$ and the $C^*$-algebras
\begin{equation}\label{B_lambda algebra defn}
 \B_\lambda:= C^*\big(\bigcup_{n\in\N} U_\lambda^{*n} \A U_\lambda^n\big),
 \end{equation}
   are essentially bigger than $\A$.
 \begin{thm}\label{stw widmo B_lambda algebry}
Let $\lambda \in (0,1]$ and let  $(\M_\lambda, \tal_\lambda)$  be the reversible extension of  $([0,1], \al_\lambda)$ associated with the set
$
Y
$ where  $
Y=
[\lambda,1]$, if $\lambda <1$, and $Y=\emptyset$, if $\lambda=1$.  The algebra $\B_\lambda$ may be identified with  $C(\M_\lambda)$ and then the partial automorphism  $\delta_\lambda:\B_\lambda\to \B_\lambda$  becomes the operator of composition with the  homeomorphism  $\tal_\lambda:\M_{\lambda}\to \tal_\lambda(\TDelta)$:
$$
\B_\lambda =C(\M_\lambda),\qquad \delta_\lambda(a)=a\circ \tal_\lambda,\quad a\in \B_\lambda.
$$
In particular operator $U_\lambda$ generates on  $\M_\lambda$ the partial homeomorphism $\tal_\lambda$.
 %Przestrzeń ideałów maksymalnych algebry $\B_\lambda$ możemy utożsamić z $\M_\lambda$, a  $\tal_\lambda$ z odwzorowaniem dualnym do endomorfizmu $\delta_\lambda:\B_\lambda\to \B_\lambda$.
% $$ (\B_\lambda,\delta_\lambda)\, \longleftrightarrow \,(\M_\lambda,\tal_\lambda).$$
\end{thm}
\begin{Proof} It follows from  Proposition \ref{otoczka logistycznych operatorow}, Theorem \ref{opis kosmiczakow} and  Definition  \ref{naturalne rozszerzeni definicja}.
\end{Proof}
The description of the extended $C^*$-dynamical systems $(\B_\lambda,\delta_\lambda)$ reduces to the description of the family of reversible topological systems
$$
\{(\M_\lambda,\tal_\lambda)\}_{\lambda \in (0,1]}.
$$
 We recall, cf. Remark \ref{mapping Phi remark}, that the mapping  $\Phi:\M_\lambda\to [0,1]$ dual to the embedding $\A\subset \B_\lambda$ is a surjection such that the diagram
$$
\begin{CD}
\M_\lambda      @>{\tal_\lambda}>>  \M_\lambda    \\
@V{\Phi}VV               @VV{\Phi}V\\
[0,1] @>{\al_\lambda}>>  [0,1]
\end{CD}
$$
commutes. We stress that  a change of the parameter value $\lambda\in (0,1]$ does not only influence the dynamics of $\tal_\lambda$ but also the topology of the space  $\M_\lambda$. The following notions of continuum theory will  be indispensable in our analysis
 \begin{defn}[see, for instance,  \cite{Nadler}]\label{definicje continuuow}
By a \emph{continuum} we mean a connected and compact  metric space.  A continuum will be called
\begin{itemize}
\item[i)] \emph{nondegenerate}, if it is not a singleton,
\item[ii)] \emph{reducible},  if it may be presented as the sum of two its proper subcontinuua,
\item[iii)] \emph{irreducible}, if it is a  nondegenerate continuum which is not  reducible,
\item[iv)] \emph{snake-like} or \emph{arc-like continuum}, if it is homeomorphic to an inverse limit of an inverse sequence with bonding maps being continuous maps of an interval, cf. \cite[12.19]{Nadler}.
\end{itemize}
\end{defn}
For all $\lambda\in (0,1]$, the spectrum $\M_\lambda$ of $\B_\lambda$ contains the snake-like continuum   $M_\infty=\underleftarrow{\,\,\lim\,}([0,1],\al_\lambda)$, and it is a general dynamical principle, discovered by    M. Barge and J. Martin \cite{Barge-Martin}, that one should expect this continuum to  be irreducible, see Theorem \ref{bifurkacji ciag pierwszy} below.  %One of the most known irreducible continua and a non-trivial example of a snake-like continuum is the Brouwer-Janiszewski-Knaster continuum presented on Fig.  \ref{buckethandle}.

\subsection{General structure of the extended systems $(\M_\lambda,\tal_\lambda)$}\label{struktura logistycznego rozszerzenia}
Let $\lambda < 1$. By  Remark \ref{naturalne rozszerzeni dyskusja} the   space $\M_\lambda$ consists of the set $M_\infty$ being the inverse limit $\underleftarrow{\,\,\lim\,}(M,\al_\lambda)$ and a countable family of \emph{arcs}
$
M_N
$, that is  sets homeomorphic to a closed interval: $M_N\cong [\lambda,1]$.
\begin{thm}\label{sprowadzenie do granicy odwrotnej}
For  $\lambda< 1$ the maximal ideal space  $\M_\lambda$ of algebra $\B_\lambda$ consists of the snake-like continuum  $M_\infty =\underleftarrow{\,\,\lim\,}(M,\al_\lambda)$ and a sequence of arcs $M_N$ such that
$$
\lim_{n\to \infty} M_N =M_\infty,
$$
where the limit is taken in  Hausdorff metric. In particular,
$$
\M_\lambda=\bigcup_{n\in \N} M_N\cup M_\infty=
\overline{\bigcup_{n\in \N} M_N}.$$
The mapping $\tal_\lambda$ generated by $U_\lambda$ on $\M_\lambda$ carries homeomorphically arc $M_N$ onto arc $M_{N+1}$, and on the continuum  $M_\infty=\underleftarrow{\,\,\lim\,}(M,\al_\lambda)$ it coincides with the homeomorphism induced by $\al_\lambda$ (Definition \ref{definicja granicy odwrotnej}).
\end{thm}
\begin{Proof} We only  need to show the equality  $
\lim_{n\to \infty} M_N =M_\infty
$ which  (see for instance \cite[Thm. 4.11]{Nadler}) is equivalent to the following two inclusions
$$
\lim \sup M_N=\{\x\in \M_\lambda: \forall_{\x\in U \textrm{ open}}\,\,
 \forall_{k\in \N}\,\, \exists_{N>k}\,\, U\cap M_N \neq \emptyset\} \subset M_\infty,
$$
$$
M_\infty \subset \lim \inf M_N=\{\x\in \M_\lambda: \forall_{\x\in U \textrm{ open}}\,\,
 \exists_{k\in \N}\,\, \forall_{N>k}\,\, U\cap M_N \neq \emptyset\}.
$$
If we assume that $\x\in \lim \sup M_N$ and  $\x \in M_{N_0}$, for certain $N_0\in \N$, then taking in the definition of $\lim \sup M_N$,  $U=M_{N_0}$ and $k=N_0$ we arrive   at a contradiction. This proves the first inclusion.
To prove the second one take  $\x=(x_0,...,x_{N_0},...) \in M_\infty$ and an open neighbourhood of $\x$ of the form
$$
U=\{\y=(y_0,...,y_{N_0},...)\in \M_\lambda: y_{N_0}\in (x_{N_0}-\varepsilon ,  x_{N_0} + \varepsilon )\}.
$$
Note that there exists $n_0\in \N$ such that  $\al_\lambda^{n_0}([\lambda,1])=[0,\lambda]$. Indeed, if
$\lambda \leq \frac{1}{2}$, then $\al_\lambda([\lambda,1])=[0,\lambda]$, and if $\lambda > \frac{1}{2}$, then for certain $n_0$, $\al_\lambda^{n_0}(\lambda) >\frac{1}{2}$, and therefore $\al_\lambda^{n_0}([\lambda,1])=[0,\lambda]$.
Putting $k={N_0}+{n_0}$ we have  $U\cap M_N \neq \emptyset$ for every $N> k$,  and  thus $\x \in  \lim \inf M_N$.
\end{Proof}
The above statement implies  that once we have at our disposal description of  $(M_\infty, \tal_\lambda)$ we may easily   describe  the whole system $(\M_\lambda, \tal_\lambda)$: it suffices to adjoin to $M_\infty$ the sequence of arcs $\{M_N\}_{N\in \N}$ converging to $M_\infty$ and prolong  $\tal_\lambda$ so that it shifts homeomorphically arcs $M_N$ towards $M_\infty$. The importance of this comment lies in that a great deal of facts concerning the systems of $(\underleftarrow{\,\,\lim\,}(M,\al_\lambda),\tal_\lambda)$   type   is known, see  \cite{Barge-Ingram}. Thus we may use them to achieve our goal.% get a completepicture  of the systems $(\M_\lambda, \tal_\lambda)$, $\lambda\in (0,1]$.
%Tym samym  dalsza część niniejszego podrozdziału w dużej mierze sprowadzać  się będzie do nowej algebraiczno-operatorowej interpretacji  faktów .

  \subsubsection{The extended system for $\lambda=1$ (B-J-K continuum).}
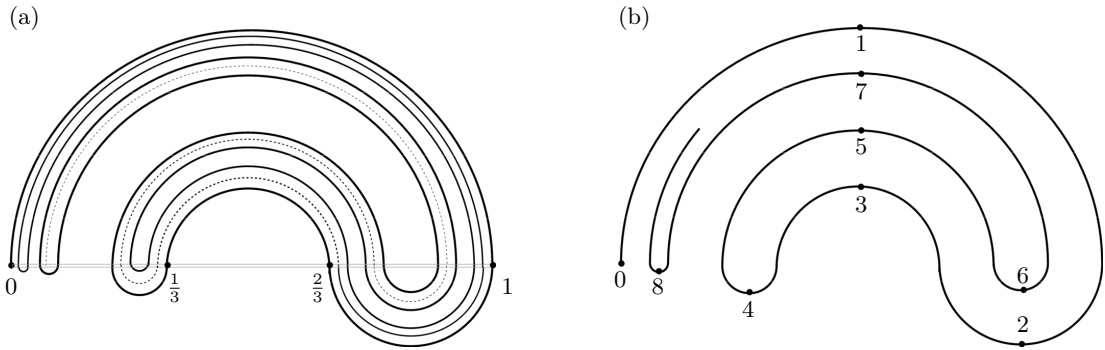
\begin{figure}[htb]
\begin{center}\setlength{\unitlength}{1mm}
\begin{picture}(110,46)(0,-2)
\scriptsize
\put(63,-3){\includegraphics[angle=0, scale=0.8]{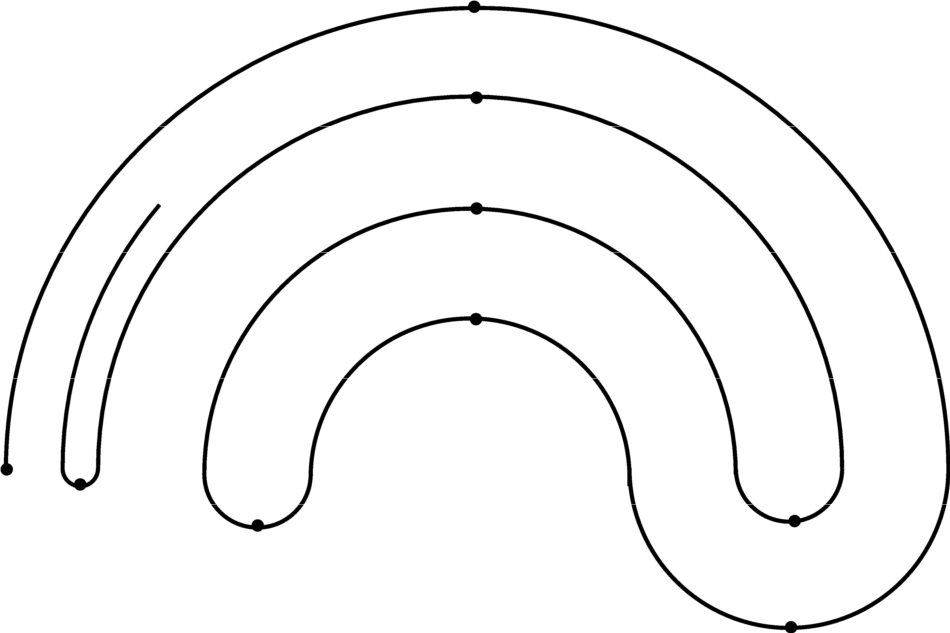}}
\put(-18,-3){\includegraphics[angle=0, scale=0.8]{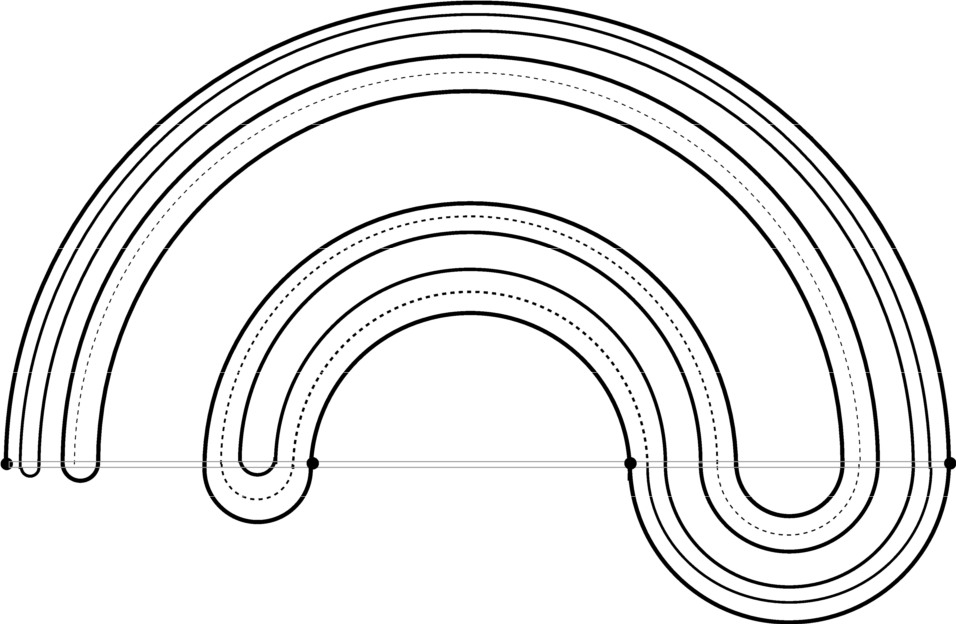}}
%\put(63,-3){\includegraphics[angle=0, scale=0.19]{buckethandle2.png}}
%\put(-18,-3){\includegraphics[angle=0, scale=0.19]{buckethandle1.png}}
\put(-18,40){(a)}\put(63,40){(b)}
\put(-18.5,4){$0$}\put(3,4){$\frac{1}{3}$}\put(22,4){$\frac{2}{3}$}\put(47.5,4){$1$}
\put(62.5,5){$0$}\put(67.5,4){$8$}\put(79.5,1){$4$}\put(116,5.5){$6$}\put(116,-1){$2$}\put(94.5,36.3){$1$}
\put(94.5,30){$7$}\put(94.5,22.7){$5$}\put(94.5,15){$3$}
\end{picture}
\end{center}
\caption{Brouwer-Janiszewski-Knaster continuum. \label{buckethandle}}
 \end{figure}

For  $\lambda=1$ the  mapping $\al_1$ is a surjection and the space  $\M_1$ coincides with the inverse limit $\underleftarrow{\,\,\lim\,}([0,1],\al_1)$ of the full logistic map $\al_1(x)=4x(1-x)$, cf. Theorem \ref{szczególne przypadki} iv). Hence   $\M_1$ is one of the most famous irreducible continuum called   \emph{Brouwer-Janiszewski-Knaster continuum},  briefly \emph{B-J-K continuum} \cite{Barge-Ingram},  \cite{Nadler}, \cite{Watkins}.  One may graphically depict  $\M_1$   by joining the points of the Cantor set with semicircles in a manner presented on Fig. \ref{buckethandle} (a).
  \\
The logistic map $\al_1(x)=4x(1-x)$ is topologically conjugate to the tent map $\al_T(x)= 1-|2x-1|$ and    B-J-K continuum is usually considered  \cite{Nadler},  \cite{Watkins} as the inverse limit $\underleftarrow{\,\,\lim\,}(M,\al_T)$.  Then there is a natural parametrization of the composant of the point $(0,0,0,...)$ of B-J-K continuum by non-negative real numbers \cite{Watkins}, see Fig. \ref{buckethandle} (b). Within this parametrization the induced homeomorphisms  $\tal_1$ and $\tal_T$ fulfill the following formulae
\begin{equation}\label{postac tal dla lambda=1}
\tal_1(t)=\begin{cases}
 2k +\al(\textbf{\{} t\textbf{\}}),&  t\in [k,k+\frac{1}{2})\\
2(k+1) -\al(\textbf{\{} t\textbf{\}}),&  t\in [k+\frac{1}{2},k+1)
\end{cases}, \qquad\tal_T(t)=2t,
\end{equation}
where $\textbf{\{} t\textbf{\}}$ denotes the fractional part of the number $t$ (obviously, the systems $(\M_1, \tal_1)$ and $(\M_1, \tal_T)$ are topologically conjugate). These comments give a certain idea about the dynamics of the system $(\M_\lambda, \tal_\lambda)$ for $\lambda=1$.

\begin{thm}\label{takie tam twierdzenie dla lambda=1}
Algebra $\B_1$ may be identified with the algebra  $C(\M_1)$ of continuous functions on B-J-K continuum $\M_1$, Fig.  \ref{buckethandle} (a).  Then the automorphism $\delta_1:\B_1\to\B_1$ becomes the operator of composition with the homeomorphism $\tal_1:\M_1\to \M_1$, which within  the parametrization presented on  Fig.  \ref{buckethandle} (b) assumes the form \eqref{postac tal dla lambda=1}.
\par
Furthermore irreducibility of $\M_1$ expresses as the following property of the algebra $\B_1$:  if $I_1$, $I_2$ are ideals in $\B_1$ such that  $\B_1/I_{i}$, for $i=1,2$, does not contain non-trivial idempotents, then
$$
I_1\cap I_2=\{0\} \,\, \, \Longrightarrow\,\,\, I_1=\{0\} \textrm{ or } I_2=\{0\}.
$$
\end{thm}
\begin{Proof}
The first part the statement follows from Theorem \ref{stw widmo B_lambda algebry}. To show the second  part  notice that for $i=1,2$ we have
$$
I_i=C_{Y_i}(\M_1),\qquad \textrm{ where }\,\, Y_i=\hull(I_i)\subset \M_1.
$$
Algebra $C(\M_1)/I_i\cong C(Y_i)$ does not contain non-trivial idempotents if and only if  $Y_i$ is connected, that is if $Y_i$ is a subcontinuum of  $\M_1$. Since  $I_1\cap I_2=C_{Y_1\cup Y_2}(\M_1)$ the condition $I_1\cap I_2=\{0\}$ is equivalent to the equality $Y_1\cup Y_2 =\M_1$ and therefore the asserted property of $\B_1$ is equivalent to irreducibility of  continuum $\M_1$.  \end{Proof}

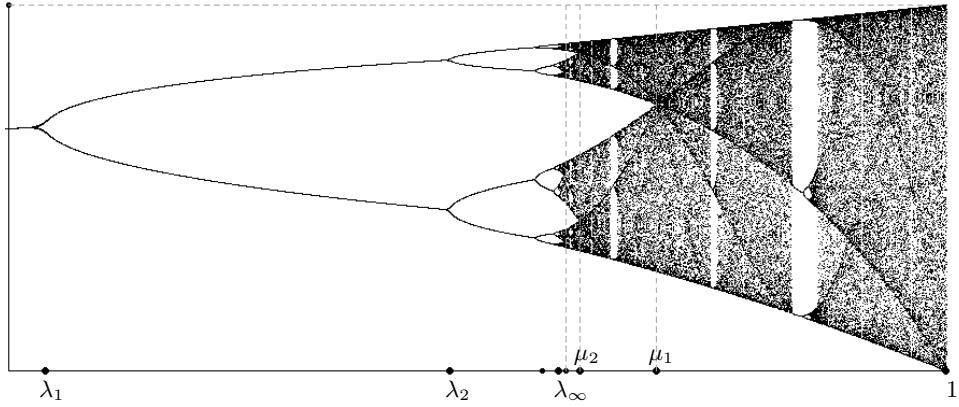
\begin{figure}[htb]
\begin{center}\setlength{\unitlength}{0.9mm}
\begin{picture}(240,53)(-6,5)
\scriptsize
\put(5,-1.5){\includegraphics[angle=0, scale=0.405]{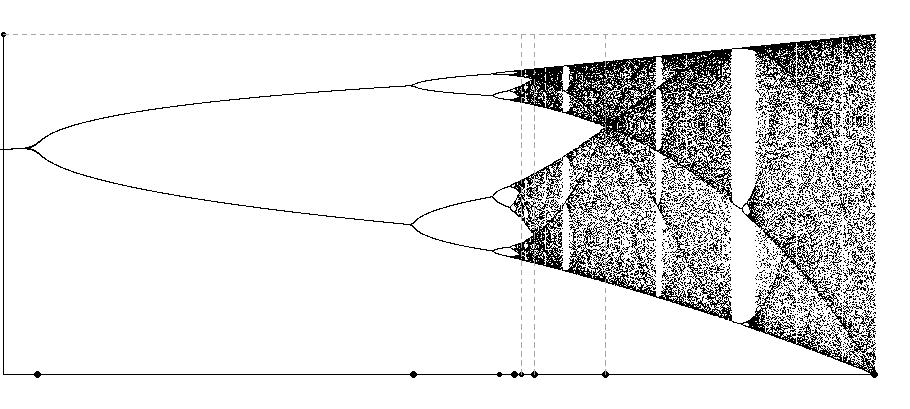}}
\put(10,0){$\lambda_1$} \put(70,0){$\lambda_2$} \put(86,0){$\lambda_\infty$} \put(89,5.5){$\mu_2$} \put(100,5.5){$\mu_1$} \put(144,0){$1$}
   \end{picture}
\end{center}
\caption{Bifurcation diagram for $\lambda \in [0,74;1]$. \label{Diagram Bifurkacyjny1}}
 \end{figure}

\subsubsection{Feigenbaum limit $\lambda_\infty$}\label{paragraf o Feigenbaumach}
Let  $\lambda_1=\frac{3}{4}$, $\lambda_2=\frac{1+\sqrt{6}}{4}\approx 0.862$, $\lambda_3$, ..., be the sequence of the parameter values $\lambda$ corresponding to the first cascade of period-doubling bifurcation of the system  $([0,1],\al_\lambda)$. This is an  increasing and the value  $\lambda_\infty=\lim_{n\to \infty}\lambda_n\approx 0.89249$ is called \emph{Feigenbaum limit} \cite{Barge-Ingram},  \cite{Col-Eck}, \cite{Devaney}. Interval of parameter splits into two parts $(0,\lambda_\infty)$ and $(\lambda_\infty,1]$  which correspond respectively to regular and chaotic behavior of the systems $([0,1],\al_\lambda)$, see Fig. \ref{Diagram Bifurkacyjny1}. This finds a splendid reflection in the structure of the algebra $\B_\lambda$.
   \begin{thm}\label{bifurkacji ciag pierwszy}
Let  $\lambda\in (0,1]$. The following conditions are equivalent
\begin{itemize}
\item[i)] $\lambda > \lambda_\infty$,
\item[ii)] spectrum  $\M_\lambda$ of  $\B_\lambda$ contains an irreducible  continuum,
\item[iii)] There exists a proper ideal $J$ in a maximal ideal in the $C^*$-algebra $\B_\lambda$ with the property  that for every two ideals $I_1,I_2$ in $\B:=\B_\lambda/J$ such that $\B/I_{1}$, $\B/I_{2}$ do not contain non-trivial idempotents we have
$$
I_1\cap I_2=\{0\} \,\, \, \Longrightarrow\,\,\, I_1=\{0\} \textrm{ or } I_2=\{0\}.
$$
\end{itemize}
  \end{thm}
 \begin{Proof} Equivalence of i) and ii) follows from Theorem  \ref{sprowadzenie do granicy odwrotnej} and \cite[Thm. 3, 4, 7]{Barge-Ingram}. In the proof of Theorem   \ref{takie tam twierdzenie dla lambda=1} we have shown that condition in item iii) is equivalent to irreducibility of  continua. Moreover,  the statement "$J$ is the proper ideal in a maximal ideal in  $\B_\lambda$" means that "the set $Z=\hull(J)$ contains more than one point". This explains  equivalence of ii) and iii).
 \end{Proof}
%Dla wartości parametru $\lambda$ leżącego w przedziale $(0, \lambda_\infty)$  sytuacja jest całkowicie zrozumiała. Wybrane wyniki dotyczące układu $(M_\infty,\tal_\lambda)$ dla $\lambda > \lambda_\infty$ przedstawimy  w podrozdziale \ref{knaster coninuua}.

  \subsection{The first cascade of bifurcation: $\lambda \in (0, \lambda_\infty$]}\label{paragraf o bifurkacjach}

 For $\lambda< \lambda_\infty$  the dynamics of  $([0,1],\al_\lambda)$ is completely understood:  $\al_\lambda$ has exactly one stable orbit which is periodic with period  $2^n$,   $n\in \N$, exactly one repelling periodic orbit with period $2^k$, for each $k=0,...,n-1$, and at most two repelling fixed points. Increasing  $\lambda$ from  $0$ to $\lambda_\infty$ the number  $n$ gradually increases - the system undergoes a \emph{period-doubling bifurcation}. In particular the period of the stable orbit of $\al_\lambda$, for $\lambda\in (0,\lambda_\infty)$, increases according to \emph{Sharkovskii's  order} \cite{Brin}, \cite{Devaney}:\label{Sharkovskii's  order}
$$
1 \triangleleft 2\triangleleft 4 \triangleleft ... \triangleleft 2^n \triangleleft  ...
$$
$$
...\triangleleft\, 2^m(2n+1) \triangleleft ...   \triangleleft  2^m\cdot 7 \triangleleft  2^m\cdot 5 \triangleleft 2^m\cdot 3 \triangleleft  ...
$$
$$
...\triangleleft\, 2(2n+1) \triangleleft ...   \triangleleft 14 \triangleleft 10 \triangleleft 6 \triangleleft  ...
$$
$$
... \triangleleft\,(2n+1) \triangleleft ...   \triangleleft 7  \triangleleft 5  \triangleleft 3.
$$
Let us now imaging that we slowly move the parameter  $\lambda$ from $0$ to $\lambda_\infty$, and simultaneously observe the maximal ideal space
$$
\M_\lambda=\bigcup_{n\in \N} M_N\cup M_\infty
$$
 of the algebra $\B_\lambda$. Let also assume that we watch the change in $\M_\lambda$ from the point of view of the initial space $M=[0,1]$ -- looking at a point  $\x=(x_0,x_1...)\in \M_\lambda$ we  read of   $x_0\in [0,1]$. Such an approach will allow us to understand in detail   how the change of the parameter $\lambda$ affects the topology of the space $\M_\lambda$.
In other words we will built an image of  $\M_\lambda$ with the help of the factor map
$$
\Phi(x_0,x_1,...)=x_0,\qquad (x_0,x_1,...)\in \M_\lambda,
$$
cf. Remark  \ref{mapping Phi remark}. Since however   $\Phi$ is not injective,  it may  "wind"  a piece of an arc $M_N$ several times  onto an interval. Thus in order to get a homeomorphic image of $\M_\lambda$  one needs to modify $\Phi$ by   "adjoining"   copies of the corresponding arcs. In this manner one  constructs a homeomorphism from $\M_\lambda$ onto a certain subset of the plane, cf. \cite{kwa2}. Here, we restrict ourselves to discussion of  the results of  that procedure, the trace of which will be seen on pictures of $\M_\lambda$ where  a point $\x\in\M_\lambda$ is labeled by its zeroth coordinate  $\Phi(\x)\in M$. The crucial role in this enterprise is played by the orbit $\{q_n\}_{n\in\N}$ of the critical point of the mapping   $\al_\lambda$.
We set
$$
q_n:=\al_\lambda^n\Big(\frac{1}{2}\Big),\qquad n\in \N.
$$
  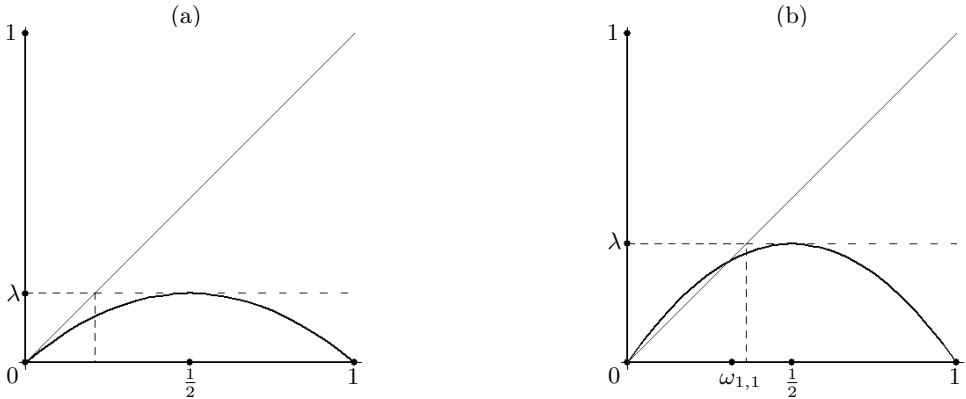
\begin{figure}[htb]
 \begin{center}\setlength{\unitlength}{1mm}
\begin{picture}(120,45)(5,60)
%\put(-8,115){a)}
%lambdy
\scriptsize
\put(20,106){(a)}
 \put(100.5,106){(b)}

  %pikczury
\put(0,60){\includegraphics[angle=0, scale=0.33]{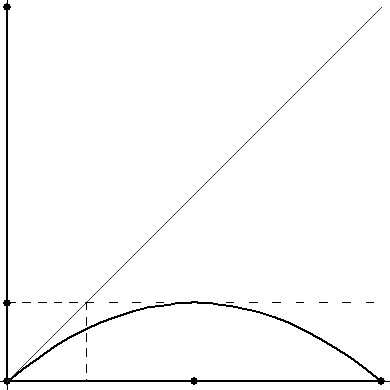}}
 \put(80,60){\includegraphics[angle=0, scale=0.33]{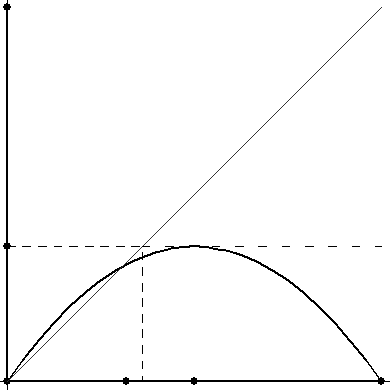}}

%punkty

%zera
\put(-1.8,58.2){$0$}
\put(78.2,58.2){$0$}
%jedynki

\put(-2,103.8){$1$}\put(43.5,58){$1$}
\put(78,103.8){$1$}\put(123.5,58){$1$}
%jedne drugie
\put(21.6,57.4){$\frac{1}{2}$}\put(101.6,57.4){$\frac{1}{2}$}\put(93,58){$\omega_{1,1}$}
%q_1 i q_2
\put(-1.8,69){$\lambda$}
\put(78.1,76){$\lambda$}
\end{picture}
\end{center}
\caption{Graph of  $\al_\lambda$  for  $0<\lambda \leq\frac{1}{4}$ (a);   $\frac{1}{4}<\lambda \leq\frac{1}{2}$ (b).\label{lambda < 1/2}}
 \end{figure}
  \begin{figure}[htb]
\begin{center}
\setlength{\unitlength}{1.12mm}
\begin{picture}(100,49)(-8,-5)
\thicklines

%mięso
\scriptsize
\put(10,45){(a)}\put(80,45){(b)}
%M_0
\qbezier(0,-5)(0,-5)(30,-5)
%M_1
\qbezier(-10,9)(0,9)(0,9)\qbezier(5,12)(4,10)(0,9)
%M_2
\qbezier(-10,19)(-8,19)(-2,19)\qbezier(2.8,21.2)(2,19.5)(-2,19)
%M_3
\qbezier(-10,26)(-8,26)(-4,26)\qbezier(0,27.8)(-1,26.5)(-4,26)
%M_4
\qbezier(-10,31)(-8,31)(-6,31)\qbezier(-3,32.6)(-3.8,31.5)(-6,31)

%punkty
\thinlines
\scriptsize
%na M_0
\put(-0.5,-8){$q_1$}
\put(9.5,-8.4){$\frac{1}{2}$}
\put(29.5,-8){$1$}
\put(0,-5){\circle*{0.8}}
\put(10,-5){\circle*{0.8}}
\put(30,-5){\circle*{0.8}}
%na M_1
\put(-12,7){$0$}
\put(-2.5,6.5){$q_1$}
\put(6,11){$q_2$}
\put(-10,9){\circle*{0.8}}
\put(0,9){\circle*{0.8}}
\put(5,12){\circle*{0.8}}
%na M_2
\put(-12,17){$0$}
\put(-4.5,16.5){$q_2$}
\put(4,20){$q_3$}
\put(-10,19){\circle*{0.8}}
\put(-2,19){\circle*{0.8}}
\put(2.8,21.2){\circle*{0.8}}
%na M_3
\put(-12,24){$0$}
\put(-5.8,23.5){$q_3$}
\put(1.7,27.2){$q_4$}
\put(-10,26){\circle*{0.8}}
\put(-4,26){\circle*{0.8}}
\put(0,27.8){\circle*{0.8}}
%na M_4
\put(-10,31){\circle*{0.8}}
\put(-6,31){\circle*{0.8}}
\put(-3,32.6){\circle*{0.8}}
%na M_\infty
\put(-12,36.5){$0$}
\put(-10,39){\circle*{0.8}}
%nieskończoność
\thinlines
\qbezier[8](-10,32)(-10,35)(-10,38)
\qbezier[8](-6.5,32)(-7.95,35)(-9.4,38)
\qbezier[8](-4,33.6)(-6.4,36)(-8.8,38.4)
%strzałki
%z M_0 na M_1
\qbezier(0.2,-4)(4,3.5)(5,10)\put(5,9.9){\vector(0,1){1}}
\qbezier(9.8,-4)(8,2)(1.4,7.6)\put(1.5,7.5){\vector(-1,1){1}}
\qbezier(29,-4)(8,-1)(-8.6,8)\put(-8.6,8){\vector(-2,1){1}}
%z M_1 na M_2
\put(-10,10.2){\vector(0,1){7.5}}
\qbezier(-0.2,10)(-0.4,14)(-1.4,17.7)\put(-1.4,17.7){\vector(-1,3){0.2}}
\put(5,13){\vector(-1,4){1.8}}
%z M_2 na M_3
\put(-10,20.2){\vector(0,1){4.6}}
\put(-2,20){\vector(-1,3){1.7}}
\put(2.4,22){\vector(-1,3){1.7}}
%z M_3 na M_4
\put(-10,27.2){\vector(0,1){2.7}}
\put(-4.3,27){\vector(-1,3){1.1}}
\put(-0.6,28.8){\vector(-1,2){1.6}}
%z M_infty na M_infty
\put(-11,41){\oval(3,3)[l]} \put(-11,41){\oval(3,3)[tr]}
\put(-9.2,42.2){\vector(-1,-4){0.6}}
%podpisy
\small
\put(-21,-5.5){$M_0$}
\put(-21,8.5){$M_1$}
\put(-21,18.5){$M_2$}
\put(-21,25.5){$M_3$}
\put(-21,30.5){$M_4$}
\put(-21,38.5){$M_\infty$}

%Napis

%\put(37.5,40){\vector(0,-3){35}} \put(39,37){p}\put(39,34){u}\put(39,31){n}\put(39,28){k} \put(39,25){t}\put(39,20){s}\put(39,17){t}\put(39,14){a}\put(39,11){ł}\put(39,8){y}
\end{picture}
%DRUGA SPRĘŻYNKA                          % D R U G A    S P R Ę Ż Y N K A
\begin{picture}(100,0)(-12,-10)
\thicklines
%M_0
\qbezier(71,-5)(71,-5)(95,-5)
%M_1
\qbezier(55,9)(72,9)(71,9)\qbezier(71,9)(72,9)(74,10.5)
%M_2
\qbezier(55,19)(70,19)(70,19)\qbezier(70,19)(71,19)(72.7,20.5)
%M_3
\qbezier(55,26)(70,26)(69,26)\qbezier(69,26)(69,26)(71.2,26.6)
%M_4
\qbezier(55,31)(69,31)(68,31)\qbezier(68,31)(68,31)(70.1,31.4)
%M_\infty
\qbezier(55,39)(65,39)(64,39)
%punkty
\thinlines
\scriptsize
%na M_0
\put(70.5,-8){$q_1$}
\put(74.5,-8.4){$\frac{1}{2}$}
\put(94.5,-8){$1$}
\put(71,-5){\circle*{0.8}}
\put(75,-5){\circle*{0.8}}
\put(95,-5){\circle*{0.8}}
%na M_1
\put(53,7){$0$}
\put(68.5,6.6){$q_1$}
\put(75,10){$q_2$}
\put(55,9){\circle*{0.8}}
\put(71,9){\circle*{0.8}}
\put(74,10.5){\circle*{0.8}}
%na M_2
\put(53,17){$0$}
\put(67.5,16.5){$q_2$}
\put(74,20){$q_3$}
\put(55,19){\circle*{0.8}}
\put(70,19){\circle*{0.8}}
\put(72.7,20.5){\circle*{0.8}}
%na M_3
\put(53,24){$0$}
\put(66.5,23.6){$q_3$}
\put(74,27.2){$q_4$}
\put(55,26){\circle*{0.8}}
\put(69,26){\circle*{0.8}}
\put(71.2,26.6){\circle*{0.8}}
%na M_4
\put(55,31){\circle*{0.8}}
\put(68,31){\circle*{0.8}}
\put(70.1,31.4){\circle*{0.8}}
%na M_\infty
\put(55,39){\circle*{0.8}}
\put(64,39){\circle*{0.8}}
\put(53,36.5){$0$}
\put(65,36.9){$\omega_{1,1}$}
%nieskończoność
\thinlines
\qbezier[8](55,32)(55,35)(55,38)
\qbezier[8](67.5,32)(66.05,35)(64.4,38)
\qbezier[8](69.6,32.6)(67.2,35.7)(64.8,38.4)

%strzałki
%z M_0 na M_1
\qbezier(71.2,-4)(73.5,3.5)(74.2,8.6)\put(74.2,8.6){\vector(0,1){1}}
\qbezier(74.8,-4)(73.5,2)(71.6,7)\put(71.6,7){\vector(-1,3){0.4}}
\qbezier(94,-4)(73,-1)(56.4,8)\put(56.4,8){\vector(-2,1){1}}
%z M_1 na M_2
\put(55,10.2){\vector(0,1){7.5}}
\qbezier(70.8,10)(70.7,16)(70.4,17.7)\put(70.4,17.7){\vector(-1,4){0.1}}
\qbezier(74,11.6)(73.8,15.1)(73,18.9)\put(73,18.9){\vector(-1,4){0.1}}
%z M_2 na M_3
\put(55,20.2){\vector(0,1){4.6}}
%\qbezier(70,20.2)(70.2,20.2)(70.2,22)
\put(70.2,20.2){\vector(-1,4){1.2}}
\put(72.5,21.2){\vector(-1,4){1.1}}
%z M_3 na M_4
\put(55,27.2){\vector(0,1){2.7}}
\put(69,26.8){\vector(-1,4){0.8}}
\put(71,27.4){\vector(-1,4){0.8}}
%z M_infty na M_infty
\put(54,41){\oval(3,3)[l]} \put(54,41){\oval(3,3)[tr]}
\put(55.8,42.2){\vector(-1,-4){0.6}}
\put(65,41){\oval(3,3)[r]}
\put(65,41){\oval(3,3)[tl]}
\put(63.2,42.2){\vector(1,-4){0.6}}
\qbezier(58,40)(59.5,43)(61,40)\put(61.2,39.8){\vector(2,-3){0.1}}
%podpisy
\small
\put(44,-5.5){$M_0$}
\put(44,8.5){$M_1$}
\put(44,18.5){$M_2$}
\put(44,25.5){$M_3$}
\put(44,30.5){$M_4$}
\put(44,38.5){$M_\infty$}

\end{picture}\end{center}
\caption{Dynamics of  $(\M_\lambda,\tal_\lambda)$ where   $0<\lambda \leq\frac{1}{4}$ (a); $\,$  $\frac{1}{4}<\lambda \leq\frac{1}{2}$ (b).\label{sprezynka <1/2}}
 \end{figure}
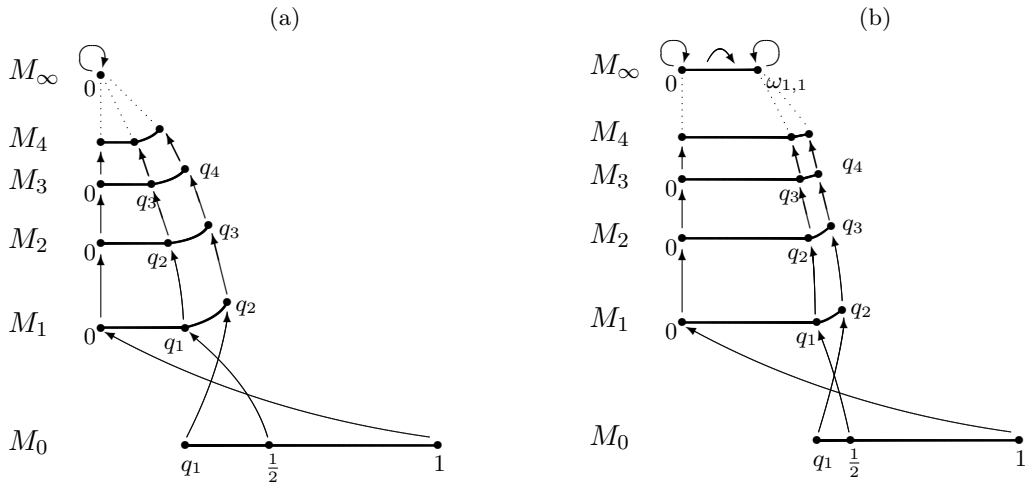
 \par
%  We recall that  the sequence  $\lambda_0=\frac{1}{4}$, $\lambda_1=\frac{3}{4}$, $\lambda_2$, $\lambda_3$, ...,   converges to the Feigenbaum limit $\lambda_\infty$.
We start with parameter   $\lambda$ taking values slightly greater than zero (smaller than $\frac{1}{4}$), see Fig. \ref{lambda < 1/2} (a). The only periodic point of the mapping  $\al_\lambda$ is a stable fixed point  $0$. In particular we have:
$$
q_1>q_2>q_3>...., \qquad \lim_{n\to \infty} q_n=0.
$$
Accordingly, the space $\M_\lambda$ compose of a singleton $M_\infty$   and a sequence of arcs $M_N$ converging to $M_\infty$,  Fig. \ref{sprezynka <1/2} (a).
While we increase the parameter  $\lambda$ the length of intervals $[0,q_N]$ increases.  Finally, when $\lambda$ surpasses   $\frac{1}{4}$ the fixed point  $0$ looses its stability transferring it onto a newly born fixed point  $\omega_{1,1} >0$:
$$
q_1>q_2>q_3>..., \qquad \lim_{n\to \infty} q_n=\omega_{1,1} >0,
$$
Fig. \ref{lambda < 1/2} (b). The set $M_\infty$ becomes an arc corresponding to the interval $[0,\omega_{1,1}]$, and it grows  as  $\lambda$ increases, Fig. \ref{sprezynka <1/2} (b).
 
   \begin{figure}[htb]
  \begin{center}\setlength{\unitlength}{1mm}
\begin{picture}(120,47)(3,0)
%\put(-8,115){a)}
%lambdy
\scriptsize
\put(20,46){(a)}
   \put(101,46){(b)}
  %pikczury
\put(0,0){\includegraphics[angle=0, scale=0.33]{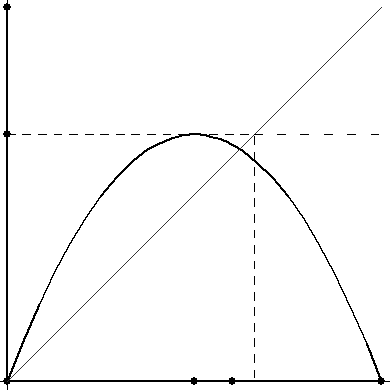}}
\put(80,0){\includegraphics[angle=0, scale=0.33]{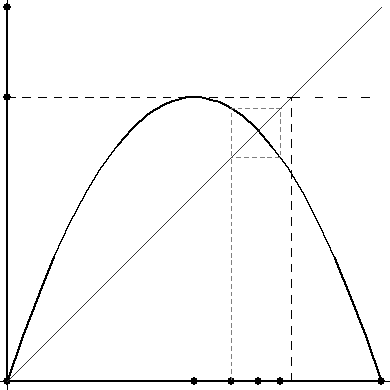}}
%punkty

%zera
\put(-1.8,-1.8){$0$}
\put(78.2,-1.8){$0$}
%jedynki
\put(-2,43.8){$1$}\put(45,-2){$1$}
\put(78,43.8){$1$}\put(124,-2){$1$}
%jedne drugie
%\put(21.6,-2.6){$\frac{1}{2}$}\put(101.6,-2.6){$\frac{1}{2}$}
\put(25.5,-2.2){$\omega_{1,1}$}
%q_1 i q_2
\put(-2.8,29){$\lambda$}
\put(77.1,33){$\lambda$}

\put(104.4,-2){$\omega_{2,2}$}\put(112,-2){$\omega_{2,1}$}
\end{picture}
\end{center}
 \caption{Graph of  $\al_\lambda$  for   $\frac{1}{2}<\lambda \leq\frac{3}{4}$ (a); $\frac{3}{4}<\lambda \leq\frac{1+\sqrt{6}}{4}$ (b).\label{lambda > 1/2}}
 \end{figure}
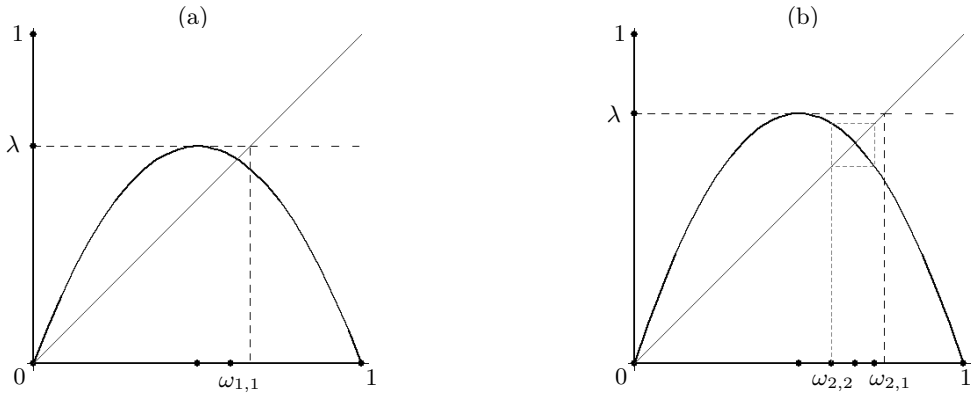
    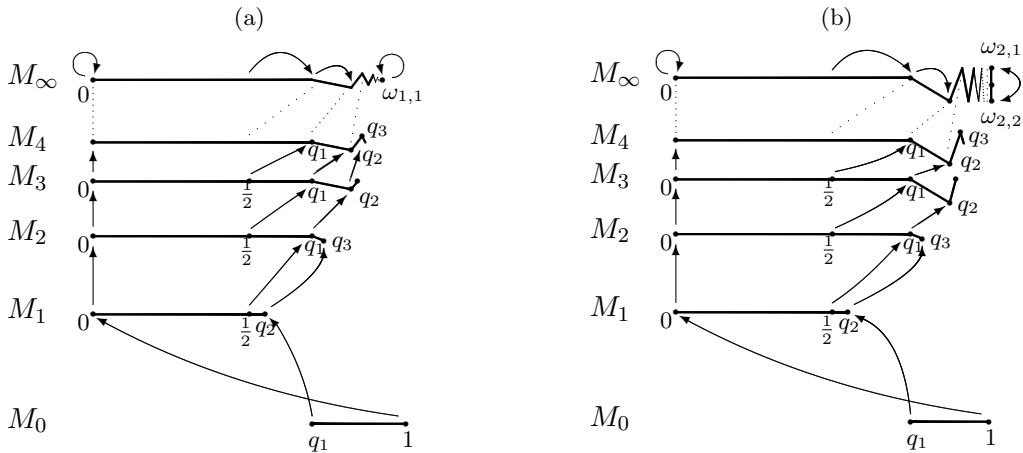
\begin{figure}[htb]
\setlength{\unitlength}{1.04mm}
\begin{center}
\begin{picture}(121,50)(-13,-2)
\thicklines

%mięso
\scriptsize
\put(8,46){(a)}\put(83,46){(b)}
%M_0
\qbezier(18,-5)(18,-5)(30,-5)
%M_1
\qbezier(-10,9)(12,9)(12,9)
%M_2
\qbezier(-10,19)(18,19)(18,19)\qbezier(19.5,18.4)(18.9,18.7)(18,19)
%M_3
\qbezier(-10,26)(18,26)(18,26)\qbezier(23,25)(20.5,25.5)(18,26)\qbezier(23,25)(23.4,25.5)(23.8,26)
%M_4
\qbezier(-10,31)(12,31)(18,31)\qbezier(23,30)(20.5,30.5)(18,31)\qbezier(23,30)(23.7,30.9)(24.4,31.8)
\qbezier(24.4,31.8)(24.6,31.4)(24.8,31)
%M_infty
\qbezier(-10,39)(12,39)(18,39)\qbezier(23,38)(20.5,38.5)(18,39)\qbezier(23,38)(23.7,38.9)(24.4,39.8)
\qbezier(24.4,39.8)(24.8,39.1)(25.2,38.4)
\qbezier(25.8,39.6)(25.5,39)(25.2,38.4)
\thinlines
\qbezier(25.8,39.6)(25.95,39.2)(26.1,38.8)
\qbezier[4](26.38,39.3)(26.24,39)(26.1,38.7)
\qbezier[2](26.38,39.3)(26.45,39)(26.54,38.9)
%punkty
\thinlines
\scriptsize
%na M_0
\put(17.5,-8){$q_1$}
\put(29.5,-8){$1$}
\put(18,-5){\circle*{0.8}}
\put(30,-5){\circle*{0.8}}
%na M_1
\put(-12,7){$0$}
\put(8.5,6){$\frac{1}{2}$}
\put(10.7,6.7){$q_2$}
\put(-10,9){\circle*{0.8}}
\put(10,9){\circle*{0.8}}
\put(12,9){\circle*{0.8}}
%na M_2
\put(-12,17){$0$}

\put(8.6,16.3){$\frac{1}{2}$}
\put(17,17.1){$q_1$}
\put(20.5,17.6){$q_3$}
\put(-10,19){\circle*{0.8}}
\put(10,19){\circle*{0.8}}
\put(18,19){\circle*{0.8}}
\put(19.5,18.4){\circle*{0.8}}

%na M_3
\put(-12,24){$0$}
\put(8.7,23.6){$\frac{1}{2}$}
\put(17.4,23.9){$q_1$}
\put(24,23.2){$q_2$}
\put(-10,26){\circle*{0.8}}
\put(10,26){\circle*{0.8}}
\put(18,26){\circle*{0.8}}
\put(23,25){\circle*{0.8}}
\put(23.8,26){\circle*{0.6}}
%na M_4
\put(-10,31){\circle*{0.8}}
\put(18,31){\circle*{0.8}}

\put(17.4,29.2){$q_1$}
\put(24.4,28.6){$q_2$}
\put(25,32.2){$q_3$}
\put(23,30){\circle*{0.8}}
\put(24.4,31.8){\circle*{0.6}}
%\put(24.8,31){\circle*{0.6}}
%na M_\infty
\put(-12,36.2){$0$}\put(27,36.6){$\omega_{1,1}$}
\put(-10,39){\circle*{0.8}}
\put(27,39){\circle*{0.8}}
%nieskończoność
\thinlines
\qbezier[8](-10,32)(-10,35)(-10,38)
\qbezier[8](10,32)(14,35)(18,38)
\qbezier[8](18,32.4)(20,34.6)(22,37.2)
\qbezier[8](22.8,30)(23.5,34.3)(24.4,38.6)
%strzałki
%z M_0 na M_1
\qbezier(18,-4)(17,2)(13.4,7.8)\put(13.5,7.6){\vector(-1,1){1}}
\qbezier(29,-4)(8,-1)(-8.6,8)\put(-8.6,8){\vector(-2,1){1}}
%z M_1 na M_2
\put(-10,10.2){\vector(0,1){7.5}}
\qbezier(10,10.2)(14,14.8)(16.5,17.8)\put(16.5,17.8){\vector(1,1){0.3}}
\qbezier(12.8,10)(19,14)(19.5,16.7)\put(19.5,16.7){\vector(0,1){0.8}}
%z M_2 na M_3
\put(-10,20.2){\vector(0,1){4.6}}
\qbezier(10.2,20)(14,22.8)(17,24.8)\put(17,24.8){\vector(1,1){0.3}}
\put(18.2,20){\vector(1,1){4.2}}
%z M_3 na M_4
\put(-10,27.2){\vector(0,1){2.7}}
\qbezier(10.2,27)(14,29)(17,30.4)\put(17,30.4){\vector(1,1){0.3}}
\put(18.2,27){\vector(3,2){3.8}}
\put(22.8,26){\vector(1,3){1.3}}
%z M_infty na M_infty
\put(-11,41){\oval(3,3)[l]} \put(-11,41){\oval(3,3)[tr]}
\put(-9.2,42.2){\vector(-1,-4){0.6}}

\put(28.2,40.8){\oval(3,3)[r]}
\put(28.2,40.8){\oval(3,3)[tl]}
\put(26.4,42){\vector(1,-4){0.6}}

\qbezier(10,40)(14,45)(17.9,40)\put(18,39.8){\vector(2,-3){0.1}}
\qbezier(18.5,39.8)(22,42)(22.8,39.4)\put(23,39){\vector(1,-3){0.1}}
%podpisy
\small
\put(-21,-5.5){$M_0$}
\put(-21,8.5){$M_1$}
\put(-21,18.5){$M_2$}
\put(-21,25.5){$M_3$}
\put(-21,30.5){$M_4$}
\put(-21,38.5){$M_\infty$}

\end{picture}
%              Napis
%\begin{picture}(0,0)(108,3) \put(44.5,45){\vector(0,-3){39}} \put(46,41){b}\put(46,37.5){i}\put(46,34){f}\put(46,30.5){u}\put(46,27){r} \put(46,23.5){k}\put(46,20){a}\put(46,16.5){c} \put(46,13){j}\put(46,9.5){a}\end{picture}
\begin{picture}(70,0)(-62,-7)
\thicklines

% druga sprężynka

%mięso
%M_0
\qbezier(20,-5)(20,-5)(30,-5)
%M_1
\qbezier(-10,9)(12,9)(12,9)
%M_2
\qbezier(-10,19)(20,19)(20,19)\qbezier(21.5,18.4)(20.9,18.7)(20,19)
%M_3
\qbezier(-10,26)(20,26)(20,26)\qbezier(25,23)(22.5,24.5)(20,26)\qbezier(25,23)(25.4,24.5)(25.8,26)
%M_4
\qbezier(-10,31)(12,31)(20,31)\qbezier(25,28)(22.5,29.5)(20,31)\qbezier(25,28)(25.7,30)(26.4,32)
\qbezier(26.4,32)(26.6,31.2)(26.8,30.6)
%M_infty
\qbezier(-10,39)(12,39)(20,39)\qbezier(25,36)(22.5,37.5)(20,39)\qbezier(25,36)(25.8,38.1)(26.6,40.2)
\qbezier(26.6,40.2)(27,38.1)(27.4,36)
\qbezier(28.1,40.2)(27.75,38.1)(27.4,36)
\thinlines
\qbezier(28.7,36)(28.3,38.1)(28.1,40.2)
\qbezier(29.1,40.2)(28.9,38.1)(28.7,36)
\qbezier[11](29.4,36)(29.25,38)(29.1,40.2)
\qbezier[5](29.8,36)(29.8,38)(29.8,40)
\thicklines
\qbezier(30.3,36)(30.3,38)(30.3,40.2)
%punkty
\thinlines
\scriptsize
%na M_0
\put(19.5,-8){$q_1$}
\put(29.5,-8){$1$}
\put(20,-5){\circle*{0.8}}
\put(30,-5){\circle*{0.8}}
%na M_1
\put(-12,7){$0$}
\put(8.3,6){$\frac{1}{2}$}
\put(10.8,6.8){$q_2$}
\put(-10,9){\circle*{0.8}}
\put(10,9){\circle*{0.8}}
\put(12,9){\circle*{0.8}}
%na M_2
\put(-12,17){$0$}

\put(8.5,16.2){$\frac{1}{2}$}
\put(19,17.1){$q_1$}
\put(22.5,17.6){$q_3$}
\put(-10,19){\circle*{0.8}}
\put(10,19){\circle*{0.8}}
\put(20,19){\circle*{0.8}}
\put(21.5,18.4){\circle*{0.8}}

%na M_3
\put(-12,24){$0$}

\put(8.5,23.6){$\frac{1}{2}$}
\put(19.4,23.9){$q_1$}
\put(26,22){$q_2$}
\put(-10,26){\circle*{0.8}}
\put(10,26){\circle*{0.8}}
\put(20,26){\circle*{0.8}}
\put(25,23){\circle*{0.8}}
\put(25.8,26){\circle*{0.6}}
%na M_4
\put(-10,31){\circle*{0.8}}
\put(20,31){\circle*{0.8}}

\put(19.4,29.2){$q_1$}
\put(26.4,27.6){$q_2$}
\put(27.4,30.8){$q_3$}
\put(25,28){\circle*{0.8}}
\put(26.4,32){\circle*{0.6}}
%\put(24.8,31){\circle*{0.6}}
%na M_\infty
\put(-12,36.2){$0$}\put(29,33.2){$\omega_{2,2}$}\put(29,42.2){$\omega_{2,1}$}
\put(-10,39){\circle*{0.8}}
\put(20,39){\circle*{0.8}}
\put(25,36){\circle*{0.8}}
\put(30.4,36){\circle*{0.8}}
\put(30.4,38.1){\circle*{0.6}}
\put(30.4,40.2){\circle*{0.8}}
%nieskończoność
\thinlines
\qbezier[8](-10,32)(-10,35)(-10,38)
\qbezier[8](10,32)(14,35)(19,38)
\qbezier[6](20.2,32)(22.2,33.6)(24.2,35.2)
\qbezier[8](24.8,30)(25.5,34.3)(26.4,38.6)
%strzałki
%z M_0 na M_1
\thinlines
\qbezier(20,-4)(19,6)(13.8,8.2)\put(13.8,8.2){\vector(-2,1){1}}
\qbezier(29,-4)(8,-1)(-8.6,8)\put(-8.6,8){\vector(-2,1){1}}
%z M_1 na M_2
\put(-10,10.2){\vector(0,1){7.5}}
\qbezier(10,10.2)(14,12.8)(18.5,17.8)\put(18.5,17.8){\vector(1,1){0.3}}
\qbezier(12.8,10)(21,14)(21.5,16.7)\put(21.5,16.7){\vector(0,1){0.8}}
%z M_2 na M_3
\put(-10,20.2){\vector(0,1){4.6}}
\qbezier(10.4,19.8)(17,22.8)(19.2,25)\put(19.2,25){\vector(1,1){0.3}}
\put(20.2,19.8){\vector(3,2){4.2}}
%z M_3 na M_4
\put(-10,27.2){\vector(0,1){2.7}}
\qbezier(10.4,27)(17,28.4)(19.1,30.2)\put(19.2,30.2){\vector(1,1){0.3}}
\put(20.6,26.6){\vector(3,1){3.7}}
%z M_infty na M_infty
\put(-11,41){\oval(3,3)[l]} \put(-11,41){\oval(3,3)[tr]}
\put(-9.2,42.2){\vector(-1,-4){0.6}}

\qbezier(32,40.2)(36,38.4)(32,36.1)
\put(32,40.1){\vector(-4,1){0.6}}
\put(32,36.18){\vector(-4,-1){0.6}}

\qbezier(10,40)(14,45)(19.8,40)\put(20,39.8){\vector(2,-3){0.1}}

\qbezier(20.8,39.8)(24,41)(24.7,38.4)\put(24.8,38.4){\vector(0,-1){1.3}}
%podpisy
\small
\put(-21,-5.5){$M_0$}
\put(-21,8.5){$M_1$}
\put(-21,18.5){$M_2$}
\put(-21,25.5){$M_3$}
\put(-21,30.5){$M_4$}
\put(-21,38.5){$M_\infty$}

\end{picture}\end{center}
 \caption{Dynamics of $(\M_\lambda,\tal_\lambda)$ where   $\frac{1}{2}< \lambda \leq \frac{3}{4}$ (a); $\frac{3}{4}<\lambda \leq\frac{1+\sqrt{6}}{4}$ (b).\label{sprezynka rysunek2}}
 \end{figure}
 
 When  $\lambda$ reaches the value   $\frac{1}{2}$  all the intervals   $[0,q_n]$, $n\in \N$, $[0,\omega_{1,1}]$ become equal to $[0,\frac{1}{2}]$ and the space  $\M_\lambda$ assumes the shape of a "regular ladder" (in which every step has the same length). As we pass  $\lambda=\frac{1}{2}$ the orbit of the critical  point looses monotonicity:
$$
q_1>q_3>q_5>... > \omega_{1,1}> ...> q_4 >q_2, \qquad \lim_{n\to \infty} q_n=\omega_{1,1},
$$
Fig. \ref{lambda > 1/2} (a). Consequently, each arc $M_N$, $N\in\N\cup\{\infty\}$, develops a "curl" at one of its endpoints: for each $1<N<\infty$ the arc $M_N$  has $N-1$ bendings, whereas $M_{\infty}$ is bended infinitely many times,   Fig. \ref{sprezynka rysunek2} (a).
 \\
\noindent   When $\lambda$ exceeds the value $\lambda_1=\frac{3}{4}$ the first bifurcation occurs -- the fixed point $\omega_{1,1}$    gives  a birth to a new stable periodic orbit
  $\{\omega_{2,1},\omega_{2,2}\}$, Fig. \ref{lambda > 1/2} (b). Then
  $$
  q_1>q_3>... >  \omega_{2,1}>\omega_{1,1}>\omega_{2,2}> ...> q_4 >q_2, \qquad \lim_{n\to \infty} q_{2n+i}=\omega_{2,i},
    $$
   $i=1,2$.   This implies that the attracting fixed point in  $M_\infty$  grows into an arc corresponding to the interval $[\omega_{2,1},\omega_{2,2}]$, see Fig. \ref{sprezynka rysunek2} (b). In other words,  $M_\infty$ becomes a \emph{$\sin(\frac{1}{x})$-continuum}\index{continuum!sin(1/x)}.
   
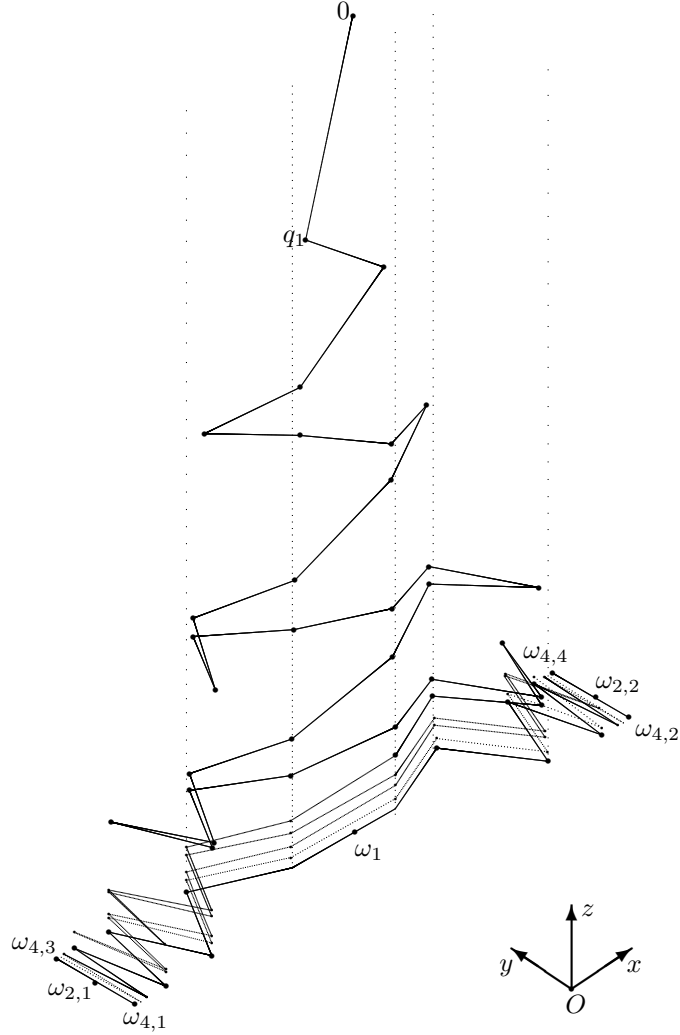
\begin{figure}[hbt]
 \begin{center}
\setlength{\unitlength}{0.8mm}
\begin{picture}(140,168)(-10,3)

\thinlines
%\scriptsize
%\put(-14,170){(a)}
%\put(-14,-15){(b)}
% Sprężynka
\footnotesize
\put(63.7,169.7){\circle*{1}}\put(61,169.2){$0$}
\qbezier(63.6,169.7)(63.6,169.7)(55.8,132.5)
\put(55.8,132.5){\circle*{1}}\put(52.,132.2){$q_1$}
\qbezier(68.8,128)(68.8,128)(55.8,132.5)
\put(68.8,128){\circle*{1}}
%\put(68.4,129){\footnotesize{-$d_0$}}
\qbezier(68.8,128)(68.8,128)(54.9,108)
\put(54.9,108){\circle*{1}}
%\put(50,109.5){\footnotesize{-$2d_0$}}
\qbezier(39,100.3)(39,100.3)(54.9,108)
\put(39,100.3){\circle*{1}}
%\put(32,102){\footnotesize{-$2d_0$-$d_1$}}
\qbezier(39,100.3)(39,100.3)(54.9,100)
\put(54.9,100){\circle*{1}}
%\put(51,101.6){\scriptsize{-$2d_0$-$2d_1$}}
\qbezier(70.05,98.6)(70.05,98.6)(54.9,100)
\put(70.05,98.6){\circle*{1}}
%\put(65,99.8){\scriptsize{-$3d_0$-$2d_1$}}
\qbezier(70.05,98.6)(70.05,98.6)(75.9,105)
\put(75.9,105){\circle*{1}}
%\put(70,106){\scriptsize{-$3d_0$-$2d_1$-$d_2$}}
\qbezier(70,92.6)(70,92.6)(75.9,105)
\put(70,92.6){\circle*{1}}
%\put(63,90){\scriptsize{-$3d_0$-$2d_1$-$2d_2$}}
\qbezier(70,92.6)(70,92.6)(54,75.9)
\put(54,75.9){\circle*{1}}
%\put(48,77){\scriptsize{-$4d_0$-$2d_1$-$2d_2$}}
\qbezier(37.1,69.6)(37.1,69.6)(54,75.9)
\put(37.1,69.6){\circle*{1}}
\qbezier(37.1,69.6)(37.1,69.6)(40.8,57.7)
\put(40.8,57.7){\circle*{1}}
\qbezier(37.1,66.54)(37.1,66.54)(40.8,57.7)
\put(37.1,66.54){\circle*{1}}
\qbezier(37.1,66.54)(37.1,66.54)(53.85,67.7)
\put(53.85,67.7){\circle*{1}}
\qbezier(70.2,71.2)(70.2,71.2)(53.85,67.7)
\put(70.2,71.2){\circle*{1}}
\qbezier(70.2,71.2)(70.2,71.2)(76.3,78.2)
\put(76.3,78.2){\circle*{1}}
\qbezier(94.55,74.7)(94.55,74.7)(76.3,78.2)
\put(94.55,74.7){\circle*{1}}
\qbezier(94.55,74.7)(94.55,74.7)(76.35,75.28)
\put(76.35,75.28){\circle*{1}}
\qbezier(70.25,63.2)(70.25,63.2)(76.35,75.28)
\put(70.25,63.2){\circle*{1}}
\qbezier(70.25,63.2)(70.25,63.2)(53.5,49.5)
\put(53.5,49.5){\circle*{1}}
\qbezier(36.5,43.75)(36.5,43.75)(53.5,49.5)
\put(36.5,43.75){\circle*{1}}
\qbezier(36.5,43.75)(36.5,43.75)(40.6,32.28)
\put(40.6,32.28){\circle*{1}}
\qbezier(23.48,35.7)(23.48,35.7)(40.6,32.28)
\put(23.48,35.7){\circle*{0.95}}
\qbezier(23.48,35.7)(23.48,35.7)(40.3,31.45)
\put(40.3,31.45){\circle*{0.9}}
\qbezier(36.5,40.95)(36.5,40.95)(40.3,31.45)
\put(36.5,40.95){\circle*{0.875}}
\qbezier(36.5,40.95)(36.5,40.95)(53.35,43.43)
\put(53.35,43.43){\circle*{0.85}}
\qbezier(70.74,51.48)(70.74,51.48)(53.35,43.43)
\put(70.74,51.48){\circle*{0.825}}
\qbezier(70.74,51.48)(70.74,51.48)(76.8,59.45)
\put(76.8,59.45){\circle*{0.8}}
\qbezier(95,56.5)(95,56.5)(76.8,59.45)
\put(95,56.5){\circle*{0.775}}
\thinlines
\qbezier(95,56.5)(95,56.5)(88.5,65.4)
\put(88.5,65.4){\circle*{0.67}}
\qbezier(95.05,55.1)(95.05,55.1)(88.5,65.4)
\put(95.05,55.1){\circle*{0.725}}
\qbezier(95.05,55.1)(95.05,55.1)(76.9,56.65)
\put(76.9,56.65){\circle*{0.7}}
\qbezier(70.74,46.78)(70.74,46.78)(76.9,56.65)
\put(70.74,46.78){\circle*{0.675}}
\qbezier[110](70.74,46.78)(62.045,41.355)(53.35,35.93)
\put(53.35,35.9){\circle*{0.65}}
\qbezier[90](36,31.5)(44.675,33.7)(53.35,35.9)
\put(36,31.5){\circle*{0.625}}
\qbezier[90](36,31.5)(38.1,26.25)(40.2,21)
\put(40.2,21){\circle*{0.6}}
\qbezier[90](40.2,21)(31.65,22.7)(23.1,24.4)
\put(23.1,24.4){\circle*{0.575}}
\qbezier[85](32.4,15.4)(27.75,19.9)(23.1,24.4)
\put(32.4,15.4){\circle*{0.575}}
\qbezier[81](32.4,15.4)(27.75,19.7)(23.2,24)
\put(23.1,24){\circle*{0.53}}
\qbezier[80](40.2,20)(31.65,22)(23.1,24)
\put(40.2,20){\circle*{0.5}}
\qbezier[78](36,30.2)(38.1,25.1)(40.2,20)
\put(36,30.2){\circle*{0.5}}
\qbezier[75](36,30.2)(44.675,32)(53.35,33.8)
\put(53.35,33.8){\circle*{0.5}}
\qbezier[90](70.74,43.6)(62.045,38.7)(53.35,33.8)
\put(70.74,43.6){\circle*{0.475}}
\qbezier[90](70.74,43.6)(70.74,43.6)(77.2,53)
\put(77.2,53){\circle*{0.45}}
\qbezier[70](95.6,50.8)(86.4,51.9)(77.2,53)
\put(95.6,50.8){\circle*{0.45}}
\qbezier[69](95.6,50.8)(92.3,55.6)(89,60.4)
\put(89,60.4){\circle*{0.4}}
\qbezier[68](104.6,54.6)(96.8,57.5)(89,60.4)
\put(104.6,54.6){\circle*{0.4}}
\qbezier[67](104.6,54.6)(96.8,57.3)(89,60)
\put(89,60){\circle*{0.4}}
\qbezier[69](95.6,49.8)(92.2,54.9)(88.8,60)
\put(95.6,49.8){\circle*{0.4}}
\qbezier[70](95.6,49.8)(86.4,50.8)(77.2,51.8)
\put(77.2,51.8){\circle*{0.4}}
\qbezier[68](70.74,41.8)(73.97,46.8)(77.2,51.8)
\put(70.74,41.8){\circle*{0.4}}
\qbezier[80](70.74,41.8)(62.045,36.7)(53.35,31.6)
\put(53.35,31.6){\circle*{0.4}}
\qbezier[63](36,27.4)(44.675,29.5)(53.35,31.6)
\put(36,27.4){\circle*{0.4}}
\qbezier[60](36,27.4)(38.1,21.6)(40.2,16.8)
\put(40.2,16.8){\circle*{0.4}}
\qbezier[60](40.2,16.8)(31.65,18.6)(23.1,20.4)
\put(23.1,20.4){\circle*{0.4}}
\qbezier[58](32.6,11.2)(27.85,15.8)(23.1,20.4)
\put(32.6,11.2){\circle*{0.4}}
\qbezier[56](32.6,11.2)(25,14.3)(17.4,17.4)
\put(17.4,17.4){\circle*{0.375}}
\qbezier[55](32.6,10.8)(25,14.05)(17.4,17.3)
\put(32.6,10.8){\circle*{0.375}}
\qbezier[53](32.6,10.8)(27.85,15.25)(23.1,19.7)
\put(23.1,19.8){\circle*{0.35}}
\qbezier[50](40.2,15.6)(31.65,17.7)(23.1,19.8)
\put(40.2,15.6){\circle*{0.325}}
\qbezier[48](36,26)(38.1,20.8)(40.2,15.6)
\put(36,26){\circle*{0.3}}
\qbezier[46](36,26)(44.675,27.9)(53.35,29.8)
\put(53.35,29.8){\circle*{0.3}}
\qbezier[44](70.74,39.6)(62.045,34.7)(53.35,29.8)
\put(70.74,39.6){\circle*{0.3}}
\qbezier[40](70.74,39.6)(74.17,44.6)(77.6,49.6)
\put(77.6,49.6){\circle*{0.3}}
\qbezier[37](96,47.4)(86.4,48.5)(77.6,49.6)
\put(96,47.4){\circle*{0.3}}
\qbezier[32](96,47.4)(92.7,52.2)(89.4,57)
\put(89.4,57){\circle*{0.3}}
\qbezier[28](105,51.4)(97.2,54.5)(89.4,57)
\put(105,51.4){\circle*{0.27}}
\qbezier[23](105,51.4)(99.4,55.6)(93.8,59.8)
\put(93.8,59.8){\circle*{0.25}}

% Podstawa Sprężynki
\thinlines

\put(14.4,13){\circle*{1}}
\put(7,14.6){$\omega_{4,3}$}
\qbezier(14.4,13)(14.4,13)(27.4,5.4)
\put(20.8,9){\circle*{1}}\put(13,7){$\omega_{2,1}$}
\put(27.4,5.4){\circle*{1}}
\put(25.5,2.1){$\omega_{4,1}$}

\qbezier[40](15.6,13.7)(22,9.75)(28.4,5.8)
\qbezier[30](15.2,13.2)(21.8,9.5)(28.4,5.8)

\put(15.7,13.7){\circle*{0.5}}
\qbezier(29.4,6.6)(29.4,6.6)(15.7,13.7)
\put(29.4,6.6){\circle*{0.6}}
\qbezier(29.4,6.6)(29.4,6.6)(17.4,14.7)
\put(17.4,14.7){\circle*{0.7}}
\qbezier(32.6,8.4)(32.6,8.4)(17.4,14.7)
\put(32.6,8.4){\circle*{0.8}}
\qbezier(32.6,8.4)(32.6,8.4)(23.1,17.4)
\put(23.1,17.4){\circle*{0.9}}
\qbezier(40.2,13.5)(40.2,13.5)(23.1,17.4)
\put(40.2,13.5){\circle*{1}}
\qbezier(40.2,13.5)(40.2,13.5)(36,24)
\put(36,24){\circle*{1}}
\qbezier(53.6,28.1)(53.6,28.1)(36,24)
\put(64,34){\circle*{1}}
\put(64,30){$\omega_1$}
\qbezier(53.6,28.1)(53.6,28.1)(70.7,37.9)
%\put(70.7,37.9){\circle*{1}}
%\put(68.2,33.8){\footnotesize{$(\rho,\infty)$}}
\qbezier(77.6,48)(77.6,48)(70.7,37.9)
\put(77.6,48){\circle*{1}}
%\put(74.6,43){$\omega_{4,1}$}
\qbezier(77.6,48)(77.6,48)(96.1,45.8)
\put(96.1,45.8){\circle*{1}}
\qbezier(89.4,55.7)(89.4,55.7)(96.1,45.8)
\put(89.4,55.7){\circle*{0.9}}
\qbezier(89.4,55.7)(89.4,55.7)(105,50.2)
\put(105,50.2){\circle*{0.8}}
\qbezier(93.8,58.7)(93.8,58.7)(105,50.2)
\put(93.8,58.7){\circle*{0.7}}
\qbezier(93.8,58.7)(93.8,58.7)(107.7,51.8)
\put(107.7,51.8){\circle*{0.6}}
\qbezier(95.5,59.8)(95.5,59.8)(107.7,51.8)
\put(95.5,59.8){\circle*{0.5}}

\qbezier[82](77,48)(77,92)(77,170)
\qbezier[120](53.6,28.1)(53.6,78.1)(53.6,158.1)
\qbezier[102](70.7,37)(70.7,87)(70.7,167)
\qbezier[48](36,24)(36,64)(36,154)
\qbezier[42](96.1,45.8)(96.1,70.8)(96.1,160.8)

\qbezier[35](95.5,59.8)(102.05,56)(108.6,52.2)
\qbezier[25](96.1,60)(102.35,56.1)(108.6,52.2)

\put(109.5,53.1){\circle*{1}}\put(110.6,50){$\omega_{4,2}$}
\qbezier(96.8,60.5)(96.8,60.5)(109.5,53.1)\put(104,56.5){\circle*{1}}\put(104,58){$\omega_{2,2}$}
\put(96.8,60.5){\circle*{1}}\put(92,63){$\omega_{4,4}$}

% UCS
\thicklines
\put(100,8){\circle*{0.7}}
\put(100,8){\vector(3,2){10}}\put(109.2,10.8){$x$}
\put(100,8){\vector(-3,2){10}}\put(88,11){$y$}
\put(100,8){\vector(0,1){14}}\put(101.6,20){$z$}
\put(99,4){\footnotesize{$O$}}

\end{picture}
\end{center}
  \caption{Snake-like continuum   $M_\infty\subset \M_\lambda$ after the second bifurcation (immersed into the $3$-dimensional space $Oxyz$).}\label{sprezynka rysunek4}
 \end{figure}

   For $\lambda$ lying approximately in the middle of the interval  $(\frac{3}{4}, \frac{1 +\sqrt{6}}{4}]=(\lambda_1,\lambda_2]$ the orbit  $\{\omega_{2,1},\omega_{2,2}\}$ is \emph{superstable} (it coincides with the critical point orbit). Lengths of the "adjoint" intervals get equal: $
  |q_{k+1}- q_k|=\omega_{2,1} -\omega_{2,1}$, $k>0$, and the  space  $\M_\lambda$ assumes  a regular shape.
  Afterwards, as we pass $\lambda=\frac{1}{2}$ the orbit of the critical point converges to the stable orbit in a more complicated manner:
  $$
  q_1>q_5>... > \omega_{2,1} >...>q_7>q_3\quad \quad q_4>q_8 >... >\omega_{2,2}> ...> q_6 >q_2.
  $$
At the level of the subspace $M_\infty\subset \M_\lambda$, see Fig. \ref{sprezynka rysunek3}, it causes  perturbations around the periodic points; the arcs of the limit bar begin to curl around their endpoints.
    \begin{figure}[htb]
\begin{center}
\setlength{\unitlength}{0.85mm}
\begin{picture}(120,65)(-25,-23)
\thinlines \scriptsize
\qbezier(0,3)(0,3)(-40,3)\put(-40,3){\circle*{1}}\put(-40,-0.5){$0$}
\qbezier(0,3)(6,10.5)(12,18)\put(0,3){\circle*{1}}\put(0,-0.5){$q_1$}
\qbezier(24,2)(18,10)(12,18)\put(12,18){\circle*{1}}
\put(13,19){$q_2$}
\qbezier(24,2)(24,2)(28,-7)\put(24,2){\circle*{1}}
\put(20,1){$q_1$}
\qbezier(32,2)(32,2)(28,-7)\put(28,-7){\circle*{1}}
\put(29,-9.5){$q_3$}
\qbezier(32,2)(32,2)(42,18.5)\put(32,2){\circle*{1}}
\put(33.5,1){$q_1$}
\qbezier(45,28)(45,28)(42,18.5)\put(42,18.5){\circle*{1}}
\put(37,19){$q_2$}
\qbezier(45,28)(45,28)(48,18.5)\put(45,28){\circle*{1}}
\put(46,29){$q_4$}
\qbezier(56,1)(56,1)(48,18.5)\put(48,18.5){\circle*{1}}
\put(50,19){$q_2$}
\qbezier(59,-8.5)(59,-8.5)(56,1)\put(56,1){\circle*{1}}
\put(51.5,1){$q_1$}
\qbezier(59,-8.5)(59,-8.5)(60.5,-14)\put(59,-8.5){\circle*{1}}
\put(54,-9.5){$q_3$}
\qbezier(62,-8.5)(62,-8.5)(60.5,-14)\put(60.5,-14){\circle*{1}}
\put(60,-17){$q_5$}
\qbezier(62,-8.5)(62,-8.5)(64,1.5)\put(62,-8.5){\circle*{1}}
\put(64,-9.5){$q_3$}
\qbezier(70,19.5)(70,19.5)(64,1.5)\put(64,1.5){\circle*{1}}
\put(65.5,1){$q_1$}
\qbezier(70,19.5)(70,19.5)(72,29)\put(70,19.5){\circle*{1}}
\put(66,19){$q_2$}\put(68,29){$q_4$}\put(74,36){$q_6$}
\qbezier(73.5,35)(73.5,35)(72,29)\put(72,29){\circle*{1}}

\put(103,20){\circle*{1}}\put(103,30){\circle*{1}}\put(103,36){\circle*{1}}
\put(103,0){\circle*{1}}\put(103,-10){\circle*{1}}\put(103,-16){\circle*{1}}

%\qbezier(-74.5,38.5)(-74.5,38.5)(-73.5,35)\put(-73.5,35){\circle*{1}}
%\qbezier(-74.5,38.5)(-74.5,38.5)(-75.5,35)\put(-75.5,35){\circle*{1}}
\qbezier(74.5,29)(74.5,29)(73.5,35)\put(73.6,35){\circle*{1}}
\qbezier(74.5,29)(74.5,29)(76,19)\put(74.5,29){\circle*{1}}
\qbezier(80,0.5)(80,0.5)(76,19.5)\put(76,19.5){\circle*{1}}
\qbezier(80,0.5)(80,0.5)(81.5,-9.5)\put(80,0.5){\circle*{1}}
\qbezier(82.5,-15.5)(82.5,-15.5)(81.5,-9.5)\put(81.5,-9.5){\circle*{1}}
\qbezier(82.5,-15.5)(82.5,-15.5)(83.25,-19)\put(82.5,-15.5){\circle*{1}}
\qbezier(84,-15.5)(84,-15.5)(83.25,-19)\put(83.25,-19){\circle*{1}}
\qbezier(84,-15.5)(84,-15.5)(85,-9.5)
\qbezier(86,0.5)(86,0.5)(85,-9.5)
\qbezier(86,0.5)(86,0.5)(88,20)
\qbezier(89,30)(89,30)(88,20)
\qbezier(89,30)(89,30)(89.75,36)
\qbezier(90.25,40)(90.25,40)(89.75,36)
\qbezier(90.25,40)(90.25,40)(91,36)
\qbezier(94,-22)(92.5,8)(91,36)
\qbezier(94,-22)(95.75,10)(96.5,42)
\qbezier[220](98.5,-23)(97.5,9.5)(96.5,42)
\qbezier[200](98.5,-23)(99.25,10)(100,43)
\qbezier[190](101,-24)(100.5,9.5)(100,43)
\qbezier[170](101,-24)(101,10)(102,44)
\qbezier[140](102,-24.5)(102,10)(102,44)

\qbezier[50](-40,-19)(50,-19)(103,-19)
\qbezier[100](-40,-16)(50,-16)(103,-16)
\qbezier[150](-40,-10)(50,-10)(103,-10)
\qbezier[200](-40,0)(50,0)(103,0)
\qbezier[200](-40,20)(50,20)(103,20)
\qbezier[150](-40,30)(50,30)(103,30)
\qbezier[100](-40,36)(50,36)(103,36)
\qbezier[50](-40,40)(50,40)(103,40)

\put(105.5,0){$q_1$}
\put(105.5,20){$q_2$}
\put(105.5,-10){$q_3$}
\put(105.5,30){$q_4$}
\put(105.5,-16){$q_5$}
\put(105.5,36){$q_6$}

\normalsize
\qbezier(103,-25)(103,10)(103,45)
\put(103,-25){\circle*{1}}\put(105,-26){$\omega_{2,1}$}
\put(103,10){\circle*{1}}\put(105,9){$\omega_{1,1}$}
\put(103,45){\circle*{1}}\put(105,45){$\omega_{2,2}$}

\end{picture}
  \end{center}
   \caption{Snake-like continuum   $M_\infty\subset \M_\lambda$ for $\lambda$ approaching  the second bifurcation.\label{sprezynka rysunek3}}
 \end{figure}
  \\
 When $\lambda$ exceeds the value  $\lambda_2=\frac{1+\sqrt{6}}{4}$ the second bifurcation occurs. The orbit  $\{\omega_{2,1},\omega_{2,2}\}$ becomes repelling and  transfers its stability  onto a newly born periodic orbit  $\{\omega_{4,1},\omega_{4,2},\omega_{4,3},\omega_{4,4}\}$ of period $4$. Consequently, periodic points of  $M_\infty$  grow into arcs corresponding to intervals $[\omega_{4,1},\omega_{4,3}]$ and $[\omega_{4,2},\omega_{4,4}]$:
 \begin{quote}
\emph{The maximal ideal space $\M_\lambda$ of algebra $\B_\lambda$ for  $\lambda\in (\lambda_2,\lambda_3]$ compose of snake-like continuum $M_\infty$ presented on Fig.  \ref{sprezynka rysunek4} and a sequence of arcs $\{M_N\}_{N\in \N}$ converging to $M_\infty$ in Hausdorff metric.}
\end{quote}

   %
 %  N I E B E Z P I E C Z   N I E B E Z P I E C Z   N I E B E Z P I E C Z  N I E B E Z P I E C Z
 %
 This process continues: as  $\lambda$ approaches approximately a middle of the interval  $(\lambda_2,\lambda_3]$ the orbit  $\{\omega_{4,1},\omega_{4,2},\omega_{4,3},\omega_{4,4}\}$ becomes superstable and the space  $M_\infty$ assumes a regular shape. Afterwards continuum  $M_\infty$ develops curls around the points corresponding to the orbit $\{\omega_{4,1},\omega_{4,2},\omega_{4,3},\omega_{4,4}\}$ until we pass   $\lambda=\lambda_3$ where each of these points grows into an arc. And so on, and so forth, cf. Fig. \ref{czwarta bifurkacja doubling}.
    \par
    In order to state the result formally we extend the sequence $\lambda_1=\frac{3}{4}$, $\lambda_2$, $\lambda_3$, ..., putting $\lambda_{-1}=0$ and  $\lambda_0=\frac{1}{4}$.
By a \emph{ray} we mean a topological space homeomorphic to $(0,1]$, and  \emph{an endpoint} of a topological space  $M$ is a point  $p\in M$ whose every neighbourhood  $U$ contains an open neighbourhood $V$  such that the boundary of $V$ is a singleton, cf. \cite{Nadler}. Subspace  $M_\infty=\underleftarrow{\,\,\lim\,}(M,\al_\lambda)\subset \M_\lambda$ for $\lambda <\lambda_\infty$ is given by the following recurrence.
\begin{thm}\emph{\cite[Thm. 3]{Barge-Ingram}}\label{bifurakcyjny opis widma thm}
If  $\lambda \in  (\lambda_n,  \lambda_{n+1}]$, $n\in \N$, then the  continuum  $\underleftarrow{\,\,\lim\,}(M,\al_\lambda)$ is the closure of a ray $R$ such that  $\underleftarrow{\,\,\lim\,}(M,\al_\lambda)\setminus R$  is the union  of two copies of $\underleftarrow{\,\,\lim\,}(M,\al_{\lambda'})$, where  $\lambda'\in (\lambda_{n-1},\lambda_n]$,  intersecting in a common endpoint.
\end{thm}
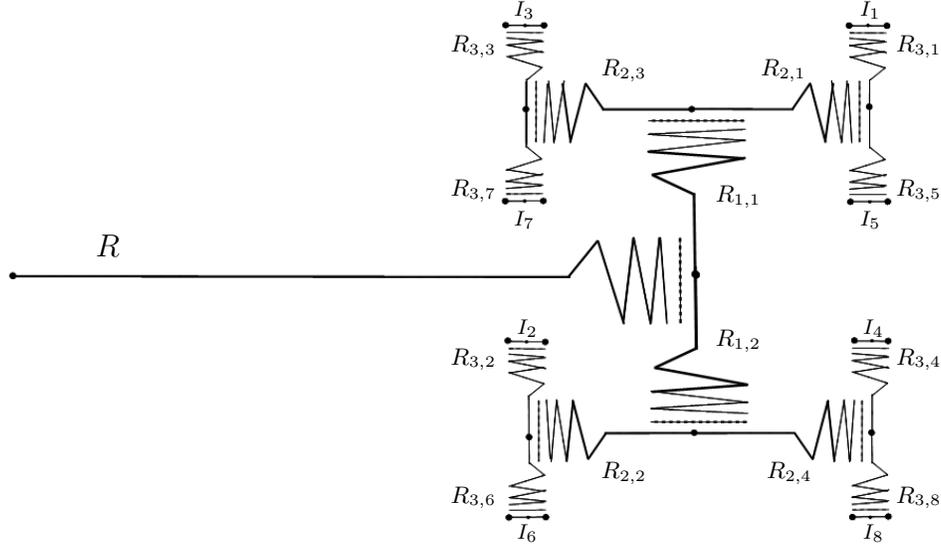
\begin{figure}[hbt]
  \begin{center}\setlength{\unitlength}{0.8mm}
\begin{picture}(130,83)(-10,3)
  \put(-24,0){\includegraphics[angle=0, scale=0.396]{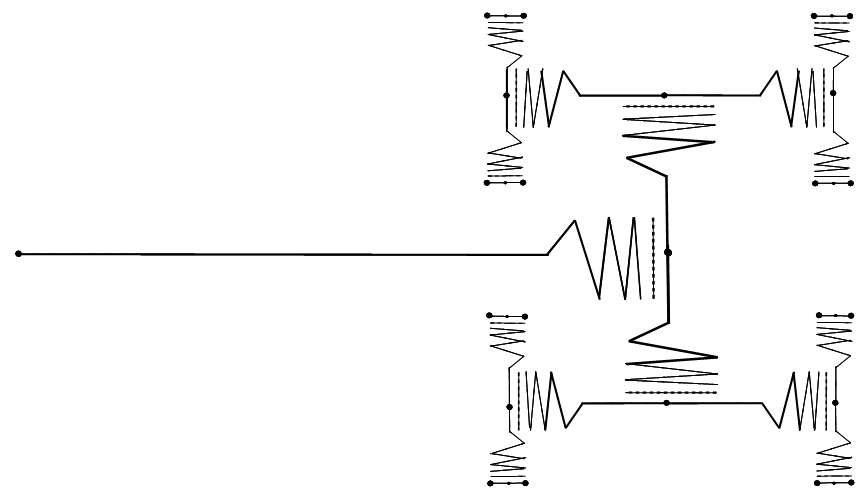}}

 \put(-7,46){$R$}
\scriptsize

\put(96,31){$R_{1,2}$}\put(96,55){$R_{1,1}$}

\put(77,76){$R_{2,3}$}\put(103.5,76){$R_{2,1}$}

\put(77,9){$R_{2,2}$}\put(104.5,9){$R_{2,4}$}

\put(126,80){$R_{3,1}$}\put(126,56){$R_{3,5}$}

\put(52,80){$R_{3,3}$}\put(52,56){$R_{3,7}$}

\put(52,28){$R_{3,2}$}\put(52,5){$R_{3,6}$}

\put(126,5){$R_{3,8}$}\put(126,28){$R_{3,4}$}

\tiny

\put(120,86){$I_{1}$}\put(120,51){$I_{5}$}
\put(62.5,86){$I_{3}$}\put(62.5,51){$I_{7}$}
\put(120.5,33){$I_{4}$}\put(120.5,-1){$I_{8}$}
\put(63,33){$I_{2}$}\put(63,-1){$I_{6}$}
\end{picture}
\end{center}
 \caption{Subspace  $M_\infty\subset \M_\lambda$ after the forth bifurcation (schematic presentation).\label{czwarta bifurkacja doubling}}
 \end{figure}
 
Now we are ready to give a full description of the system   $(\M_\lambda, \tal_\lambda)$  for $\lambda$ lying in the interval $(0,\lambda_\infty)=\bigcup_{n=-1}^\infty(\lambda_n,\lambda_{n+1}]$.  If $\lambda\in (0,\frac{3}{4}]$, such system is presented on Fig.  \ref{sprezynka <1/2}, \ref{sprezynka rysunek2}.
\begin{thm}
\label{bifurakcyjny opis widma thm2}
Let   $U_\lambda$ be the operator given by \eqref{T_lambda operator defn} and  $\B_\lambda$ the $C^*$-algebra given by \eqref{B_lambda algebra defn}. If  $\lambda \in  (\lambda_n,  \lambda_{n+1}]$, $n>0$, then
\begin{itemize}
\item[i)] the maximal ideal space  $\M_\lambda$ of the algebra  $\B_\lambda$ compose of a snake-like  continuum $M_\infty$ and a sequence of arcs $\{M_N\}_{N\in \N}$ converging to $M_\infty$, where
$$
M_\infty=R\cup (R_{1,1}\cup R_{1,2})\cup ... \cup (R_{n-1,1}\cup...\cup R_{n-1,2^{n-1}})\cup (I_1 \cup ... \cup I_{2^{n-1}})
$$
is the sum of  $2^n-1$ rays $R$, $R_{k,i}$, $k=1,...,n-1$, $i=1,...,2^{k}$, and $2^{n-1}$ arcs $I_i$, $i=1,...,2^{n-1}$, cf. Fig. \ref{czwarta bifurkacja doubling}.
The closure of $R$ gives $M_\infty$ and
$$
\overline{R_{k,i}}=\bigcup_{j=0}^{n-k-1}\bigcup_{l=0}^{2^j-1}  R_{k+j,i+ l\cdot 2^k} \cup \bigcup_{l=0}^{2^{n-k-1}-1}  I_{i+ l\cdot 2^k}, \qquad k=1,...,n-1, i=1,...,2^{k}.
$$
\item[ii)]  Partial homeomorphism  $\tal_\lambda$ generated by $U_\lambda$ on $\M_\lambda$ carries  $M_N$ onto  $M_{N+1}$, $N\in \N$, and  $\tal_\lambda:M_\infty \to M_\infty$ is a homeomorphism that preserves $R$, permutes cyclically  arcs $I_i$, and  the rays $R_{k,i}$ (for each fixed  $k=1,...,n-1$):
$$
\tal_\lambda(I_1)=I_2,\,...,\,\tal_\lambda(I_{2^n})=I_1,\qquad  \tal_\lambda(R_{k,1})=R_{k,2},\,...,\,\tal_\lambda(R_{k,2^k})=R_{k,1}.
$$
\end{itemize}
Moreover, all the rays $R$, $R_{k,i}$ and  arcs $I_i$ are pairwise disjoint except of the following intersections
$$
R_{k,i}\cap R_{k,2^{k-1}+i}=\{\tomega_{2^{k-1},i}\}, \qquad  k=1,...n-1,\,\,i=1,...2^{k-1},
$$
that form periodic orbits   $\{\tomega_{2^{k},1}, ... \tomega_{2^{k},2^{k}}\}$ with period $2^k$, $k=1,...,n-2$.
Middle-points of arcs $I_i$ form a periodic orbit  $\{\tomega_{2^{n-1},1}, ... \tomega_{2^{n-1},2^{n-1}}\}$
with period $2^{n-1}$, and  endpoints of arcs $I_i$ form a periodic orbit
$\{\tomega_{2^{n},1}, ... \tomega_{2^{n},2^{n}}\}$
with period $2^{n}$.
\end{thm}
\begin{Proof}
For the description of the space $\M_\lambda$ apply Theorem \ref{sprowadzenie do granicy odwrotnej},   \ref{bifurakcyjny opis widma thm}. The dynamics of $\tal_\lambda$ on $M_\infty$ may deduced from the proof of  \cite[Thm. 3]{Barge-Ingram}, see also \cite[6]{Barge-Ingram}.
\end{Proof}
 
The infinite sequence of period-doubling  bifurcation  leaves the following imprint   on the structure of  $\M_\lambda$ for  $\lambda$ attaining the  Feigenbaum limit.
\begin{thm}
For  $\lambda=\lambda_\infty$ the maximal ideal space  $\M_\lambda$ of the algebra $\B_\lambda$ possess the property that  every nondegenerate subcontinuum of $\M_\lambda$  is reducible. Moreover $M_\lambda$ contains only three topologically different nondegenerate subcontinua: arcs, copies of the space $M_\infty$ and  sums of two copies of $M_\infty$ intersecting in the common endpoint.
\end{thm}
\begin{Proof}
Apply   Theorem \ref{sprowadzenie do granicy odwrotnej} and  \cite[Thm. 7]{Barge-Ingram}.
\end{Proof}
 \subsection{The parameter values corresponding to chaotic dynamics}\label{knaster coninuua}

 For   $\lambda > \lambda_\infty$ the dynamics of the system  $([0,1],\al_\lambda)$ is chaotic and the available  knowledge concerning mappings $\{\al_\lambda\}_{\lambda >\lambda_\infty}$ is far from being complete.
However we are able to present a number of results that shed much light onto the  structure of the considered systems. For example, we already know that   $\M_\lambda$ contains irreducible continua  (Theorem \ref{bifurkacji ciag pierwszy}), and for  $\lambda=1$, $\M_\lambda$  is a  B-J-K continuum.   We start with the case when B-J-K continua are the only irreducible subcontinua of $\M_\lambda$.
  \subsubsection{Cascade of B-J-K continuum doubling}\label{musze miec tu etykietke}
Let us consider the sequence   $\mu_0=1$, $\mu_1$, $\mu_2$, ...,    of parameter values, in which  the bifurcation diagram  splits into two "copies" of itself,  Fig.  \ref{Diagram Bifurkacyjny1}. It is a decreasing sequence converging to $\lambda_\infty$,  and formally $\mu_n$ could be defined as the solution of the following equation
$$
\al_\lambda^{2^n}(\lambda)=(\textrm{the largest fixed point of } \al_\lambda^{2^{n-1}}),
$$
 cf.  \cite{Barge-Ingram}, \cite{Col-Eck}. The sequence  $\{\mu_n\}_{n\in \N}$ admits an inductive procedure surprisingly similar to the one presented in Theorem  \ref{bifurakcyjny opis widma thm}.
\begin{thm}\emph{\cite[Thm. 6]{Barge-Ingram}}\label{twierdzenie Barga-ingrama}
For $n>0$,  $\underleftarrow{\,\,\lim\,}(M,\al_{\mu_n})$ is the closure of a ray $R$ such that  $\underleftarrow{\,\,\lim\,}(M,\al_{\mu_n})\setminus R$   is the union of two copies of the space $\underleftarrow{\,\,\lim\,}(M,\al_{\mu_{n-1}})$ intersecting in a common endpoint.
\end{thm}
 
  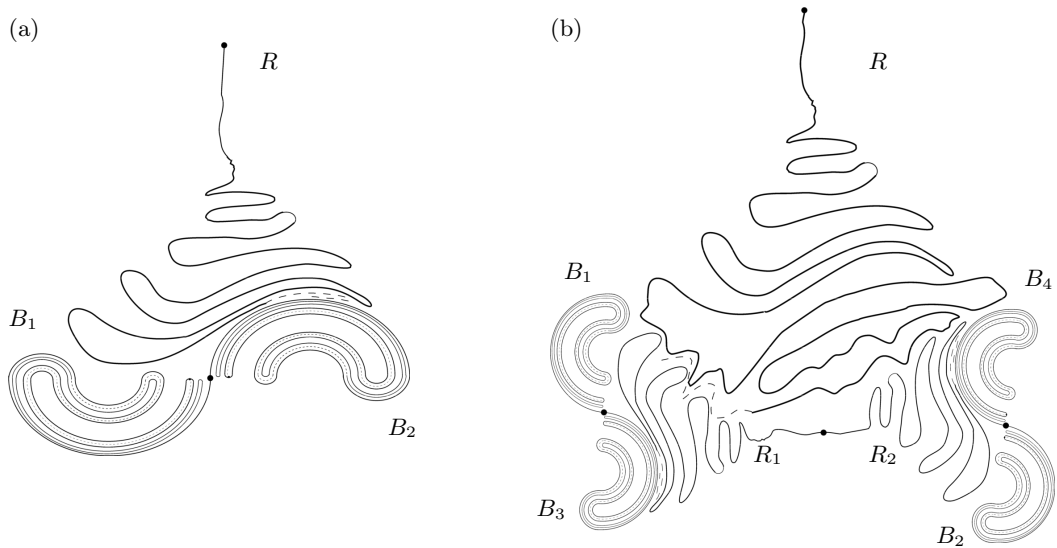
\begin{figure}[htb]
\begin{center}\setlength{\unitlength}{0.9mm}
\begin{picture}(150,75)(2,2)
%\put(-8,115){a)}
%lambdy
  %pikczury

 % \put(80,0){\includegraphics[angle=0, scale=0.17]{bucketx4.png}}
 % \put(0,12){\includegraphics[angle=0, scale=0.18]{bucketx2.png}}
    \put(80,0){\includegraphics[angle=0, scale=0.7]{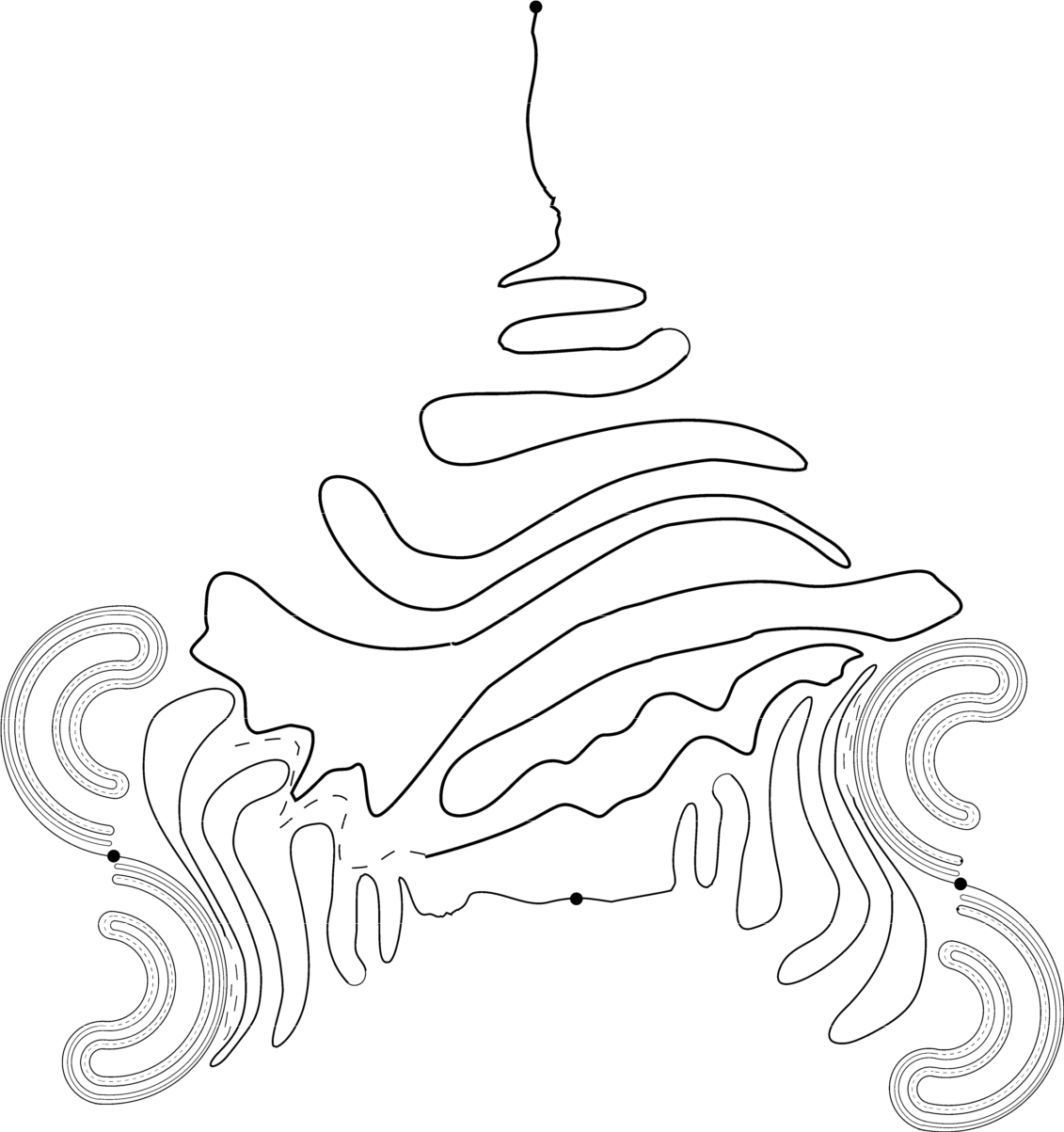}}
  \put(0,13){\includegraphics[angle=0, scale=0.7]{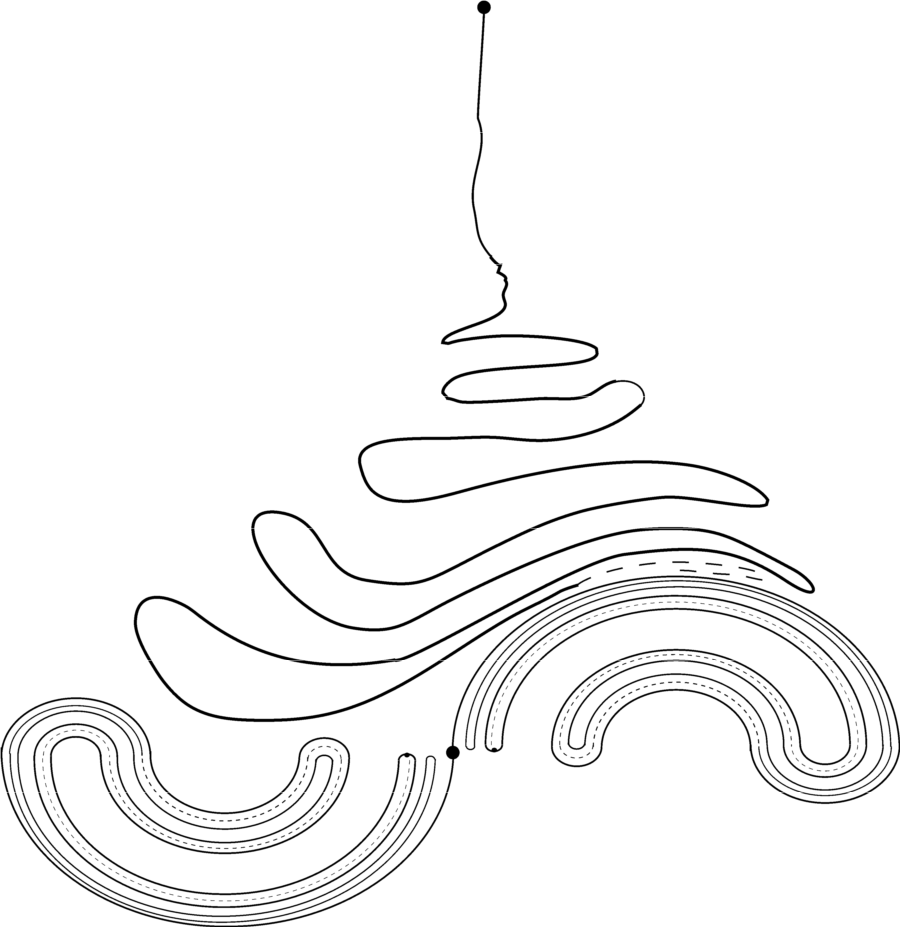}}
 %popisy
\scriptsize
\put(37,70){$R$}\put(127,70){$R$}
\put(80,75){(b)}\put(0,75){(a)}

\put(0,32){$B_1$}\put(56,16){$B_2$}

\put(82,39){$B_1$}\put(78,4){$B_3$}\put(137,0){$B_2$}\put(150,38){$B_4$}

\put(127,12){$R_2$}\put(110,12){$R_1$}
\end{picture}
\end{center}
 \caption{Subspace  $M_\infty\subset \M_\lambda$: for  $\lambda=\mu_1$ (a); for $\lambda=\mu_2$ (b).}\label{continua knastera doubling}
 \end{figure}
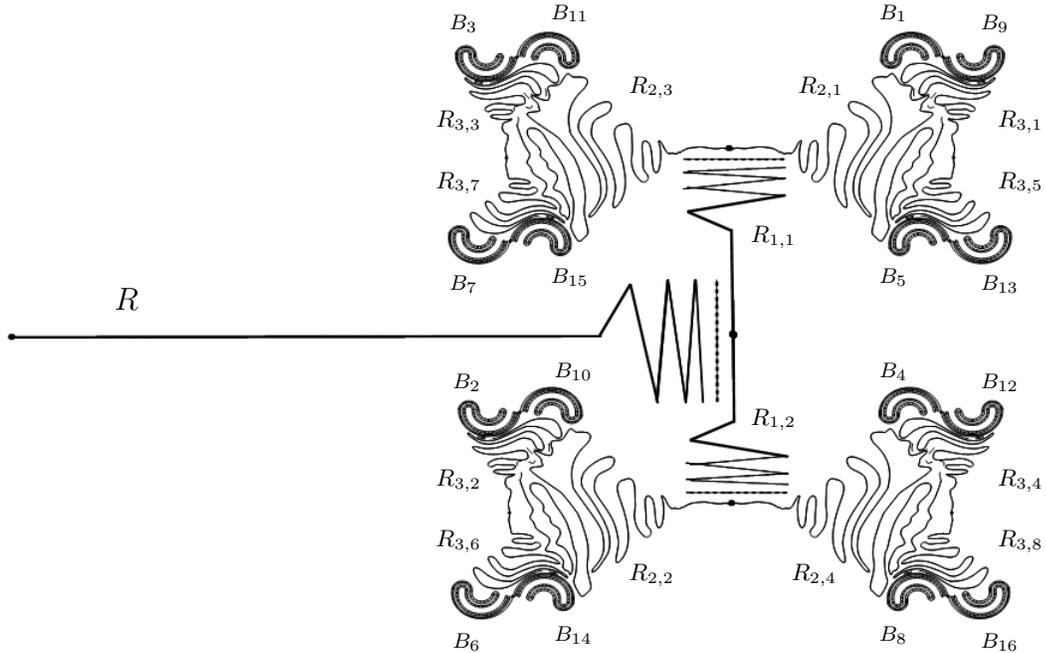
\begin{figure}[hbt]
\begin{center}\setlength{\unitlength}{0.9mm}
\begin{picture}(120,96)(3,-3)
%\put(-8,115){a)}
%lambdy
  %pikczury

  \put(-12,0){\includegraphics[angle=0, scale=0.45]{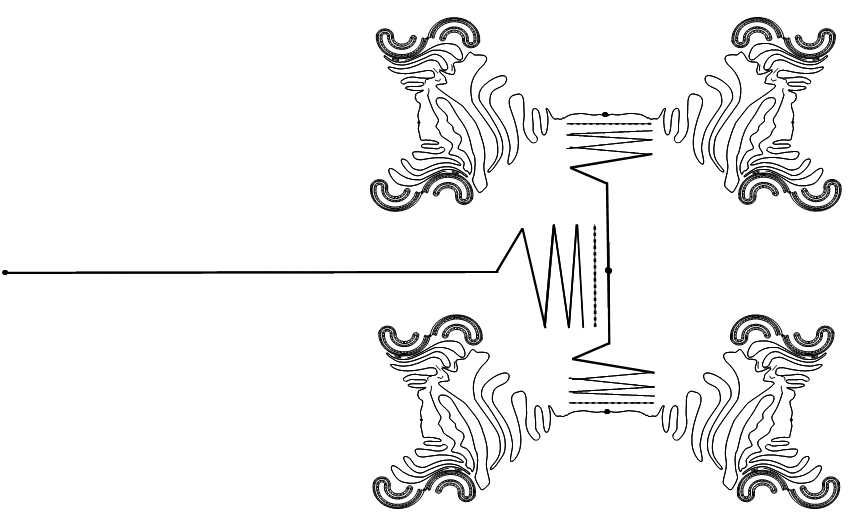}}

%    \put(80,0){\includegraphics[angle=0, scale=0.7]{bucketx4.png}}
%  \put(0,13){\includegraphics[angle=0, scale=0.7]{bucketx2.png}}
 %popisy
 \put(4,48){$R$}
\scriptsize

\put(98,31){$R_{1,2}$}\put(98,58){$R_{1,1}$}

\put(80,80.1){$R_{2,3}$}\put(105,80.1){$R_{2,1}$}

\put(80,8){$R_{2,2}$}\put(104,8){$R_{2,4}$}

\put(134.5,75){$R_{3,1}$}\put(134.5,66){$R_{3,5}$}

\put(51.5,66){$R_{3,7}$}\put(51.5,75){$R_{3,3}$}

\put(51.5,22){$R_{3,2}$}\put(51.5,13){$R_{3,6}$}

\put(134.5,13){$R_{3,8}$}\put(134.5,22){$R_{3,4}$}

\tiny

\put(117,91){$B_{1}$}\put(132,89.4){$B_{9}$}
\put(117,52){$B_{5}$}\put(132,51){$B_{13}$}

\put(53.5,89.4){$B_{3}$}\put(68.5,91){$B_{11}$}
\put(53.5,51){$B_{7}$}\put(68.5,52){$B_{15}$}

\put(117,38){$B_{4}$}\put(132,36.4){$B_{12}$}
\put(117,-1){$B_{8}$}\put(132,-2){$B_{16}$}

\put(54,36.4){$B_{2}$}\put(69,38){$B_{10}$}
\put(54,-2){$B_{6}$}\put(69,-1){$B_{14}$}
%\put(63,98){$I_{3}$}\put(63,61){$I_{7}$}
%\put(118,36){$I_{4}$}\put(118,-1){$I_{8}$}
%\put(63,36){$I_{2}$}\put(63,-1){$I_{6}$}
\end{picture}
\end{center}
 \caption{Subspace $M_\infty\subset \M_\lambda$ for $\lambda=\mu_4$ (schematic presentation).\label{czwarta bifurkacja BJK}}
 \end{figure}

  The above statement together with Theorems \ref{sprowadzenie do granicy odwrotnej}, \ref{takie tam twierdzenie dla lambda=1}  says that for $\lambda=\mu_1$  the space $\M_\lambda$ compose of a sequence of arcs $\{M_N\}_{N\in\N}$ converging to the snake-like continuum   $M_\infty$ which is the union of two copies $B_1$, $B_2$ of   B-J-K continuum  and a ray $R$, see Fig. \ref{continua knastera doubling} (a). Partial homeomorphism  $\tal_\lambda$ transforms  $M_N$ onto $M_{N+1}$, $N\in \N$, and the dynamics of  $\tal_\lambda$ on $M_\infty$ is as follows: the endpoint of the ray $R$ is a fixed point, the remaining points of  $R$  slide on $R$  toward the continua $B_1$, $B_2$.  Points from   $B_1$  are carried onto $B_2$ and vice versa. The intersection  $B_1\cap B_2$ consists of a fixed point.
 \\
Analogously, for $\lambda=\mu_2$, see  Fig. \ref{continua knastera doubling} (b),  $M_\infty\subset \M_\lambda$ compose of four copies  $B_1$, $B_2$, $B_3$, $B_4$ of B-J-K continuum  and three arcs  $R$, $R_1$, $R_2$. The endpoint of  $R$ is a fixed point, and the remaining points  of $R$ move toward arcs $R_1$, $R_2$. Points from $R_1$ are carried onto $R_2$ and  vice versa. The intersection  $R_1\cap R_2$ consists of a fixed point. Continua $B_i$  are cyclically permuted
$$
\tal_\lambda(B_1)=B_2,\quad \tal_\lambda(B_2)=B_3,\quad \tal_\lambda(B_3)=B_4,\quad \tal_\lambda(B_4)=B_1,
$$
and  $(B_1\cap B_3) \cup (B_2\cap B_4)$ constitute a periodic orbit of period  $2$.
\\
In general, for  $\lambda= \mu_n$ we have the description of   $(\M_\lambda,\tal_\lambda)$ which differs from the one presented in  Theorem \ref{bifurakcyjny opis widma thm2} only in that the arcs $I_i$ are replaced with  B-J-K continua, cf. Fig. \ref{czwarta bifurkacja BJK}.

\begin{thm}\label{bifurakcyjny opis knastrow thm2}
Let   $U_\lambda$ be the operator given by \eqref{T_lambda operator defn} and  $\B_\lambda$ the $C^*$-algebra given by \eqref{B_lambda algebra defn}. If  $\lambda =\mu_n$,  $n\in \N_+$, then
\begin{itemize}
\item[i)] the maximal ideal space  $\M_\lambda$ of algebra  $\B_\lambda$ compose of a snake-like  continuum $M_\infty$ and a sequence of arcs $\{M_N\}_{N\in \N}$ converging  to $M_\infty$, where
$$
M_\infty=R\cup (R_{1,1}\cup R_{1,2})\cup ... \cup (R_{n-1,1}\cup...\cup R_{n-1,2^{n-1}})\cup (B_1 \cup ... \cup B_{2^{n}})
$$
is the sum of  $2^n-1$ rays $R$, $R_{k,i}$, $k=1,...,n-1$, $i=1,...,2^{k}$, and $2^{n}$ B-J-K continua  $B_i$, $i=1,...,2^{n}$, cf. Fig. \ref{czwarta bifurkacja BJK}. The closure of $R$ gives $M_\infty$ and
$$
\overline{R_{k,i}}=\bigcup_{j=0}^{n-k-1}\bigcup_{l=0}^{2^{j-1}}  R_{k+j,i+ l\cdot 2^k} \cup \bigcup_{l=0}^{2^{n-k}-1}  B_{i+ l\cdot 2^k}, \qquad k=1,...,n-1, i=1,...,2^{k},
$$

\item[ii)]  Partial homeomorphism  $\tal_\lambda$ generated by $U_\lambda$ on $\M_\lambda$ carries  $M_N$ onto  $M_{N+1}$, $N\in \N$, and on  $\tal_\lambda:M_\infty\to M_\infty$ is a homeomorphism that preserves $R$, permutes cyclically  continua $B_i$, and  the rays $R_{k,i}$ (for fixed  $k=1,...,n-1$):
$$
\tal_\lambda(B_1)=B_2,\,...,\,\tal_\lambda(B_{2^n})=B_1,\qquad  \tal_\lambda(R_{k,1})=R_{k,2},\,...,\,\tal_\lambda(R_{k,2^k})=R_{k,1}.
$$
\end{itemize}
Moreover, all the rays $R$, $R_{k,i}$ and  continua $B_i$ are pairwise disjoint except of the following intersections
$$
R_{k,i}\cap R_{k,2^{k-1}+i}=\{\tomega_{2^{k-1},i}\}, \qquad i=\{1,...2^{k-1}\},
$$
$$
B_i\cap B_{2^{n-1}+i}=\{\tomega_{2^{n-1},i}\} \qquad i=\{1,...2^{n-1}\},
$$
which form periodic orbits   $\{\tomega_{2^{k},1}, ... \tomega_{2^{k},2^{k}}\}$ with period $2^k$, $k=1,...,n-1$.
\end{thm}
\begin{Proof} It follows from Theorems \ref{sprowadzenie do granicy odwrotnej}, \ref{takie tam twierdzenie dla lambda=1} and inductively applied Theorem  \ref{twierdzenie Barga-ingrama}; see also  the proof of \cite[Thm. 6]{Barge-Ingram}.
\end{Proof}
\subsection{Windows of stable periodic orbits of odd period}\label{ustep o nieparzystych okanach}
%Opiszemy teraz zachowanie układu $(\M_\lambda,\tal_\lambda)$ dla wartości parametru $\lambda$
\begin{figure}[hbt]
\begin{center}\setlength{\unitlength}{0.9mm}
\begin{picture}(240,53)(-7,1.5)
\scriptsize
\put(5,-1.5){\includegraphics[angle=0, scale=0.405]{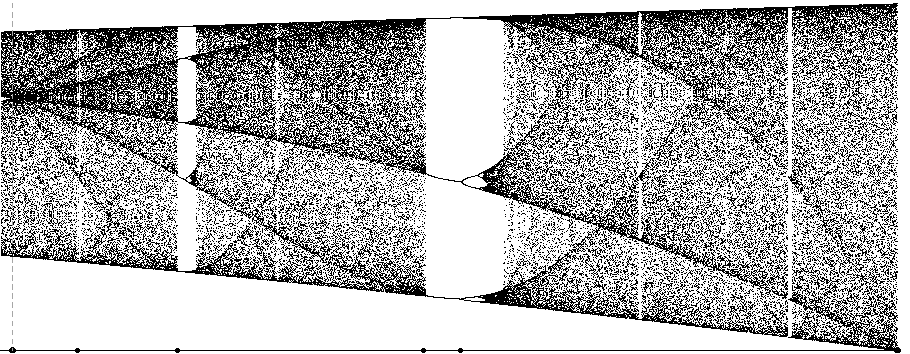}}
\put(6.2,-4){$\mu_1$}
%\put(70,-4){$0.9564$} \put(32,-4){$0.9345$}  \put(16,-4){$0.9253$}  \put(143.5,-4){$1$}
\put(71,-4){$\eta_1$} \put(78.5,-4){$\nu_1$}\put(32,-4){$\eta_2$}  \put(16,-4){$\eta_3$}  \put(148,-4){$1$}
   \end{picture}
\end{center}
\caption{Bifurcation diagram for $\lambda \in [0,91; 1]$. \label{bifurkacja od mu_1}}
 \end{figure}
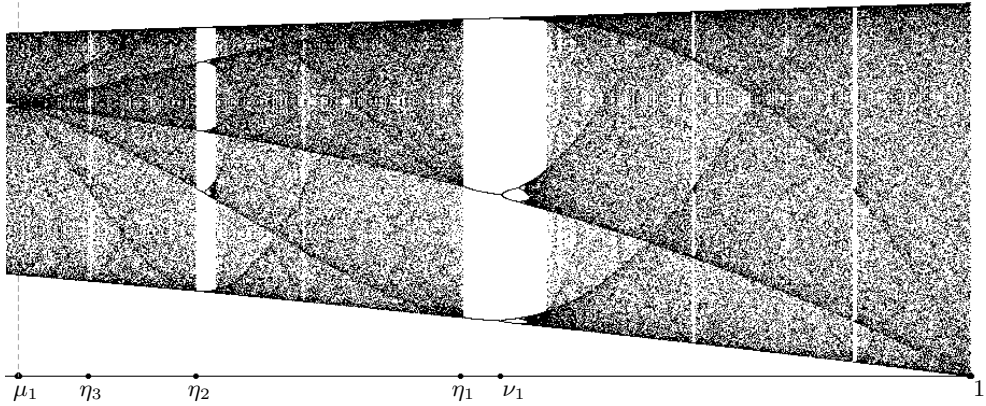
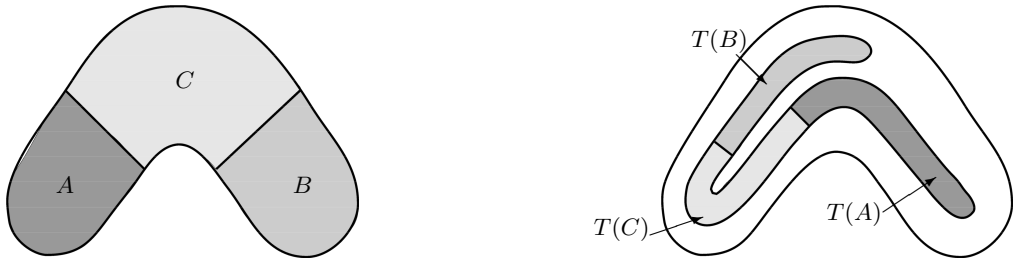
\begin{figure}[htb]
\begin{center}\setlength{\unitlength}{1mm}
\begin{picture}(240,32)(5,2)

\put(13,0){\includegraphics[angle=0, scale=0.42]{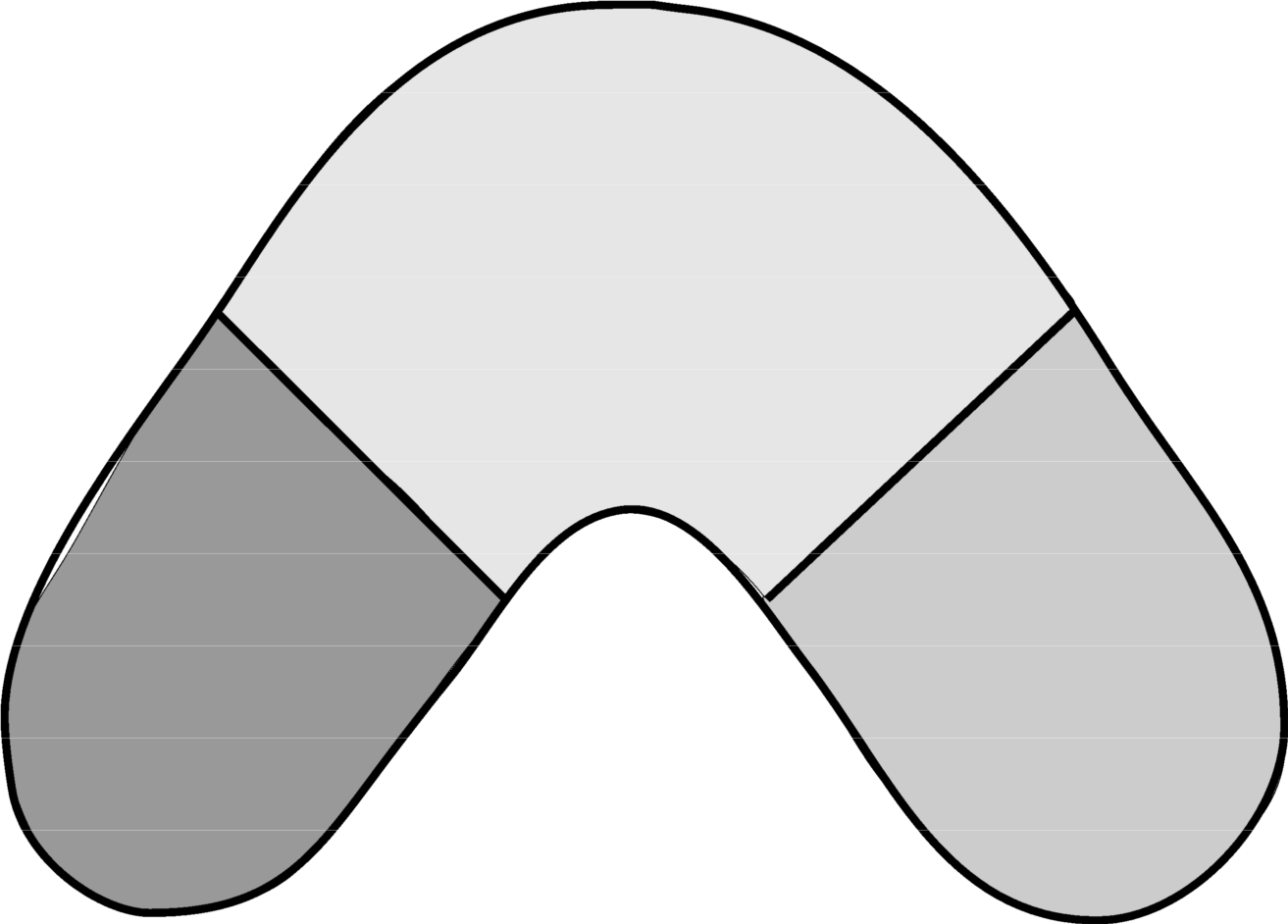}}
  \put(100,0){\includegraphics[angle=0, scale=0.42 ]{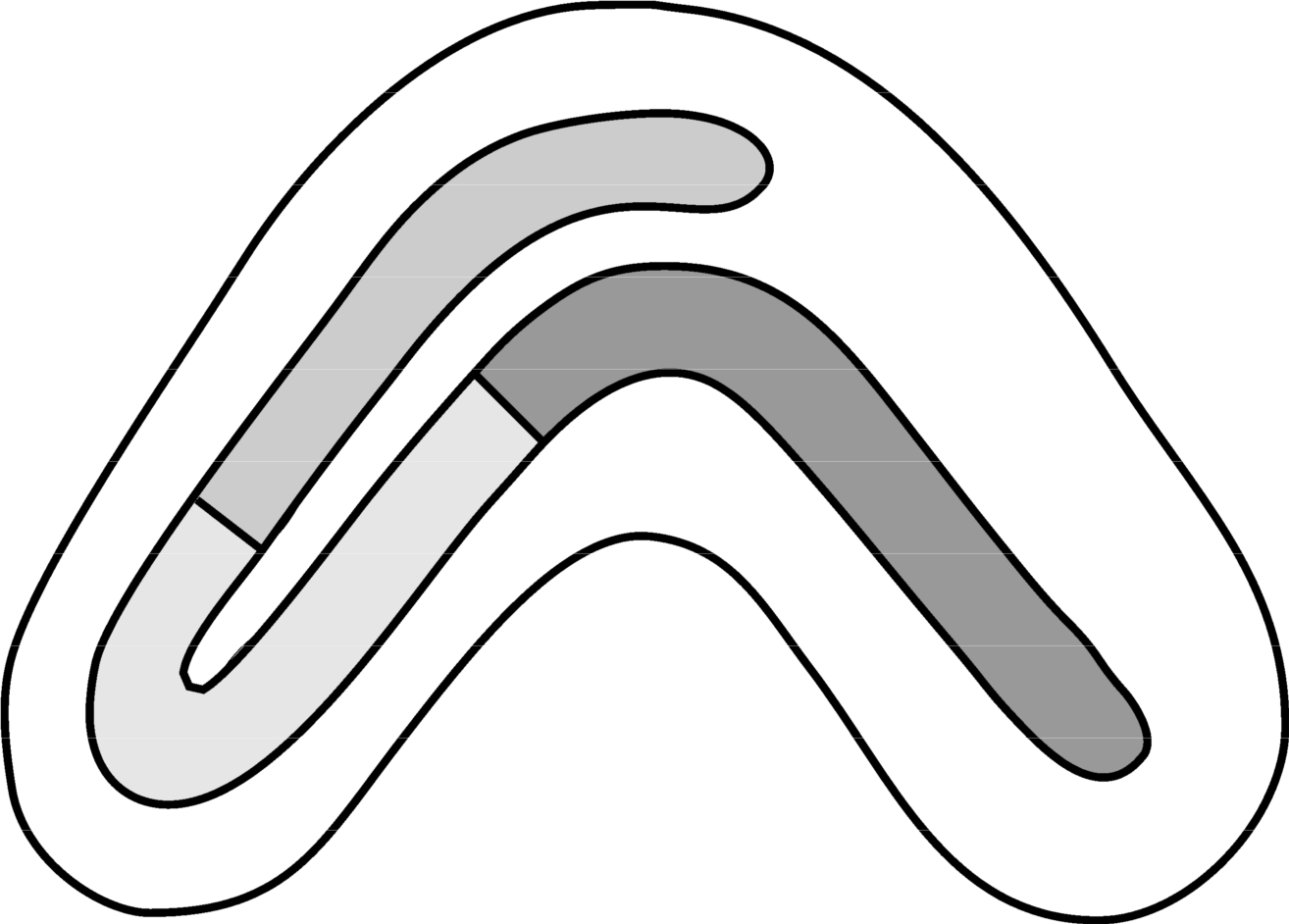}}
%\put(13,0){\includegraphics[angle=0, scale=0.1]{3continuum1.png}}   \put(100,0){\includegraphics[angle=0, scale=0.1]{3continuum2.png}}

  \scriptsize
%\put(35.5,29){c}\put(128,31){(b)}

  \put(19.5,8.5){$A$}
  \put(35.4,22.5){$C$}\put(51,8.5){$B$}

  \put(104,28){$T(B)$}\put(110,27){\vector(1,-1){4}}
  \put(122,5){$T(A)$}\put(129.5,6){\vector(3,2){7}}
    \put(91,3){$T(C)$}\put(98.5,3){\vector(3,1){7}}
  \end{picture}
\end{center}
\caption{Continuum $C_3$ as an attractor. }\label{3continuum rys}
 \end{figure}

%leżących w przedziałach zwanych oknami stabilności rodziny $\{\al_\lambda\}_{\lambda\in (0,1]}$.
Let  $(\eta_n, \nu_n]$, $n>0$, be the interval of  parameter values for $\lambda$ where $\al_\lambda$  has its first stable orbit of period $2n+1$. The sequences  $\eta_n$, $\nu_n$  converge decreasingly to   $\mu_1$, see Fig. \ref{bifurkacja od mu_1}. Significantly,  passing with  $\lambda$ from  $(\eta_n, \nu_n]$ to $(\eta_{n+1}, \nu_{n+1}]$ the period of the stable periodic orbit of $\al_\lambda$ decreases according to Sharkovskii's order (page \pageref{Sharkovskii's  order}).
\begin{thm}\label{nie mam labela wiec klade szpadela}
Let   $U_\lambda$ be the operator given by \eqref{T_lambda operator defn} and  $\B_\lambda$ the $C^*$-algebra given by \eqref{B_lambda algebra defn}. If  $\lambda \in (\eta_n, \nu_n]$, $n>0$, then
\begin{itemize}
\item[i)] the maximal ideal space $\M_\lambda$ of algebra $\B_\lambda$ compose of a snake-like continuum $M_\infty$ and a sequence of arcs  $\{M_N\}_{N\in \N}$ converging  to $M_\infty$, and
$$
M_\infty=R\cup C_{2n+1},
$$
 where $R$ is a ray,  $C_{2n+1}$ is an irreducible  continuum   with exactly \mbox{$2n+1$} endpoints and whose only proper nondegenerate subcontinua are  arcs. Furthermore,  $C_{2n+1}\cap R=\emptyset$ and $\overline{R}=R\cup C_{2n+1}$.
\item[ii)] Partial homeomorphism  $\tal_\lambda$ generated by $U_\lambda$ on $\M_\lambda$ carries  $M_N$ onto  $M_{N+1}$, for $N\in \N$, and  invariates both $R$ and $C_{2n+1}$. The endpoint of $R$ is a fixed point and the remaining points of  $R$ move toward  $C_{2n+1}$. The endpoints of   $C_{2n+1}$ form a periodic orbit with period  $2n+1$.
\end{itemize}
\end{thm}
\begin{Proof}
By Theorem \ref{sprowadzenie do granicy odwrotnej} it suffices to inductively apply   \cite[Thm. 8]{Barge-Ingram}.
\end{Proof}
For $\lambda \in (\eta_1, \nu_1]$ the continuum $C_{3}\subset \M_\lambda$ is considered to be the simplest example of an irreducible  continuum, cf. \cite{Nadler}.   It may be obtained  as an attractor of a continuous injective map $T:\Omega\to \Omega$ defined on a compact subset  $\Omega=A\cup B \cup C$ of $\R^2$ which acts according to  Fig.  \ref{3continuum rys}.
Actually, the subsystem $ (C_{3}, \tal_\lambda)$ of $(\M_\lambda,\tal_\lambda)$ is  topologically conjugate to the system $(\Lambda, T)$, where $\Lambda=\bigcap_{n\in \N} T^n(\Omega)$.

 \subsubsection{Cascades of bifurcation  that follow the windows of stability}
\begin{figure}[htb]
\begin{center}\setlength{\unitlength}{1mm}
\begin{picture}(240,50)(11,1)

\put(13,-1.5){\includegraphics[angle=0, scale=0.25]{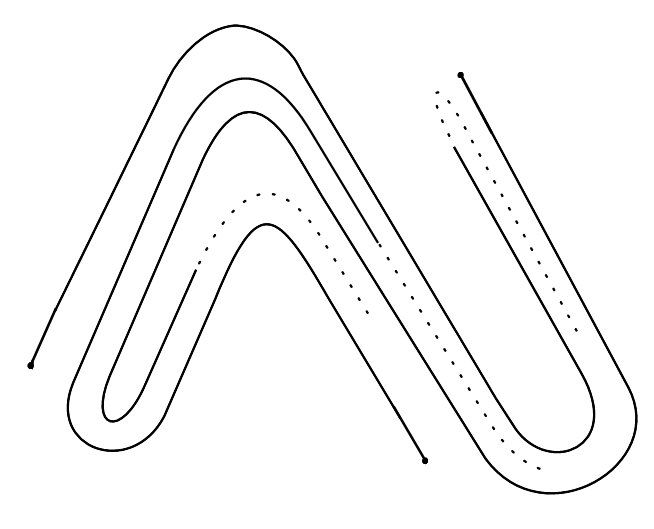}}
  \put(100,-2){\includegraphics[angle=0, scale=0.25]{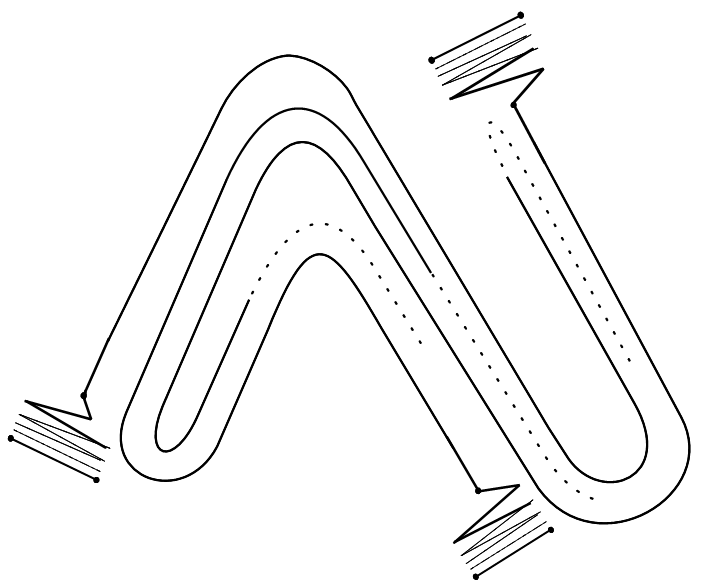}}

 \scriptsize
\put(32.5,48.5){(a)}\put(124,48.5){(b)}

% \put(18.5,32.5){$C_3$}\put(110,32.5){$C_3$}
  \put(12,9.5){$\tomega_{3,1}$}\put(94.8,7.5){$\tomega_{6,1}$}\put(107,3){$\tomega_{6,4}$}

  \put(48,0){  $\tomega_{3,2}$} \put(140,-5.5){$\tomega_{6,2}$}\put(149.3,-1.3){$\tomega_{6,5}$}

 \put(52,40.5){$\tomega_{3,3}$}  \put(133,46){$\tomega_{6,3}$}\put(144.5,50){$\tomega_{6,6}$}

 \put(100,16.5){$R_{0,1}$}\put(133,5){$R_{0,2}$}\put(151,44.5){$R_{0,3}$}
 % \put(104,28){$T(B)$}\put(110,27){\vector(1,-1){4}}
 % \put(122,5){$T(A)$}\put(129.5,6){\vector(3,2){7}}
  %  \put(91,3){$T(C)$}\put(98.5,3){\vector(3,1){7}}
  \end{picture}
\end{center}
\caption{Continua $C_3$ and $C_6$. \label{3 i 6continuum rys}}
 \end{figure}
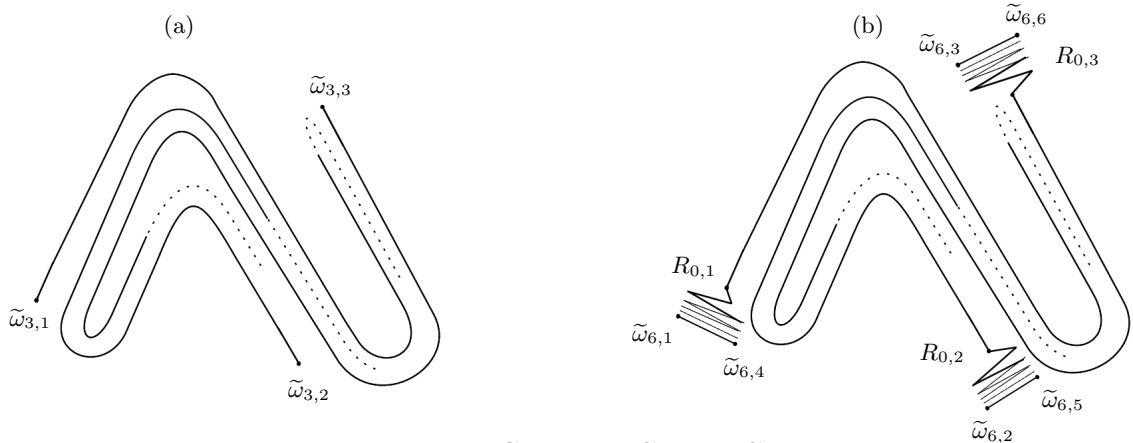

After each of the intervals  $(\eta_n, \nu_n]$, cf. Fig. \ref{bifurkacja od mu_1}, there occurs a cascade of period-doubling  bifurcations. For instance, increasing the parameter value $\lambda$  in the window  $(\eta_1, \nu_1 ]$ of the stable orbit of period three, we observe that  the  continuum $C_3\subset \M_\lambda$ curls around its endpoints $\{\tomega_{3,1},\tomega_{3,2},\tomega_{3,3}\}$, Fig. \ref{3 i 6continuum rys} (a). When we pass $\lambda=\nu_1$ endpoints of $C_3$ grow into arcs:  continuum $C_3$ turns into irreducible continuum $C_6$ with $6$ endpoints, Fig. \ref{3 i 6continuum rys} (b). Afterwards  continuum $C_6$  begin to curl around its endpoints which finally grow into $6$ arcs whose endpoints are  endpoints of an irreducible continuum $C_{12}$. And so on, and so forth. In particular, after four such bifurcations we get a continuum $C_{64}$ which arises from  continuum $C_3$ by replacing each of its endpoint with a copy of the continuum presented on Fig. \ref{czwarta bifurkacja doubling}.
\\
Similar phenomena occur after every  window of stability  $(\eta_n,\nu_n]$, $n>0$.
In order to get a formal description, for each $n>0$, we  denote by   $\lambda^{(n)}_1=\nu_n$, $\lambda^{(n)}_2$, $\lambda^{(n)}_3$, ...   the sequence of parameter values  $\lambda$ that correspond to the cascade of period-doubling bifurcation of the stable orbit appearing immediately after  $\lambda_0^{(n)}:=\eta_n$.
  \begin{thm}
\label{bifurakcyjny opis widma thm4}
Let   $U_\lambda$ be the operator given by \eqref{T_lambda operator defn} and  $\B_\lambda$ the $C^*$-algebra given by \eqref{B_lambda algebra defn}. If  $\lambda \in  (\lambda_m^{(n)},  \lambda_{m+1}^{(n)}]$, for $n>0$,  $m \geq 0$, then
\begin{itemize}
\item[i)] the maximal ideal space $\M_\lambda$ of algebra $\B_\lambda$ compose of the snake-like continuum $M_\infty$ and the sequence of arcs $\{M_N\}_{N\in \N}$ converging  to  $M_\infty$ where
 $$
 M_\infty = R\cup C_{2^m(2n+1)}
 $$
is the union of a ray $R$ and an irreducible  continuum  $C_{2^m(2n+1)}$ with exactly  $2^m(2n+1)$ endpoints. Furthermore $\overline{R}=R\cup C_{2^m(2n+1)}$ and
$$
C_{2^m(2n+1)}= C_{2n+1}\cup \bigcup_{k=0}^{m-1}\,\,\bigcup_{i=1}^{2^k(2n+1)}R_{k,i} \cup \bigcup_{i=1}^{2^{m-1}(2n+1)}I_i
$$
is the union of $(2^m-1)(2n+1)$ rays  $R_{k,i}$  and $2^{m-1}(2n+1)$ arcs $I_i$.
The closure of $C_{2n+1}$ coincides with $C_{2^m(2n+1)}$ and
$$
\overline{R_{k,i}}=\bigcup_{j=0}^{m-k-1}\bigcup_{l=0}^{2^j-1}  R_{k+j,i+ l\cdot 2^k(2n+1)} \cup \bigcup_{l=0}^{2^{m-k-1}}  I_{i+ l\cdot 2^k(2n+1)},
$$
for  $i=1,...,2^{k}(2n+1)$, $k=0,...,m-1$.
\item[ii)]  Partial homeomorphism  $\tal_\lambda$ generated by $U_\lambda$ on $\M_\lambda$ carries  $M_N$ onto  $M_{N+1}$, $N\in \N$, and  $\tal_\lambda: M_\infty\to M_\infty$ is a homeomorphism that preserves $R$ and $C_{2n+1}$, permutes cyclically the arcs $I_i$, and (for each fixed $k=0,...,m-1$) the rays $R_{k,i}$.
\end{itemize}
Moreover, all the rays $R$, $R_{k,i}$, all the arcs $I_i$ and the continuum $C_{2n+1}$ are pairwise disjoint except the following intersections
$$
C_{2n+1}\cap R_{0,i}, \qquad i=1,...,2n+1,
$$
$$
R_{k,i}\cap R_{k,i+2^{k-1}(2n+1)}=\{\tomega_{2^{k-1}(2n+1),i}\}, \qquad i=1,...,2^{k-1}(2n+1), \,\,k= 1,...,m-1.
$$
The sets  $\{\tomega_{2^{k}(2n+1),1}, ... \tomega_{2^{k}(2n+1),2^{k}(2n+1)}\}$  form periodic orbits with periods $2^{k}(2n+1)$ for $k=0,...,m-2$.
The middles of the arcs $I_i$ form a periodic orbit with period $2^{m-1}(2n+1)$ and their endpoints form a periodic orbit with period $2^{m}(2n+1)$.

\end{thm}
\begin{Proof}
In view of Theorem \ref{sprowadzenie do granicy odwrotnej} it suffices to repeat the argument from the proof of \cite[Thm. 8]{Barge-Ingram}, see also remarks preceding  \cite[Thm. 9]{Barge-Ingram}, as well as an introduction to  \cite[6]{Barge-Ingram}.
\end{Proof}

 \section{Reversible extensions of  homeomorphisms of a circle}\label{homeomorphisms of a circle section}
 The characteristic feature of the $C^*$-method  developed in this paper is that  it leads from irreversible dynamics  to reversible dynamics. Therefore, it may seem  surprising that applying it to (already) reversible systems one may also get nontrivial results.  Clearly, all the interesting phenomena, arising in   this case, are related to the freedom of choice of the set $Y$, equivalently the  ideal $J$, see paragraph \ref{podrozdzial o zbiorzez Y}.  As we  show below such considerations  arise naturally in investigation of compressions of unitary operators.
 We start with  the general structure of dynamical systems we will here deal with.
 \begin{prop}\label{opis rozszerzenia ukladu odwracalnego}
If  $\al:M\to M$ is a homeomorphism and $(\M,\tal)$ is a reversible extension of $(M,\al)$ associated with a set  $Y\subset M$,  then $\M$ may be treated as a closed subset
$$
\M=\bigcup_{N\in \N}M_N \cup M_\infty \subset \overline{\N}\times M
$$
 of the product space, where $\overline{\N}=\N\cup \{\infty\}$ is the one point compactification of the discrete space $\N$,
$$
M_\infty=\{\infty\}\times M, \qquad M_N=\{N\}\times \al^N(Y), \,\,\, N \in \N,
$$
 and the partial partial homeomorphism $\tal:\M\to \M$  acts according to the formula
 $$
 \tal(N,x)=(N+1,\al(x)),\qquad \tal(\infty,x)=(\infty,\al(x)).
 $$
 \end{prop}
 \begin{Proof} We define a homeomorphism  $\Psi$  of  $\M=\bigcup_{N\in \N}M_N \cup M_\infty$ onto a closed subspace of  $\overline{\N}\times M$ using  the factor map  $\Phi:\M\to M$, see \eqref{rzut na zerowo wspolrzedno}, by the formulae
 $$
 \Psi(\x):=(\infty, \Phi(\x)),\,\, \textrm{ for }\,\, \x\in M_\infty,\qquad \Psi(\x):=(N, \Phi(\x)),\,\, \textrm{ dla }\,\, \x \in M_N,\,\,\, N \in \N.
$$
Since  $\al$ is a homeomorphism, one readily sees that
$$
\Psi:\bigcup_{N\in \N}M_N \cup M_\infty \to \bigcup_{N\in \N }\{N\}\times \al^N(Y) \cup \{\infty\}\times M
$$
is a   homeomorphism and  identifying  $\M$ with $\Psi(\M)$ the assertion follows.
 \end{Proof}
\subsection{Compression of unitaries generating homeomorphisms of a circle}\label{operatory homeomorfizmy okręgu subsection}
\begin{figure}[htb]
 \begin{center}
\setlength{\unitlength}{1.mm}
\begin{picture}(125,61)(0,0)
\put(-11,-1){\includegraphics[angle=0, scale=0.32]{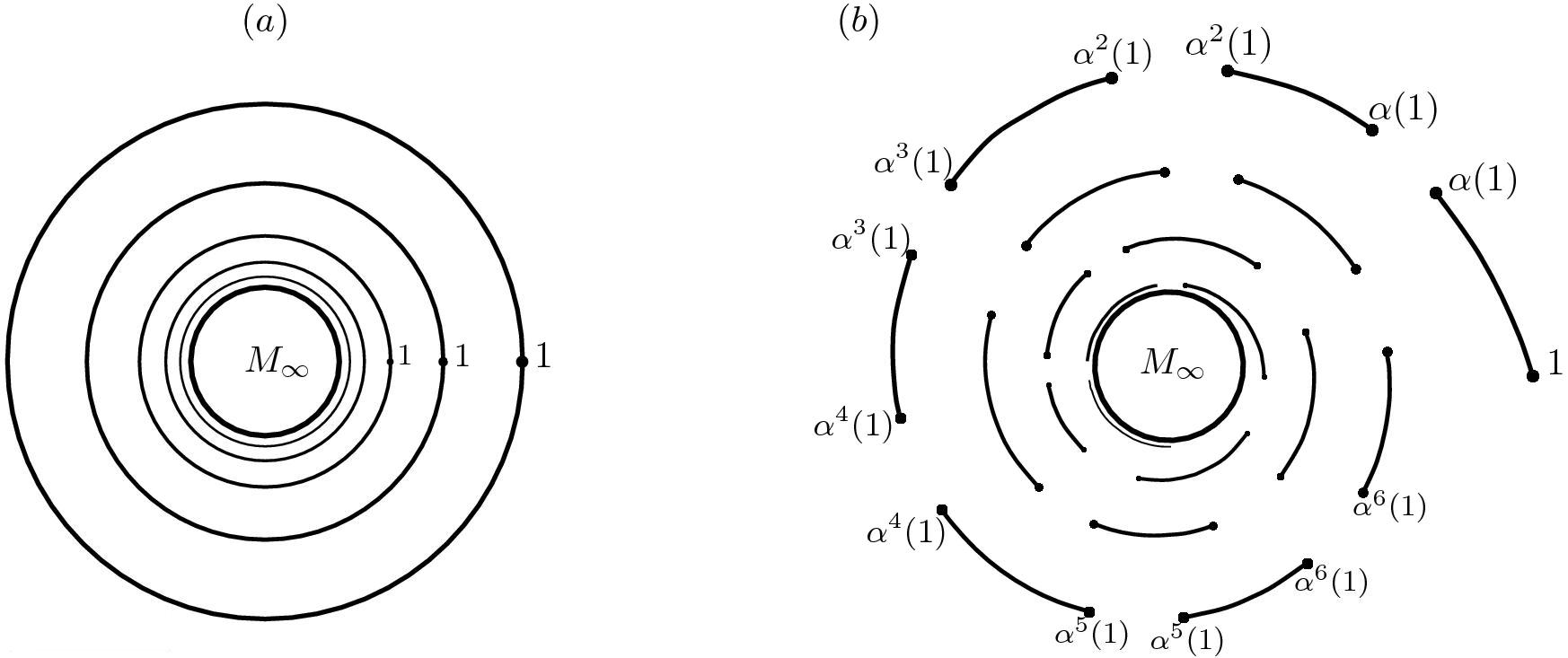}}
 \end{picture}
  \end{center}

  \caption{Spectrum of  $\B$  related to a homeomorphism of the circle.\label{kolka3 rys}}
 \end{figure}

Let  $\al:S^{1}\to S^{1}$  be an orientation preserving homeomorphism of the circle and let $\gamma:\R\to \R$ be its  \emph{lift} to $\R$, i.e. a continuous mapping satisfying
$$
\al(e^{2\pi i t})=e^{2\pi i \gamma(t)}, \qquad t\in [0,1].
$$
We recall that $\gamma$ is an increasing homeomorphism such that  $
\gamma(t+1)=\gamma(t)+1$, $t \in \R$, determined by $\al$ up to a translation by  an integer constant.
We define a unitary operator $\U \in L(\H)$  on the space  $\mathbb{H}=L^2(\R)$ by the formula
$$
(\U f)(t)=\sqrt{|\gamma'(t)|}f(\gamma(t)),
$$
which (by monotonicity of  $\gamma$) make sense for almost all $t$ in $\R$.
In particular,   $
(\U^{*}f)(t)=\sqrt{|(\gamma^{-1})'(t)|}f(\gamma^{-1}(t))
$.
 We let    $\AA \subset L(\H)$ be an algebra of operators of multiplication by continuous periodic functions with period $1$: $
\AA\cong C(S^1)$. Clearly
 \begin{equation}\label{relacje przed kompresją}
 \U\, \AA\, \U^*\subset \AA, \qquad \U^*\,   \AA\, \U\subset \AA,
 \end{equation}
 and we have
\begin{prop}
The operators $\U$ and  $\U^*$ generate  on the maximal ideal space of $\AA$ (identified with $S^1$) the systems  $(S^{1},\al)$ and   $(S^{1},\al^{-1})$, respectively.
\end{prop}
Let us now consider compressions  of the introduced  objects to the space $H:=L^2([0,\infty))$ naturally treated as a subspace of $\H=L^2(\R)$. %(We may likewise treat $L(H)$ as a   full corner subalgebra of $L(\H)$). 
Namely, we denote by  $P:\H \to H$ the projection from  $\H$ onto  $H$ and we put
$$
U:=P\, \U P, \qquad \A:=P\, \AA\, P.
$$
Then the algebra $\A \subset L(H)$ is isomorphic  to $C(S^1)$ and   $U\in L(H)$ is a partial isometry such that $U^*=P \U^*  P\in L(H)$.
\begin{prop}\label{stwierdzenie o kompresjach}
Within the above notation the following possibilities may occur.
\begin{itemize}
\item[i)] If  $\gamma(0)> 0$, then  $U$ is a non-invertible coisometry in $L(H)$,
$$
U\, \A \, U^* \subset \A,\qquad  U^* \,\A\, U \nsubseteq \A
,$$
operator $U$ generates on the spectrum of $\A$ the system $(S^1,\al)$ and \\
$
\hull(U^*U\A\cap \A)=\left\{e^{2\pi i t}: t \in [0,\gamma(0)]\right\},
$
\item[ii)] If  $\gamma(0)= 0$, then $U$ is unitary.
$$
U\, \A\, U^{*} \subset \A,\qquad
U^* \, \A \, U\subset \A,
$$
 $U$ and $U^*$ generate respectively the systems $(S^{1},\al)$ and   $(S^{1},\al^{-1})$.
\item[iii)] If $\gamma(0) < 0$,  then $U$ is non-invertible isometry:
$$
U^* \, \A \, U \subset \A,\qquad  U \,\A\,  U^* \nsubseteq \A,
$$
 operator $U^*$ generates on the spectrum of $\A$ the system $(S^1,\al^{-1})$ and \\
$
\hull(U^*U\A\cap \A)=\left\{e^{2\pi i t}: t \in [0,\gamma^{-1}(0)]\right\}.
$
\end{itemize}
\end{prop}
\begin{Proof}
If  $\gamma(0)>0$, then
$$
(Uf)(t)=\sqrt{|\gamma'(t)|}f(\gamma(t)),\quad (U^*f)(t)=\begin{cases}
\sqrt{|(\gamma^{-1})'(t)|}f(\gamma^{-1}(t)),& t \in [\gamma(0),\infty) \\
0, & t \in [0,\gamma(0)).
\end{cases}
$$
In particular, the projection   $U^*U$ is the operator of multiplication by the characteristic function  of $[\gamma(0),\infty)$ and thereby $
\hull(U^*U\A\cap \A)=\left\{e^{2\pi i t}: t \in [0,\gamma(0)]\right\}
$.
Item  iii) is  obtained   by reversing the roles of  $\gamma$ and $\gamma^{-1}$ in item i), and item  ii) is straightforward.
\end{Proof}
We see that in a process of compression one may lose one of  the relations \eqref{relacje przed kompresją}. However, according to our results  one may always retrieve what is lost by passing to a bigger algebra. To fix attention  let us from now on assume that
$$
\gamma(0) >0
$$
(the case when $\gamma(0)<0$ is completely analogous). We put
$$
\B=C^*\left( \bigcup_{n=0}^\infty  U^{*n} \A U^n \right),
$$
which by Theorem \ref{takie sobie stw2},  is the smallest $C^*$-algebra containing  $\A$  such that
$$
U\, \B\, U^* \subset \B, \qquad U^*\, \B \, U \subset \B.
$$
%Furthermore, we have
\begin{thm}\label{opis ukladu dualnego do homeomorfizmow okregu}
 The spectrum of  $\B$ assumes one of the forms:
\begin{itemize}
\item[i)] If $\gamma(0)\geq 1$, then  $
\M=\overline{\N}\times S^1$, see Fig. \ref{kolka3 rys} (a).
\item[ii)] If $\gamma(0)\in (0, 1)$, then $\M=\bigcup_{N\in \N}M_N \cup M_\infty \subset \overline{\N}\times S^1$, where $M_\infty=\{\infty\}\times S^1$ is a circle, and  the sets $M_N$, $N\in \N$, are arcs:
$$
M_N=\{N\}\times [\al^N(1),\al^{N+1}(1)],
$$
where $[\al^N(1),\al^{N+1}(1)]$ stands for an arc on $S^1$ with the origin $\al^N(1)$ and ending $\al^{N+1}(1)$, see Fig. \ref{kolka3 rys} (b).
\end{itemize}
Under the identification  $\B=C(\M)$ the operator  $U$ generates on $\M$ the  mapping  $\tal$ that acts according to the formula
 $$
 \tal(N,x)=(N+1,\al(x)),\qquad \tal(\infty,x)=(\infty,\al(x)).
 $$
\end{thm}
\begin{Proof}
By  Proposition \ref{stwierdzenie o kompresjach} and Theorem \ref{opis kosmiczakow} the system   $(\M,\tal)$ is the reversible extension of $(S^1,\al)$ associated to the set $
Y=\left\{e^{2\pi i t}: t \in [0,\gamma(0)]\right\}
$
which either is a circle $S^1$, when  $\gamma(0)\geq 1$, or an arc    with the origin  $1=e^{2\pi i 0}$ and ending $\al(1)=e^{2\pi i \gamma(0)}$. Thus it suffices  to apply Proposition \ref{opis rozszerzenia ukladu odwracalnego}.
\end{Proof}

It follows that the  algebra $\B$  depends not only on  the homeomorphism  $\al$ but also on the choice of the lift  $\gamma$. In particular, if  $\gamma(0) \geq 1$, then independently of $\al$
$$
\B\cong C(\overline{\N}\times S^1).
$$
If however $\gamma(0)\in (0,1)$,  the structure  of $\B$ is uniquely determined by the orbit of the point $1\in S^1$. We will discuss this issue in detail below.

 \subsection{Classification of the extended algebras  via  rotation numbers}
 \label{kto nie lubi liczb obrotu ten niech nie czyta}
We recall \cite{Brin}, \cite{Devaney}, that if  $\gamma:\R\to \R$ is a lift  of a homeomorphism  $\al:S^1\to S^1$,  the limit $\lim_{n\to \infty}\frac{\gamma(t)}{n}$  always exits, does not depend on  $t\in \R$ and its fractional part
\begin{equation}\label{rotation number euqation}
\tau(\al):= \mantysap \, \lim_{n\to \infty}\frac{\gamma(t)}{n} \mantysak \in [0,1)
\end{equation}
  depend only on $\al$. The quantity  $\tau(\al)$  \label{rotation map defn}  is called the  \emph{rotation number} of  $\al$. Its  role is explained by the following  theorem due to H. Poincare, cf. \cite{Brin}, \cite{Devaney}. The symbol $\Theta_{\tau}$, $\tau\in [0,1)$, stands for the rotation by $2\pi\tau$: $\Theta_{\tau}(z)=z \cdot e^{2\pi i\tau}$, $z\in S^1$.
\begin{thm}[Poincare classification]\label{twierdzenie o Poincare machen}\index{twierdzenie!Poincar\'e}
Let $\al:S^1\to S^1$ be an orientation preserving homeomorphism of the circle.
\begin{itemize}
\item[1)] If $\tau(\al)=\frac{m}{n}$, where  $m$ and $n$ are coprime, then all periodic points for $\al$ are of period $n$ (and there exists at least one such a point).
\item[2)] If $\tau(\al)\notin \Q$, then $\al$ does not possess periodic points and the set  $\Omega(\al)$ consisting of the accumulation points of an arbitrary orbit   $\{\al(x)\}_{n\in \Z}$ does not depend on a choice of  $x\in S^1$.  There are two possible subcases
\begin{itemize}
\item[a)] When $\Omega(\al)=S^1$, that is when $\al$ is topologically transitive, then the system $(S^1,\al)$ is topologically conjugated to $(S^1,\Theta_{\tau(\al)})$.
\item[b)] When  $\al$ is not topologically transitive, then $\Omega(\al)$ is a perfect nowhere dense subset of $S^1$ and there exits a continuous surjection  $\phi:\Omega(\al)\to S^1$ 
 which is a semiconjugacy from   $(\Omega(\al),\al)$ onto  $(S^1,\Theta_{\tau(\al)})$.
\end{itemize}
\end{itemize}
\end{thm}
We  use the above theorem to classify the algebras  $B$ described in item ii)  of Theorem  \ref{opis ukladu dualnego do homeomorfizmow okregu}, that is we  assume throughout this subsection that  the lift  $\gamma:\R\to\R$ satisfies
\begin{equation}\label{warunek na gamma dla klasyfikacji}
 0<\gamma(0) <1.
 \end{equation}
Such lift always exists, provided $1\in S^1$ is not a fixed point of $\al:S^1\to S^1$. Since $\gamma$  is uniquely determined by $\al$, so is the operator $U$ and the algebra $\B\cong C(\M)$,  defined in the previous subsection.  Thus,  it makes sense to adopt the following notation
 $$
 \B_\al:=\B\,\,\, \textrm{  and  }\,\,\, \M_{\al}:=\M.
 $$
The space  $M_\al$ compose of a circle $M_\infty$ and a sequence of arcs $\{M_N\}_{N\in \N}$.
   \begin{thm}\label{twierdzenie o klasyfikacji poincare algebr operatorowych}
In  the  situation under  consideration the following cases may occur:
\begin{itemize}
\item[1)] If  $\tau(\al)=\frac{m}{n}$ where  $m$ and $n$ are coprime, then the limit points of the endpoints of arcs $\{M_N\}_{N\in\N}$ form a subset  of $M_\infty$ with cardinality  $n$,  and
$$
\B_\al \cong \B_{\Theta_{\frac{m}{n}}}.
$$
\item[2)] If  $\tau(\al)\notin \Q$,  then the two subcases are possible:
\begin{itemize}
\item[a)] $\al$ is topologically transitive,  and then $M_\infty$ is the set of limit points of the endpoints of arcs $\{M_N\}_{N\in\N}$ and
$$
\B_\al \cong \B_{\Theta_{\tau(\al)}}.
$$
\item[b)]   $\al$ is not topologically transitive, and then  the set of limit points of the endpoints of arcs  $\{M_N\}_{N\in\N}$ form a perfect nowhere dense subset of  $M_\infty$. In particular,
$$
\B_\al \neq \B_{\Theta_{\tau}},\qquad \tau \in [0,1).
$$
\end{itemize}
\end{itemize}

 \end{thm}
 \begin{Proof} The set of limit points of the endpoints of arcs $\{M_N\}_{N\in\N}$  coincides with the set  $\{\infty\}\times\Omega(\al)\subset M_\infty$ where $\Omega(\al)$ is the set of limit points of the orbit $\{\al^N(1)\}_{N\in \N}$. By Theorem \ref{twierdzenie o Poincare machen} we only need to consider  the cases listed in the assertion.
 \par
 1) The set $\Omega(\al)$ consists of  $n$ points that form a periodic orbit of $\al$. More precisely, there are  $n$ points $x_0$, $x_1$, ..., $x_{n-1}\in S^1$ enumerated according to the orientation and such that
 $$
\lim_{N\to \infty}\al^{Nn+k}(1)=x_{km}, \quad \qquad k=0,...,m-1.,
 $$
cf.  \cite{Brin}. Thus it follows that
 $$
\lim_{N\to \infty}M_{Nn+k}=\{\infty\}\times [x_{km\,\, (mod\,\,n)}, x_{(k+1)m\,\, (mod\,\,n)}], \qquad k=0,...,n-1,
 $$
that is the sequence of arcs  $\{M_{Nm+k}\}_{N\in\N}$ converge in Hausdorff metric to the arc on  $M_\infty$ with the origin $(\infty, x_{kn\,\, (mod\,\,m)})$ and ending  $(\infty,x_{(k+1)n\,\, (mod\,\,m)})$.
Let  $\phi:\M_\al \to \M_{\Theta_{\frac{m}{n}}}$ be the mapping that   acts  "linearly" according to the scheme:
$$
\{\infty\}\times [x_k,x_{k+1 (mod n)}]\stackrel{\phi}{\longmapsto} \{\infty\}\times \left[e^{2\pi i \frac{k}{n}},e^{2\pi i \frac{k+1}{n}}\right], \qquad k=0,...,n-1,
$$
$$
\{N\}\times [\al^{N}(1),\al^{N+1}(1)]\stackrel{\phi}{\longmapsto} \{N\}\times \left[e^{2\pi i \frac{N}{n}},e^{2\pi i \frac{N+1}{n}}\right],\qquad N\in \N.
$$
It is evident  that  $\phi$ is a homeomorphism and hence the algebras $\B_\al=C(\M_\al)$ and $\B_{\Theta_{\frac{m}{n}}}=C(\M_{\Theta_{\frac{m}{n}}})$ are isomorphic.
\par
2a) We have $\{\infty\}\times \Omega(\al)=\{\infty\}\times S^1=M_\infty$  and there exists a  homeomorphism $\phi:S^{1}\to S^1$ such that
$$
\begin{CD}
S^1 @>{\al}>>  S^1  \\
@V{\phi}VV               @VV{\phi}V\\
S^1 @>{\Theta_{\tau(\al)}}>>  S^1
\end{CD}
$$
is commutative.  Furthermore, $\phi$ may be arranged so  that $\phi(1)=1$, see the proof of \cite[Thm. 7.1.9]{Brin}. It follows that the mapping
$
id \times \phi:\M_\al\longrightarrow  \M_{\Theta_{\tau(\al)}}$:
$$
(id \times \phi)(N,x)=(N,\phi(x)),\qquad  (N,x)\in M_N,\,\, N\in \overline{\N},
$$
is a homeomorphism. Hence $\B_\al\cong \B_{\Theta_{\tau(\al)}}$.
 \par
 2b) Since   $\{\infty\}\times \Omega(\al)$ is  a perfect nowhere dense subset of   $M_\infty$ the space $\M_\al$ is not homeomorphic to any of the spaces  $\M_{\Theta_\tau}$, $\tau \in [0,1)$. Equivalently  $
\B_\al \neq \B_{\Theta_{\tau}}$,  for all $\tau \in [0,1)$.
 \end{Proof}
If $\phi$ is a topological conjugacy between $(S^1,\al)$ and $(S^1,\beta)$, then either $\tau(\al)=\tau(\beta)$  (when $\phi$ is orientation preserving) or  $\tau(\al)+\tau(\beta)=1$ (when  $\phi$ changes the orientation), so the rotation number is "almost an invariant" for homeomorphisms of the circle.   For the algebras  $\B_\al$ the rotation number is an invariant  \emph{sensu stricto}.
\begin{thm}\label{liczba obrotu jako niezmiennik algebr}
If algebras $\B_\al$ and $\B_{\beta}$ are isomorphic, then $\tau(\al)=\tau(\beta)$.
\end{thm}
\begin{Proof}
Suppose that $\B_\al$ and $\B_{\beta}$ are isomorphic. There exists a homeomorphism $\phi:\M_\al \to \M_{\beta}$, where $\M_\al=\bigcup_{N\in\overline{\N}} M_{N}$ and $\M_{\beta}=\bigcup_{N\in\overline{\N}} M_{N}'$ are maximal ideal spaces of  $\B_\al$ and  $\B_{\beta}$ respectively. Clearly,    $\phi$ necessarily  carries the arcs $\{M_N\}_{N\in\N}$ onto arcs $\{M_N'\}_{N\in\N}$,  the circle $M_\infty$ onto circle $M_\infty'$, and the set  $\Omega(\al)\subset M_\infty$ of limit points of endpoints of $\{M_N\}_{N\in\N}$ onto the set  $\Omega(\beta)\subset M_\infty'$ of limit points of endpoints of the arcs $\{M_N'\}_{N\in\N}$.
For the simplicity of notation we  adopt the identification $M_\infty=M_\infty'=S^1$. We claim that $\phi$ establishes conjugacy between the systems $(\Omega(\al),\al)$ and $(\Omega(\beta),\beta)$ or $(\Omega(\al),\al)$ and $(\Omega(\beta),\beta^{-1})$ depending on whether  $\phi:M_\infty\to M_\infty'$ preserves or changes the orientation.
Once we  prove this, the standard argument   give us  that either   $\tau(\al)=\tau(\beta)$ or  $\tau(\al)+\tau(\beta^{-1})=1$, where  in the letter case   we get  $\tau(\al)=\tau(\beta)$ since $\tau(\beta^{-1})=1 -\tau(\beta)$.
\\
To prove our claim we fix a sequence $\{\al^{N_k}(1)\}_{k\in \N}$  converging to an arbitrarily chosen point $x_0\in \Omega(\al)$.  Then the sequence $\{\al^{N_k+1}(1)\}_{k\in \N}$ converges to  $\al(x_0)$, and the sequence of arcs $\{M_{N_k}\}_{N\in \N}$ converges  (in Hausdorff metric)  to the arc $[x_0, \al(x_0)]$.
In the case   $\phi:M_\infty \to M_\infty'$ preserves the orientation, almost all arcs from the sequence $\{M_{N_k}\}_{N\in \N}$  are mapped in accordance with (the natural) orientation onto almost all arcs of the sequence $\{\phi(M_{N_k})\}_{N\in \N}$. Hence
$$
[\phi(x_0),\phi(\al(x_0))] =\phi([x_0, \al(x_0)])=\phi(\lim_{k\to \infty}M_{N_k})
=[\phi(x_0),\beta(\phi(x_0))].
$$
Thus  $\phi(\al(x_0))=\beta(\phi(x_0))$ and consequently  $\phi$ conjugates the systems $(\Omega(\al),\al)$,  $(\Omega(\beta),\beta)$.
In the case $\phi:M_\infty \to M_\infty'$ changes the orientation, arguing similarly as above one gets
$$
[\phi(\al(x_0)),\phi(x_0)]=\phi([x_0, \al(x_0)])
=\lim_{k\to \infty}\phi(M_{N_k})=[\beta^{-1}(\phi(x_0)),\phi(x_0)].
$$
Hence $\phi(\al(x_0))=\beta^{-1}(\phi(x_0))$ and consequently  $\phi$ conjugates the systems $(\Omega(\al),\al)$ and   $(\Omega(\beta),\beta^{-1})$. This proves our claim.
 \end{Proof}
Applying the classical result of A. Denjoy \cite{Denjoy},  \cite{Brin}, \cite{Devaney} which states that every diffeomorphism of the circle  with finite variation and irrational rotation number is topologically transitive  we get that in the class of  algebras  $\B_\al$ associated with such diffeomorphisms  the rotation number is not only an invariant but actually a numerical equivalent.
 \begin{thm}
If one of the homeomorphisms  $\al$, $\beta$ is a diffeomorphism with finite variation, then the algebras $\B_\al$ and $\B_{\beta}$ are isomorphic if and only if  $\tau(\al)=\tau(\beta)$.
 \end{thm}
 \begin{Proof}
Apply Theorems \ref{twierdzenie o klasyfikacji poincare algebr operatorowych}, \ref{liczba obrotu jako niezmiennik algebr} and Denjoy Theorem.
 \end{Proof}

\end{document}